\documentclass[a4paper,12pt,onesided]{article}
 \usepackage[all]{xy}
 \usepackage{amsfonts}
 \usepackage{amsmath}
  \usepackage{amsthm}
  \usepackage{mathpazo}

 \newtheorem{mth}{Theorem}[section]
 \newtheorem{mlem}[mth]{Lemma}
 \newtheorem{mprop}[mth]{Proposition}
 \newtheorem{mcor}[mth]{Corollary}

 \theoremstyle{definition} 
 \newtheorem{mdef}[mth]{Definition}

 \newcommand{\mK}{\mathbb{K}}
 \newcommand{\mR}{\mathbb{R}}
 \newcommand{\mC}{\mathbb{C}}
 \newcommand{\mZ}{\mathbb{Z}}
 \newcommand{\mN}{\mathbb{N}}
 \newcommand{\Mod}{\mathcal{M}}
 \newcommand{\mul}{m}
 \newcommand{\comul}{\Delta}
 \newcommand{\counit}{\varepsilon}
 \newcommand{\unit}{\eta}
 \newcommand{\act}{\rho}
 \newcommand{\coact}{\rho}
 \newcommand{\sw}[1]{{}_{(#1)}}
 \newcommand{\co}{\text{co}}
 \newcommand{\ot}{\otimes}
 \newcommand{\can}{\text{\rm can}}
 \newcommand{\trans}{\tau}
 \newcommand{\twa}{{}^{[1]}}
 \newcommand{\twb}{{}^{[2]}}

\newcommand{\eqmap}[3]{\ar@<1ex>[#1]^-{#2}  \ar@<-1ex>[#1]_-{#3}}
\newcommand{\twosid}[3]{\ar@<1ex>@{<-}[#1]^{#2}  \ar@<-1ex>[#1]_{#3}}

\newcommand{\cp}{{}^{\text{\bf c}}}
\newcommand{\mIm}{\text{Im}}
\newcommand{\ffun}{\vartheta}
\newcommand{\spec}{\sigma}
\newcommand{\Bnd}{\mathbf{B}}
\newcommand{\Hil}[1]{\mathcal{#1}}
\newcommand{\pchoose}[3]{
\left[
\begin{array}{c}
#1\\
#2
\end{array}
\right]_{#3}
}
\newcommand{\cm}{\kappa}
\newcommand{\uwa}[1]{{}^{[1]_{#1}}}
\newcommand{\uwb}[1]{{}^{[2]_{#1}}}
\newcommand{\transc}{{\tau\cp}}

\newcommand{\mdeg}{\text{deg}}
\newcommand{\Span}{\text{Span}}

\title{Locally coalgebra-Galois extensions}
\author{Bartosz Zieli\'nski\footnote{
Department of Theoretical Physics II, University of \L{}\'od\'z,
Pomorska 149/153 90-236 \L{}\'od\'z, Poland;
bzielinski@uni.lodz.pl.
}}

\begin{document}

\maketitle

\begin{abstract}
\noindent The paper introduces the notion of a {\em locally coalgebra-Galois extension} and,
as its special case, a {\em locally cleft extension}, generalising concepts from \cite{CalMat:Con}.
The necessary and sufficient conditions for a locally coalgebra-Galois extension to be a (global)
coalgebra-Galois extension are stated. As an important special case, it is proven, 
that under not very restrictive conditions
the gluing of two locally cleft extensions is a globally coalgebra-Galois extension.
As an example, the quantum lens space of positive charge is constructed by gluing of two quantum solid tori.
\end{abstract}

\section{Introduction}
Constructing new topological spaces by gluing together several known
 ones, or studying the properties
of a given space by presenting it as a patching of topological spaces of a simpler structure,
is the standard method in the classical topology which was frequently adapted to the
noncommutative geometry.
Examples include the quantum real projective sphere (\cite{HaMaSz:Pro}),
and the Podle\'s sphere (defined in \cite{Pod:Sph}, 
it was proven in \cite{Sheu:Pod} and \cite{HaMaSz:Pro}
that it is $C^\ast$-isomorphic with gluing of two quantum discs),
the Matsumoto sphere and the quantum lens space (\cite{Mat:Sph} and \cite{MatTom:Lens}).

The basic idea about covering quantum spaces by quantum subsets stems from the observation that
an ideal of the algebra of functions on a given quantum space can be interpreted as
consisting of those functions, which assume the value zero on some quantum subset.
Then the quotient algebra can be viewed as the algebra of functions on this quantum subset.

Suppose that the algebra of functions on some quantum space has a family of ideals,
which intersect at zero. Dually it means that the corresponding quantum subsets cover 
the whole of quantum space.

The covering and gluing of $C^\ast$-algebras and the notion of a 
locally trivial quantum principal bundles
in the context of $C^\ast$-algebras were introduced in \cite{BuKo:Pri}. The purely algebraic 
theory of covering and gluing  of algebras and differential algebras was presented in \cite{CalMat:Glu},
which was followed by the purely algebraic definition of locally trivial quantum principal
bundles in \cite{CalMat:Con}. An example of a locally trivial principal bundle was produced 
in \cite{CalMat:Con} and further elaborated in \cite{HaMaSz:LocHopFib}.
The theory resembling locally cleft extensions, based on sheaf theory was independently
developed in \cite{Pfl:Fib}.

In what follows we introduce the concept of {\em locally coalgebra-Galois extensions} which, 
without going into technical details, are the algebras and comodules which have covers such that each
of the quotient spaces are coalgebra-Galois extensions. Of a particular interest are  the
conditions which ensure that a locally coalgebra-Galois extension is  a (global) coalgebra-Galois
extension. Later we concentrate on the special case, in which all the quotient spaces of the cover
of a locally Galois extension are cleft, which can be considered as the natural generalisation
of locally trivial quantum principal bundles introduced in \cite{CalMat:Con}.
 
\section{Basic definitions and notation}

We work over the general commutative ring $\mK$.

\paragraph{Categories of modules and comodules.}
Let $C$ and $H$ be coalgebras, and $A$ and $B$ algebras. We denote by ${}_A\Mod$,
$\Mod_B$, ${}^H\Mod$, $\Mod^C$, ${}_A\Mod_B$, ${}^H\Mod^C$, ${}_A\Mod^C$, etc, respectively, 
the category of left $A$-modules, right $B$-modules, left $H$-comodules, right $C$-comodules,
$(A,B)$-bimodules, $(H,C)$-bicomodules,
left $A$-modules and right $C$-comodules such that
the $C$-coaction  commutes with the $A$-action, etc.
Unadorned $\Mod$ will denote category of $\mK$-modules. All algebraic objects considered belong
to this category. 

\paragraph{Comultiplication, coaction and the Sweedler notation.} 
Suppose $C$ and $H$ are coalgebras, $M\in {}^H\Mod$ and $N\in\Mod^C$.
We denote the comultiplication by $\comul:C\rightarrow C\ot C$, the left $H$-coaction by
${}^H\coact:M\rightarrow H\ot M$, the right $C$-coaction by 
$\coact^C:N\rightarrow N\ot C$. Occasionally, to avoid confusion, we indicate the module being 
coacted on, writing $\coact^C_N$ for $\coact^C:N\rightarrow N\ot C$, and similarly for 
left coactions. 
We also use the Sweedler notation: $\comul(c)=c\sw{1}\ot c\sw{2}$, ${}^H\coact(m)=m\sw{-1}\ot m\sw{0}$,
$\coact^C(n)=n\sw{0}\ot n\sw{1}$ for all $c\in C$, $m\in M$, $n\in N$, and the summation is implicitly understood.

\paragraph{Antipodes, counits and units.} Let $C$ be a coalgebra. We denote by 
$\counit^C:C\rightarrow \mK$ the counit of $C$. 
If there is no danger of confusion, we write $\counit$ for $\counit^C$.
If $C$ is a Hopf algebra, then the antipode of $C$ is
denoted by $S:C\rightarrow C$. If $A$ is an algebra, in most cases we
use the symbol $1_A$ or simply $1$ for the unit of $A$. 
Occasionally, we may need to write the unit explicitly as a map
$\unit:\mK\rightarrow A$, $k\mapsto k1_A$. Whenever it does not cause any ambiguity, we
identify in notation the ground ring $\mK$ with the 
subalgebra $\mK 1_A$ of $A$, i.e., depending on context,
for any $k\in\mK$, $k$ may mean also $k1_A$ for some algebra $A$.

\paragraph{Entwining structures} Let $C$ be a coalgebra and $A$ an algebra. 
Entwining structures were introduced in \cite{BrzMaj:CoalBun}
as a very general 
way of linking algebra
structure on $A$ and coalgebra structure on $C$. 
A good introduction to entwining strtuctures can be found in \cite{BrzWis:Cor}.
Throughout this paper by entwining map we mean a right-right version
$\psi:C\ot A\rightarrow A\ot C$, the corresponding entwining structure is denoted
$(A,C)_\psi$.

We use the following summation notation for an entwining map $\psi$:
\[
\psi(c \ot a)=a_\alpha\ot c^\alpha, \ 
\text{for all }c\in C,\ a\in A, 
\]
where  
small Greek  
letters are used for implicit summation indices.

With this summation notation, the bow-tie diagram condition (cf. \cite{BrzWis:Cor}) can be explicitly written as
\begin{subequations}
\label{bowtie}
\begin{gather}
1_\alpha\ot c^\alpha=1\ot c,\label{bowtiea}\\
a_\alpha\counit^C(c^\alpha)=a\counit^C(c),\label{bowtieb}\\
(aa')_\alpha\ot c^\alpha=a_\alpha a'{}_\beta\ot c^{\alpha\beta},\label{bowtiec}\\
a_\alpha\ot c^\alpha\sw{1}\ot c^\alpha\sw{2}=a_{\beta\alpha}\ot c\sw{1}{}^{\alpha}\ot c\sw{2}{}^\beta,\label{bowtied}
\end{gather}
for all $a,a'\in A$, $c\in C$.
\end{subequations}

Let $(A,C)_\psi$ be an entwining structure, and let $P$ be an algebra and an $(A,C)_\psi$-entwined module.
An algebra extension $B\subseteq P$ is called an {\em $(A,C)_\psi$-extension} if and only if
$B=P^{\co C}$. Of particular interests are $(A,C)_\psi$-extensions $B\subseteq A$. Such extensions
are denoted by $A(B,C,\psi)$. In this case, if there exists a grouplike element $e\in C$
such that,  for all $a\in A$, $\coact^C(a)=\psi(e\ot a)$ then $A_e(B,C,\psi)$ is called an $e$-copointed
$(A,C)_\psi$-extension.

\paragraph{The canonical map and quantum principal bundles}
Suppose that $C$ is a coalgebra and $P$ is an algebra and a right $C$-comodule.
Let $B=P^{\co C}$ be a subalgebra of coinvariants of right $C$-coaction.
To set the notation we recall the definition of the canonical map,
\begin{equation*}
\can^C_P:P\ot_B P\rightarrow P\ot C,\ \ p\ot_Bp'\mapsto pp'\sw{0}\ot p'\sw{1}.
\end{equation*}
If the canonical map is a bijection then $P(B)^C$ is called a $C$-coalgebra Galois extension
of $B$. If, in addition $C$ is a Hopf algebra and $P$ is a $C$-comodule algebra,
then $P(B)^C$ is called a $C$-Hopf Galois extension.
We recall the definition of the translation map,
\begin{equation*}
\trans^C_P:C\rightarrow P\ot_BP,\ \ c\mapsto(\can^C_P)^{-1}(1\ot c),
\end{equation*}
for all $c\in C$.
We use an explicit 
`Sweedler like' notation for the translation map,
\begin{equation*}
\trans^C_P(c)=c\twa\ot_B c\twb, \text{ for all } c\in C,
\end{equation*}
where an implicit summation is understood. The translation map has a number of 
useful properties (cf.  34.4 \cite{BrzWis:Cor}). In particular,
for all $c\in C$, $p\in P$,
\begin{gather}
 c\twa c\twb\sw{0}\ot c\twb\sw{1}=1_P\ot c, \label{transdef}\\
 p\sw{0}p\sw{1}\twa\ot_B p\sw{1}\twb=1_P\ot_B p,\label{trpp}
 \end{gather}
Entwinings and $C$-coalgebra Galois extensions are closely related. If $P(B)^C$ is a $C$-coalgebra
Galois extension, then the map (cf. Theorem 2.7 \cite{BrzHaj:Coext})
\begin{equation}
\psi_\can:C\ot P\rightarrow P\ot C,\ \ 
c\ot a\mapsto\can^C_P(\trans^C_P(c)p)\label{canent}
\end{equation}
is the unique entwining such that $P$ is an entwined module. 
In particular, if $P(B)^C$ is a $C$-Hopf
 Galois extension, then $\psi_\can(c\ot p)=p\sw{0}\ot cp\sw{1}$.

\paragraph{Multiplication of subsets} Suppose that $A$ is an algebra, $M\in{}_A\Mod$ and 
$\act_A:A\otimes M\rightarrow M$ denotes the left $A$-action. Let
$B\subseteq A$, $N\subseteq M$ be subsets. Unless otherwise stated, in what follows,
we denote $BN\equiv\coact_A(B\ot N)$. We use similar convention for right modules.

\subsection{Cleft extensions}

An interesting class of coalgebra-Galois extensions is provided by cleft extensions.

\begin{mdef} Let $C$ be a coalgebra, $P$ be an algebra and a right $C$-comodule. 
A convolution invertible and  a right $C$-colinear map 
 $\gamma:C\rightarrow P$  is called a {\em cleaving map}.
 A $C$-coalgebra Galois extension $P(B)^C$ such that
  there exists a cleaving map $\gamma:C\rightarrow P$
 is called a {\em cleft coalgebra Galois extension}
  and is denoted by $P(B)^C_\gamma$. Similarly a {\em cleft $(P,C)_\psi$-extension} 
 $P_\gamma(B,C,\psi)$ is a $(P,C)_\psi$-extension with a cleaving map $\gamma$.  
 \end{mdef}
 Observe that if $P(B)^C$ is cleft, then the inverse of the canonical map has the form
 \begin{equation}
 (\can^C_P)^{-1}(p\ot c)=p\gamma^{-1}(c\sw{1})\ot_B\gamma(c\sw{2}),
 \text{ for all } c\in C, p\in P,\label{cleftcaninv}
 \end{equation}
 where $\gamma^{-1}$ means the convolution inverse.
   
 \begin{mprop}(Proposition 2.3 \cite{Brz:ModAs}.)\label{brzcl}
 Let $C$ be a coalgebra, $P$ be a right comodule and $B=P^{\co C}$. If there exists a cleaving
 map $\gamma:C\rightarrow P$, then the following are equivalent:
 \begin{enumerate}
 \item $P(B)^C$ is a $C$-coalgebra Galois extension.
 \item There exists an entwining $\psi$ such that $P(B,C,\psi)$ is a $(P,C)_\psi$-extension.
 \item For all $p\in P$, $p\sw{0}\gamma^{-1}(p\sw{1})\in B$.
 \end{enumerate}
 If any of the above conditions hold, then $P\simeq B\ot C$ in ${}_B\Mod^C$, where the isomorphism 
 $\theta_\gamma:P\rightarrow B\ot C$ and its inverse $\theta_\gamma^{-1}:B\ot C\rightarrow P$
 are given explicitly by
 \begin{gather}
 \theta_\gamma(p)=p\sw{0}\gamma^{-1}(p\sw{1})\ot p\sw{2},\nonumber\\
 \theta_\gamma^{-1}(b\ot c)=b\gamma(c).\label{clthe}
 \end{gather}
 \end{mprop}

\begin{mlem}\label{gaugelem}
Suppose that $P(B)^C_\gamma$ is a cleft $C$-coalgebra Galois extension, and that 
$\gamma:C\rightarrow P$ is a cleaving map on $P$. The map
\begin{equation}
\gamma'=m\circ(\Gamma\ot \gamma)\circ\comul,
\end{equation}
where $\Gamma:C\rightarrow B$ is a convolution invertible map 
called a {\em gauge transformation}, is also a cleaving map on $P$, and any other cleaving map on $P$
has this form. 
\end{mlem}

\section{Covering  of modules and algebras} 
\label{section1c4}
In this and the next section we recall basic definitions and theorems from
\cite{CalMat:Glu}. Note that the covering and gluing of modules was actually introduced in
\cite{CalMat:Assoc}.

In what follows, all algebraic objects, unless specified otherwise, are $\mK$-modules,
where $\mK$ is a unital commutative ring such that $\mK\ni2\neq 0$, $\mK\ni3\neq 0$,
and any $\mK$-module $M$ considered is such that $2M=M$, $3M=M$.

\begin{mdef}
Let $A$, $B$ be algebras and let $C$ be a coalgebra. Suppose that
$M$ is an $(A,B)$-bimodule (resp. a right $C$-comodule, an algebra, 
an algebra and a right $C$-comodule, etc.) Let $I$ be a finite index set, and let $(J_i)_{i\in I}$
be a family of sub-bimodules (resp. $C$-sub-comodules, ideals, ideals which are also right 
$C$-sub-comodules, etc.) of $M$, such that
\begin{equation}
\bigcap_{i\in I}J_i=\{0\}.
\end{equation}
Then the family  $(J_i)_{i\in I}$ is called a {\em cover} or a {\em covering} of $M$.
\end{mdef}
In what follows we will only consider finite covers, i.e., in the statement
`$(J_i)_{i\in I}$ is a covering' it should be implicitly understood that the index set $I$ is finite.
 
Observe that the quotient modules
\begin{equation}
M_i=M/J_i,\ \ \ M_{ij}=M/(J_i+J_j), \ \ \ M_{ijk}=M/(J_i+J_j+J_k), \ \ \ \ldots,
\end{equation}
are $(A,B)$-bimodules (resp.\  $C$-comodules, algebras, algebras and right 
$C$-co\-mo\-dules, etc.),
and hence, for all $i,j,k\in I$, the canonical  surjections
\begin{gather}
\pi_i:M\rightarrow M_i,\ m\mapsto m+J_i;\ \ \ \pi_{ij}:M\rightarrow M_{ij},\ m\mapsto m+J_i+J_j;\nonumber\\
\pi_{ijk}:M\rightarrow M_{ijk},\ m\mapsto m+J_i+J_j+J_k;\ \ \ \ldots;\nonumber\\
\pi^i_j:M_i\rightarrow M_{ij},\ m+J_i\mapsto m+J_i+J_j;\nonumber\\
\pi^i_{jk}:M_i\rightarrow M_{ijk},\ m+J_i\mapsto m+J_i+J_j+J_k;\ \ \ \ldots
\end{gather}
are morphisms in the respective categories. Note that, by our assumption at the beginning of this section, 
about the ground ring $\mK$,
for all $i,j,k\in I$,
 $M_{ii}=M_i=M_{iii}$, $M_{iij}=M_{ij}$, 
$M_{ij}=M_{ji}$, etc., 
and also
\begin{gather}
\pi_{ii}=\pi_{i},\ \pi_{ij}=\pi_{ji},\ \pi^i_i=M_i,\ \pi^i_j\circ\pi_i=\pi_{ij},\  
\pi^i_{jk}\circ\pi_i=\pi_{ijk}, \text{ etc.,}
\end{gather}

A module 
\begin{equation}
M\cp=\{(m_i)_{i\in I}\in\bigoplus_{i\in I}M_i\ |\ \forall_{i,j\in I}\ \pi^i_j(m_i)=\pi^j_i(m_j)\}
\end{equation}
 is called a {\em covering completion} of $M$.
 Observe that $M\cp=\ker\Psi_M$, where 
 \begin{equation}
 \Psi_M:\bigoplus_{i\in I}M_i\rightarrow\bigoplus_{i,j\in I}M_{ij},\ \ 
 (m_i)_{i\in I}\mapsto (\pi^i_j(m_i)-\pi^j_i(m_j))_{i,j\in I}.\label{Psid}
 \end{equation}
 The map
 \begin{equation}
 \cm_M:M\rightarrow M\cp,\ \ 
 m\mapsto (\pi_i(m))_{i\in I}\label{covmap}
 \end{equation}
 is clearly injective (as $\ker\cm_M=\bigcap_{i\in I}\ker\pi_i=\bigcap_{i\in I}J_i=\{0\}$). If $\cm_M$ 
 is also surjective, then the cover $(J_i)_{i\in I}$ of $M$ is called {\em a complete cover}. 
 
 Note
 that the definitions of  the module $M\cp$ and the map $\cm_M$ make   sense even if the family
 $(J_i)_{i\in I}$ is not a cover. Accordingly, in what follows, we shall make use of the term 
 covering completion $M\cp$ and the map $\cm_M$ also when $\bigcap_{i\in I}J_i\neq\{0\}$.
 In fact, $\cm_M$ is injective if and only if  $(J_i)_{i\in I}$ is a cover.
 
 If $M$ is an algebra with unit  (and $J_i$, $i\in I$, are ideals), then $M\cp$ is an algebra, with unit
 $(1_{M_i})_{i\in I}$ and
 \begin{equation} 
 (m_i)_{i\in I}\cdot(n_j)_{j\in I}=(m_in_i)_{i\in I},\ \ 
 \text{for all } (m_i)_{i\in I}, (n_j)_{j\in I}\in M\cp.
 \end{equation}
 The map $\cm_M$ is then an algebra morphism. Similarly, if $C$ is a coalgebra, flat as a 
 $\mK$-module,
 and $M$ is a $C$-comodule (and $(J_i)_{i\in I}$ is a family of sub-comodules), then $M\cp$ is naturally
 a $C$-comodule with the coaction
 \begin{equation}
 \coact^C:M\cp\rightarrow M\cp\ot C,\ \ 
 (m_i)_{i\in I}\mapsto (m_i\sw{0}\ot m_i\sw{1})_{i\in I}.\label{cpco}
 \end{equation}
 With respect to this  coaction  $\cm_M$ is  a right $C$-colinear map.
 
 The following two propositions give criterions for a cover to be a complete one.
 
 \begin{mlem} (Proposition 1, \cite{CalMat:Glu}.) Let $M$ be a $\mK$-module, and let
 $J_1,J_2\subseteq M$ be $\mK$-sub-modules. Then the map $\cm_M:M\rightarrow M\cp$
 defined by (\ref{covmap}) is surjective. In particular, any covering by two subspaces is complete.
 \end{mlem}
 
 \begin{mlem}(Proposition 3, \cite{CalMat:Glu}.)\label{calmat3}
 Let $M$ be a $\mK$-module and let $(J_i)_{i\in I}$ be a covering of $M$.
 Assume that the index set is $I=\{1,2,\ldots,n\}$ and that, for all
 $k\in I$, the submodules $M_i$ satisfy
 \begin{equation}
 \bigcap_{i\in\{1,2,\ldots,k-1\}}(J_i+J_k)=\left(\bigcap_{i\in\{1,2,\ldots,k-1\}}J_i\right)+J_k. 
 \label{seminet}
 \end{equation}
 Then the covering $(J_i)_{i\in I}$ is complete.
 \end{mlem}
 
 The condition (\ref{seminet}) is not necessary. Note however that the closed ideals of a
 $C^\ast$-algebra form a net with respect to intersection and addition, which is stronger
 condition than (\ref{seminet}). Moreover, a similar but weaker than (\ref{seminet}) condition  is
 a necessary condition.
 \begin{mlem}(Proposition 4, \cite{CalMat:Glu}.)
 Let $M$ be a $\mK$-module, and let $(J_i)_{i\in I}$ be a complete covering of $M$.
 Then, for all $k\in I$,
 \begin{equation}
 \bigcap_{i\neq k}(J_i+J_k)=\left(\bigcap_{i\neq k}J_i\right)+J_k.  
 \end{equation}
 \end{mlem}
 
 \section{Gluing of modules and algebras}\label{glusect}
 
 Let $M_i$, $M_{ij}$, $i,j\in I$, be a finite family of modules, and let 
 $\pi^i_j:M_i\rightarrow M_{ij}$ be a family of surjective homomorphisms such that
 $M_{ii}=M_i$, $M_{ij}=M_{ji}$ and $\pi^i_i=M_i$, for all $i,j\in I$. Then the module
 \begin{equation}
 \bigoplus_{\pi^i_j}M_i=\left\{(m_i)_{i\in I}\in \bigoplus_{i\in I}M_i\ |\ \forall_{i,j\in I}\pi^i_j(m_i)
 =\pi^j_i(m_j)\right\}\label{gluingdef}
 \end{equation}
 is called a {\em gluing of the modules $M_i$ with respect to $\pi^i_j$} (Definition 3 \cite{CalMat:Glu}).
 
 Similarly as in the case of covering completions,  if the modules $M_i$, $M_{ij}$, $i,j\in I$ are
 (unital) algebras
 and maps $\pi^i_j$ are (unital) algebra maps, then gluing $\bigoplus_{\pi^i_j}M_i$ is an algebra. 
 If $C$ is a coalgebra flat as a $\mK$-module, 
 and modules $M_i$, $M_{ij}$ are right $C$-comodules and  the maps
 $\pi_j^i$ are right $C$-colinear, then $\bigoplus_{\pi^i_j}M_i$ is naturally a right $C$-comodule.
 
 \begin{mprop}(Proposition 8 \cite{CalMat:Glu}.)
 Suppose that $M=\bigoplus_{\pi^i_j}M_i$. For all $i\in I$, let 
 \begin{equation}
 p_i:\bigoplus_{\pi^k_l}M_k\rightarrow M_i,\ 
 (m_j)_{j\in I}\mapsto m_i. \label{sempr}
 \end{equation}
 Then $(\ker(p_i))_{i\in I}$ is a complete covering of $M$.
 \end{mprop}
 
 The maps $p_i$, $i\in I$, defined above are not in general surjective. The reason,
 given in  \cite{CalMat:Glu}, is that our definition of gluing (\ref{gluingdef}) does not exclude self-gluing. 
 The following proposition gives sufficient condition for surjectivity of maps $p_i$, $i\in I$.
 \begin{mprop}(Proposition 9 \cite{CalMat:Glu}.)\label{chinsurpro}
 Let $M=\bigoplus_{\pi^i_j}M_i$. Assume that the epimorphisms $\pi^i_j:M_i\rightarrow M_{ij}$ have
 the following properties.
\begin{equation}\label{calmat1as}
\text{For all } i,j,k\in I,\ \ 
\pi^i_j(\ker\pi^i_k)=\pi^j_i(\ker\pi^j_k).
\end{equation}
Define isomorphisms
\begin{gather}
\theta^{ij}_k:M_i/(\ker\pi^i_j+\ker\pi^i_k)\rightarrow M_{ij}/\pi^i_j(\ker\pi^i_k),\nonumber\\
m_i+\ker\pi^i_j+\ker\pi^i_k\mapsto \pi^i_j(m_i)+\pi^i_j(\ker\pi^i_k).\label{thetaijk}
\end{gather}
Then assume that the isomorphisms
\begin{equation}
\phi_{ij}^k=(\theta^{ij}_k)^{-1}\circ\theta^{ji}_k:M_j/(\ker\pi^j_i+\ker\pi^j_k)\rightarrow 
M_i/(\ker\pi^i_j+\ker\pi^i_k)\label{phiijkdef}
\end{equation}
satisfy
\begin{equation}
\phi_{ik}^j=\phi_{ij}^k\circ\phi_{jk}^i, \text{ for all }i,j,k\in I.\label{phiijkcoc}
\end{equation}
Let $I=\{1,2,\ldots, n\}$. If $n>3$, assume that, for all
$1\leq k<n$ and $1\leq i<k$,
\begin{equation}\label{covdeepcap}
\bigcap_{1\leq j\leq i}(\ker\pi^{k+1}_j+\ker\pi^{k+1}_{i+1})
=\left(\bigcap_{1\leq j\leq i}\ker\pi^{k+1}_j\right)+\ker\pi^{k+1}_{i+1}.
\end{equation}
Then, for all $i\in I$, the maps $p_i$ (\ref{sempr}) are surjective.
 \end{mprop}
 \noindent {\bf Remark. }Note that, for all $i,j,k\in I$, $m_j\in M_j$, 
 $\phi_{ij}^k(m_j+\ker(\pi^j_i)+\ker(\pi^j_k))=m_i+\ker(\pi^i_j)+\ker(\pi^i_k)$, where
 $m_i$ is any element of $M_i$ such that $\pi^i_j(m_i)=\pi^j_i(m_j)$.
 
 \section{Locally $C$-coalgebra Galois extensions}\label{locgalsect}
 
 In what follows we shall frequently use the following two simple observations.
 
 \begin{mlem}\label{idlemot}
 Suppose that  $A$ and $B$ are algebras, and that $\pi:A\rightarrow B$ is a surjective algebra
 morphism. Take any $N\in {}_B\Mod_B$.
Clearly $N\in{}_A\Mod_A$ with left and right $A$-actions defined by
  $a\cdot n \cdot a'=\pi(a)n\pi(a')$,
 for all $a,a'\in A$, $n\in N$. Then we can identify $N\ot_B N$ with $N\ot_A N$.
 \end{mlem}
 
 \begin{mlem}\label{surgallem}
 Suppose that $P(B)^C$ is a $C$-coalgebra Galois extension. Let $A$ be an algebra and a right 
 $C$-comodule, and suppose that $\pi:P\rightarrow A$ is an algebra and a right 
 $C$-comodule morphism.
 If $\pi(B)=A^{\co C}$, then $A(\pi(B))^C$ is a $C$-coalgebra Galois extension.
 \end{mlem}
 \begin{proof}
 It is clear that the map
 \begin{equation}
 \trans^C_A=(\pi\ot_B\pi)\circ\trans^C_P:C\rightarrow A\ot_B A\simeq A\ot_{\pi(B)} A\label{transcf}
 \end{equation}
 is the translation map on $A$, where $\trans^C_P$ is the translation map on $P$,
 and we used the identification of $A\ot_BA$ with $A\ot_{\pi(B)}A$ (Lemma~\ref{idlemot}).
 \end{proof}
 
 Observe that it is not true in general that for an arbitrary coalgebra $C$ and algebras and right 
 $C$-comodules $P$ and $A$, such that there exists an algebra surjection $\pi:P\rightarrow A$,
 we have $A^{\co C}=\pi(P^{\co C})$. As an example take $C$ being a commutative
 Hopf algebra generated by a single primitive element $x$, i.e.,
 $\comul(x)=1\ot x+x\ot 1$ and $\counit(x)=0$. Let $P$ be a free commutative algebra
 generated by two elements $b$ and $a$, and let us define a right $C$-coaction 
 $\coact^C:P\rightarrow P\ot C$ 
 as an algebra map defined by an algebra extension of the relations
 \begin{equation}
 \coact^C(b)=b\ot 1,\ \  \coact^C(a)=a\ot 1+b\ot x.
 \end{equation}
 It is clear that $P^{\co C}$ is a subalgebra of $P$ generated by the element $b$. 
 Let $A=P/(Pb)$, i.e., $A$ is a free 
 commutative algebra generated by $a+Pb$, and 
 let $\pi:P\rightarrow A$ be the canonical surjection on the quotient space. The map $\pi$ is clearly
 right $C$-colinear, and moreover, $C$ acts trivially on $A$, i.e, $A^{\co C}=A$. However,
 $\pi(P^{\co C})=\mK$.
 
 Unfortunately $P$ is not a $C$-Hopf Galois extension, and we were unable either to prove that
 $A^{\co C}=\pi(P^{\co C})$
 when $P$ is a $C$-coalgebra Galois extension,  nor to produce
 a counterexample. The following lemma, however, shows that, under a not very restrictive condition,
 $A^{\co C}=\pi(P^{\co C})$ when $P$ is a cleft extension.
 \begin{mlem}\label{cleftcoinvslem}
 Suppose that $P(B)^C$ is a cleft $C$-coalgebra Galois extension, 
 and that $\pi:P\rightarrow A$ is a surjective algebra and
 right $C$-comodule morphism. Let $\gamma_P:C\rightarrow P$ be a cleaving map in $P$.
  If $\pi(1\sw{0}\gamma_P^{-1}(1\sw{1}))$ has a right inverse in $\pi(B)$, then 
 $A(\pi(B))^C$ is a cleft $C$-coalgebra Galois extension.
 \end{mlem}  
 \begin{proof}
 The map $\gamma_A=\pi\circ\gamma_P:C\rightarrow A$ is right $C$-colinear as the composition of
  two $C$-colinear maps, and it is convolution invertible, with the convolution inverse given by
 $\gamma_A^{-1}=\pi\circ\gamma_P^{-1}$. Moreover, since $\pi$ is surjective, for all $a\in A$,
 there exists  $p\in P$ such that $a=\pi(p)$, and then
 \begin{equation*}
 a\sw{0}\gamma_A^{-1}(a\sw{1})=\pi(p\sw{0}\gamma_P^{-1}(p\sw{1}))\in \pi(B)\subseteq A^{\co C}.
 \end{equation*}
 Therefore, by Proposition~\ref{brzcl}, it remains to prove that $A^{\co C}=\pi(B)$.
 Consider the map
 \begin{equation}
 \xymatrix{
 B\ot C\ar[rrr]^{b\ot c\mapsto b\gamma_P(c)} &&& P\ar[rr]^\pi && A
 }
 \end{equation}
 which is surjective by Proposition~\ref{brzcl}. Therefore, in particular, for all $s\in A^{co C}$,
 there exists $\sum_i b_i\ot c_i\in B\ot C$, such that 
 $s=\pi(\sum_ib_i\gamma_P(c_i))=\sum_i\pi(b_i)\gamma_A(c_i)$.
 The $C$-coinvariants of $A$ are characterised by the property $\coact^C(s)=s\coact^C(1)$, 
 therefore,
 \begin{equation*}
 \sum_i\pi(b_i)\gamma_A(c_i)1\sw{0}\ot 1\sw{1}=\sum_i\pi(b_i)\gamma_A(c_i\sw{1})\ot c_i\sw{2}.
 \end{equation*}
 Applying $\mul\circ(P\ot\gamma_A^{-1})$ to both sides of the above equation yields
 \begin{equation*}
 \sum_i\pi(b_i)\gamma_A(c_i)1\sw{0}\gamma_A^{-1}(1\sw{1})=\sum_i\pi(b_i)\counit(c_i),
 \end{equation*}
 hence, if $1\sw{0}\gamma_A^{-1}(1\sw{1})$ has a right inverse $R\in\pi(B)$, then
 \begin{equation*}
 s=\sum_i\counit(c_i)\pi(b_i)R\in\pi(B),
 \end{equation*}
 which ends the proof.
 \end{proof}
 
 Note that if $P(B)^C_{e,\gamma_P}$ is an $e-copointed$ cleft $C$-coalgebra Galois extension, then
 $1\sw{0}\gamma^{-1}_P(1\sw{1})=\gamma^{-1}_P(e)$, which is invertible in $B$ with
 $(\gamma^{-1}_P(e))^{-1}=\gamma_P(e)\in B$.
 
 \begin{mdef}
 A pair $(P(B)^C, (J_i)_{i\in I})$ is called a {\em locally $C$-coalgebra Galois extension}
 if the following conditions are satisfied.
 \begin{enumerate}
 \item $P$ is an algebra and a right $C$-comodule and $B=P^{\co C}$.
 \item The family $(J_i)_{i\in I}$ of ideals and right $C$-subcomodules of 
 $P$ is a complete cover of the algebra $P$.
 \item For all $i\in I$,  $\pi_i(B)=P_i{}^{\co C}$, and $P_i(\pi_i(B))^C$ is a $C$-coalgebra Galois 
 extension.
 \item For all $i,j\in I$, $\pi_{ij}(B)=P_{ij}^{\co C}$.
 \end{enumerate}
 \end{mdef}
 Note that while $(B\cap J_i)_{i\in I}$ is a cover of $B$ it does not need to be a complete cover. 
 Indeed, in general $B\cap J_i+B\cap J_j\neq B\cap(J_i+J_j)$, therefore
 $\bar{\pi}^i_j\neq\left.\pi^i_j\right|_{B/(B\cap J_i)}$, where
 $\bar{\pi}^i_j:B/(B\cap J_i)\rightarrow B/(B\cap J_i+B\cap J_j)$,
 $b+B\cap J_i\mapsto b+B\cap J_i+B\cap J_j$.
  
 The following lemma is very  technical and apparently obvious. 
  However,  we shall  make use of it  several times in critical places and we want to
 state it explicitly.
 \begin{mlem}\label{idlem}
 Let $I$, $J$ be index sets and suppose that $C$, $M_i$, $i\in I$, $N_j$, $j\in J$,  are $\mK$-modules.
 There is a well known canonical identification
 \begin{equation*}
 \vartheta_M:(\bigoplus_{i\in I}M_i)\ot C\rightarrow \bigoplus_{i\in I}(M_i\ot C),\ 
 (m_i)_{i\in I}\ot c\mapsto (m_i\ot c)_{i\in I},
 \end{equation*}
 and similarly  
 $\vartheta_N:(\bigoplus_{j\in J}N_j)\ot C\simeq \bigoplus_{j\in J}(N_j\ot C)$.
 Let $F^i_j:M_i\rightarrow N_j$, $i\in I$, $j\in J$ be  a family of $\mK$-linear morphisms. Define maps
 \begin{gather}F:\bigoplus_{i\in I}M_i\rightarrow \bigoplus_{j\in J}N_j,\ \ 
 (m_i)_{i\in I}\mapsto (\sum_{i\in I}F^i_j(m_i))_{j\in J},\\
 G:\bigoplus_{i\in I}(M_i\ot C)\rightarrow\bigoplus_{j\in J}(N_{j}\ot C),\  
 (m_i\ot c_i)_{i\in I}\mapsto (\sum_iF^i_j(m_i)\ot c_i)_{j\in J}. 
 \end{gather}
Then $G\circ\vartheta_M=\vartheta_N\circ(F\ot C)$.
 \end{mlem}
  
 \begin{mlem}
 Suppose that $(P(B)^C, (J_i)_{i\in I})$ is a locally $C$-coalgebra Galois extension. 
 For all $i\in I$, denote 
 by $\trans_i:C\rightarrow P_i\ot_B P_i$ the translation map in $P_i$.
 Then, for all $i,j\in I$,
 \begin{equation}
 (\pi^i_j\ot_B\pi^i_j)\circ\trans_i=(\pi^j_i\ot_B\pi^j_i)\circ\trans_j.\label{treqcov}
 \end{equation}
 \end{mlem}
 \begin{proof}
 By Lemma~\ref{surgallem}, both sides of (\ref{treqcov}) are translation maps in 
 $P_{ij}$. But the translation map, if it exists, is unique, hence the equality.
 \end{proof}
 
 We use an indexed summation notation for the translation map. For all $i\in I$ and
 $c\in C$,
 \begin{equation}
 \trans_i(c)=c\uwa{i}\ot_B c\uwb{i},
 \end{equation}
 implicit summation (not over $i$ though!) is understood.
 
 \begin{mprop}
 Let $(P(B)^C, (J_i)_{i\in I})$ be a locally $C$-coalgebra Galois extension and suppose
 that $C$ is flat as a $\mK$-module. Then $P$ is a $(P,C)_\psi$-entwined module with
 \begin{equation}
 \psi:C\ot P\rightarrow P\ot C,\ \ c\ot p\mapsto(\cm^{-1}_P\ot C)((\psi_i(c\ot\pi_i(p)))_{i\in I}),
 \end{equation} 
 where, for all $i\in I$, $\psi_i$ is the canonical entwining on $P_i$ (\ref{canent}).
 \end{mprop}
 \begin{proof}
 First we prove that the map $\psi$ is well defined. Using (\ref{treqcov}) we show that,
 for all $i,j\in I$, $c\in C$, $p\in P$,
 \begin{equation}
  (\pi^i_j\ot C)\circ\psi_i(c\ot\pi_i(p))=(\pi^j_i\ot C)\circ\psi_j(c\ot\pi_j(p)).\label{psipi}
 \end{equation}
 Indeed,
  \begin{multline*}
 (\pi^i_j\ot C)\circ\psi_i(c\ot\pi_i(p))
 =(\pi^i_j\ot C)(c\uwa{i}\coact^C(c\uwb{i}\pi_i(p)))\\
 =\pi^i_j(c\uwa{i})\coact^C(\pi^i_j(c\uwb{i})\pi^i_j(\pi_i(p)))
 =\pi^j_i(c\uwa{j})\coact^C(\pi^j_i(c\uwb{j})\pi^j_i(\pi_j(p)))\\
 =(\pi^j_i\ot C)(c\uwa{j}\coact^C(c\uwb{j}\pi_j(p)))
 =(\pi^j_i\ot C)\circ\psi_j(c\ot\pi_j(p)).
 \end{multline*}
 Define the map
 \begin{equation*}
 \bar{\psi}:C\ot P\rightarrow \bigoplus_{i\in I}(P_i\ot C),\ \ 
 c\ot p\mapsto (\psi_i(c\ot\pi_i(p)))_{i\in I}.
 \end{equation*}
By (\ref{psipi}) and Lemma~\ref{idlem}, 
$\mIm(\bar{\psi})\subseteq\ker(\Psi_P\ot C)=\ker(\Psi_P)\ot C=P\cp\ot C$, where
$\Psi_P$ is as in (\ref{Psid}), and we used the flatness of $C$ and the definition of $P\cp=\ker(\Psi_P)$.
Hence the map $\psi=(\cm^{-1}_P\ot C)\circ\bar{\psi}$ is well defined.

In order to distinguish between different entwining maps we use 
indexed summation notation $\psi_i(c\ot \pi_i(p))=\pi_i(p)_{\alpha_i}\ot c^{\alpha_i}$,
for all $i\in I$, $p\in P$, $c\in C$, 
together with the usual notation $\psi(c\ot p)=p_\alpha\ot c^\alpha$, an implicit summation understood. 
We need to check whether $\psi$ satisfies conditions (\ref{bowtie}).
Indeed, for any $c\in C$,
\begin{equation*}
\psi(c\ot 1_P)=(\cm^{-1}_P\ot C)((\psi_i(c\ot 1_{P_i}))_{i\in I})(\cm^{-1}_P\ot C)((1_{P_i})_{i\in I}\ot c)=1_P\ot c,
\end{equation*}
where we used that $\cm_P$ (hence $\cm^{-1}_P$) is a unital map. Similarly, for all
$c\in C$, $p\in P$,
\begin{multline*}
(P\ot \counit)\circ\psi(c\ot p)=(\cm^{-1}_P\ot\counit)((\psi_i(c\ot \pi_i(p)))_{i\in I})\\
=(\cm^{-1}_P\ot C)(((P_i\ot\counit)\circ\psi_i(c\ot\pi_i(p)))_{i\in I})
=\cm^{-1}_P((\pi_i(p))_{i\in I}\counit(c))=p\counit(c).
\end{multline*}
Observe that, for all $i\in I$,
\begin{equation}
(\pi_i\ot C)\circ\psi=\psi_i\circ(C\ot\pi_i).
\end{equation}
Explicitly, for all $c\in C$, $p\in P$, $i\in I$, 
$\pi_i(p_\alpha)\ot c^\alpha=\pi_i(p)_{\alpha_i}\ot c^{\alpha_i}$.
Hence, 
for all $c\in C$, $p,p'\in P$,
\begin{multline*}
(pp')_\alpha\ot c^\alpha=(\cm^{-1}_P\ot C)((\pi_i(pp')_{\alpha_i}\ot c^{\alpha_i})_{i\in I})\\
=(\cm^{-1}_P\ot C)((\pi_i(p)_{\alpha_i}\pi_i(p')_{\beta_i}\ot c^{\alpha_i\beta_i})_{i\in I})\\
=(\cm^{-1}_P\ot C)((\pi_i(p_\alpha)\pi_i(p'{}_\beta)\ot c^{\alpha\beta})_{i\in I})
=p_\alpha p'{}_\beta\ot c^{\alpha\beta}.
\end{multline*}
Similarly, for all $p\in P$, $c\in C$,
\begin{multline*}
p_\alpha\ot c^\alpha\sw{1}\ot c^\alpha\sw{2}
=(\cm^{-1}_P\ot C\ot C)((\pi_i(p_\alpha)\ot c^\alpha\sw{1}\ot c^\alpha\sw{2})_{i\in I})\\
=(\cm^{-1}_P\ot C\ot C)((\pi_i(p)_{\alpha_i}\ot c^{\alpha_i}\sw{1}\ot c^{\alpha_i}\sw{2})_{i\in I})\\
=(\cm^{-1}_P\ot C\ot C)((\pi_i(p)_{\alpha_i\beta_i}\ot c\sw{1}{}^{\beta_i}\ot c\sw{2}{}^{\alpha_i})_{i\in I})\\
=(\cm^{-1}_P\ot C\ot C)((\pi_i(p_{\alpha\beta})\ot c\sw{1}{}^\beta\ot c\sw{2}{}^\alpha)_{i\in  I})
=p_{\alpha\beta}\ot c\sw{1}{}^\beta\ot c\sw{2}{}^\alpha.
\end{multline*}
It remains to prove that $P$ is an entwined module. Indeed, using the coaction (\ref{cpco}), for all
$p,p'\in P$,
\begin{multline*}
\coact^C(pp')=\coact^C(\cm^{-1}_P\circ\cm_P(pp'))
=(\cm^{-1}_P\ot C)\circ\coact^C\circ\cm_P(pp')\\
=(\cm^{-1}_P\ot C)((\coact^C(\pi_i(p)\pi_i(p')))_{i\in I})
=(\cm^{-1}_P\ot C)((\pi_i(p)\sw{0}\psi_i(\pi_i(p)\sw{1}\ot \pi_i(p')))_{i\in I})\\
=(\cm^{-1}_P\ot C)((\pi_i(p\sw{0})\psi_i(p\sw{1}\ot \pi_i(p')))_{i\in I})
=p\sw{0}(\cm^{-1}_P\ot C)((\psi_i(p\sw{1}\ot \pi_i(p')))_{i\in I})\\
=p\sw{0}\psi(p\sw{1}\ot p').
\end{multline*}
 \end{proof}
 
 Although a locally coalgebra Galois extension is built out of Galois extensions it is not necessarily 
  a (global) coalgebra Galois extension. The next theorem, which is the main result of this section,
  gives (sufficient and necessary in the case of a flat coalgebra) conditions for when a
  locally coalgebra Galois extension is a global coalgebra Galois extension.
  
 \begin{mth}\label{thelocgalm}
 Let $(P(B)^C, (J_i)_{i\in I})$ be  a locally $C$-coalgebra Galois extension.
 \begin{enumerate}
\item If $(\ker(\pi_i\ot_B\pi_i))_{i\in I}$ is a complete covering of $P\ot_B P$
 and 
 \begin{equation}
 \ker(\pi_{ij}\ot_B\pi_{ij})=\ker(\pi_i\ot_B\pi_i)+\ker(\pi_j\ot_B\pi_j),\label{covsum}
 \end{equation}
  then $P(B)^C$ is a $C$-coalgebra Galois extension. 
  \item Suppose that the coalgebra $C$ 
  is flat as a $\mK$ module. The family $(\ker(\pi_i\ot_B\pi_i))_{i\in I}$ is a 
  cover of $P\ot_B P$ if and only if $\can^C_P$ is injective.
  \item If  $P(B)^C$ is a 
  $C$-coalgebra Galois extension and $C$ is flat as a $\mK$-module, 
  then the condition (\ref{covsum}) is satisfied and $(\ker(\pi_i\ot_B\pi_i))_{i\in I}$
  is a complete covering of $P\ot_BP$. 
  \end{enumerate}
 \end{mth}
 \begin{proof}
  Observe that, for all $i,j\in I$, modules
 $P_i\ot_BP_i$ and $P_{ij}\ot_BP_{ij}$, as the images of the surjective maps
 $\pi_i\ot_B\pi_i$ and $\pi_{ij}\ot_B\pi_{ij}$ respectively, can  be identified with the respective
 quotient spaces $P\ot_BP/\ker(\pi_i\ot_B\pi_i)$ and 
 $P\ot_BP/\ker(\pi_{ij}\ot_B\pi_{ij})$. Under this identification the maps $\pi_i\ot_B\pi_i$ and
 $\pi_{ij}\ot_B\pi_{ij}$ can be viewed as quotient maps.
 
 \noindent {\bf 1).} Suppose first that $(\ker(\pi_i\ot_B\pi_i))_{i\in I}$ is a complete cover of $P\ot_BP$ and  that relation (\ref{covsum})
  is satisfied. The condition (\ref{covsum}) means that
 $P_{ij}\ot_BP_{ij}$ can be identified with
 $P\ot_BP/(\ker(\pi_i\ot_B\pi_i)+\ker(\pi_j\ot_B\pi_j))$. The
  map $\pi^i_j\ot_B\pi^i_j$ is surjective and
  $(\pi^i_j\ot_B\pi^i_j)\circ(\pi_i\ot_B\pi_i)=\pi_{ij}\ot_B\pi_{ij}$. Therefore, the map $\pi^i_j\ot_B\pi^i_j$
  can be viewed as a quotient map
 \begin{gather}
 \pi^i_j\ot_B\pi^i_j:P_i\ot_BP_i\rightarrow P_{ij}\ot_BP_{ij},\nonumber\\ 
 x+\ker(\pi_i\ot_B\pi_i)\mapsto x+\ker(\pi_i\ot_B\pi_i)+\ker(\pi_j\ot_B\pi_j).
 \end{gather}
 It follows that the covering completion of $P\ot_BP$ can be equivalently defined as 
 \begin{equation}
 (P\ot_BP)\cp=\{(x_i)_{i\in I}\in\bigoplus_{i\in I}P_i\ot_BP_i\ |\ 
 \forall_{i,j\in I}(\pi^i_j\ot_B\pi^i_j)(x_i)=(\pi^j_i\ot_B\pi^j_i)(x_j)\},\label{cpnewdef}
 \end{equation}
 and then, by assumption, the map (cf. (\ref{covmap}))
 \begin{equation}
 \cm_{P\ot_BP}:P\ot_BP\mapsto (P\ot_BP)\cp,\ \ 
 x\mapsto ((\pi_i\ot_B\pi_i)(x))_{i\in I}
 \end{equation}
 is bijective. 
Define a map 
\begin{equation}
\transc:C\rightarrow (P\ot_BP)\cp,\  \ 
c\mapsto (\trans_i(c))_{i\in I},
\end{equation}
where $\trans_i:C\mapsto P_i\ot_BP_i$ is the  translation map on $P_i$, $i\in I$. 
Equation (\ref{treqcov})
ensures that this map has image in $(P\ot_B P)\cp$. We claim that the map
\begin{equation}\label{caninvlocgalexplfor}
(\can^C_P){}^{-1}:P\ot C\rightarrow P\ot_BP,\ \ 
p\ot c\mapsto p\cm^{-1}_{P\ot_BP}\circ\transc(c)
\end{equation}
is the inverse of the canonical map $\can^C_P:P\ot_BP\rightarrow P\ot C$ of $P$. Indeed, denote
\begin{equation}
\bar{\can}:(P\ot_BP)\cp\rightarrow\bigoplus_{i\in I}P_i\ot C,\ \ 
(p_i\ot_Bq_i)_{i\in I}\mapsto (p_iq_i\sw{0}\ot q_i\sw{1})_{i\in I}.\label{barcan}
\end{equation}
It is easy to see that
\begin{equation}
(\cm_P\ot C)\circ\can^C_P=\bar{\can}\circ\cm_{P\ot_B P},\label{cankcan}
\end{equation}
and therefore, for all $p\in P$ and $c\in C$,
\begin{multline*}
\can^C_P\circ(\can^C_P){}^{-1}(p\ot c)=(\cm^{-1}_P\ot C)\circ\bar{\can}\circ\cm_{P\ot_BP}(
p\cm^{-1}_{P\ot_BP}\circ\transc(c))\\
=(\cm^{-1}_P\ot C)\circ\bar{\can}((\pi_i(p)\trans_i(c))_{i\in I})
=(\cm^{-1}_P\ot C)((\pi_i(p)\ot c)_{i\in I})=p\ot c.
\end{multline*}
Similarly, for all $p,q\in P$,
\begin{multline*}
(\can^C_P){}^{-1}\circ\can^C_P(p\ot_B q)
=(\can^C_P){}^{-1}(pq\sw{0}\ot q\sw{1})
=pq\sw{0}\cm^{-1}_{P\ot_BP}\circ\transc(q\sw{1})\\
=\cm^{-1}_{P\ot_BP}((\pi_i(p)\pi_i(q)\sw{0}\trans_i(\pi_i(q)\sw{1}))_{i\in I})
=\cm^{-1}_{P\ot_BP}(\pi_i(p)\ot_B\pi_i(q))=p\ot_Bq,
\end{multline*}
where in the fourth equality we used (\ref{trpp}), and for the third equality we observed
that 
$p\cm^{-1}_{P\ot_BP}((x_i)_{i\in I})=\cm^{-1}_{P\ot_BP}((\pi_i(p)x_i)_{i\in I})$,
for all $p\in P$, $(x_i)_{i\in I}\in (P\ot_BP)\cp$.

\noindent {\bf 2.)} 
Suppose that $C$ is flat as  a $\mK$-module. 
It is clear that
\begin{equation*}
\bar{\can}((P\ot_BP)\cp)\subseteq \ker(\Psi_P\ot C)=P\cp\ot C,
\end{equation*} 
where the last equality follows
by the flatness of $C$ and the definition of $P\cp$ (\ref{cpnewdef}) and $\Psi_P$ (\ref{Psid}). Define 
\begin{equation}
\can\cp:(P\ot_BP)\cp\rightarrow P\cp\ot C,\ \ 
x\mapsto\bar{\can}(x).
\end{equation}
It is clear that $\can\cp$ is invertible with the inverse
\begin{equation}
(\can\cp){}^{-1}:P\cp\ot C\rightarrow (P\ot_BP)\cp,\ \ 
(p_i)_{i\in I}\ot C\rightarrow ((\can_{P_i}^C)^{-1}(p_i\ot c))_{i\in I}.
\end{equation}
Using (\ref{cankcan}), and noticing that $\cm_P\ot C$ is bijective, it is easy to see
that $\cm_{P\ot_BP}$ is injective if and only if $\can^C_P$ is. But the injectivity
of $\cm_{P\ot_BP}$ is equivalent to  $(\ker(\pi_i\ot_B\pi_i))_{i\in I}$ being a cover.

\noindent {\bf 3.} For brevity, we denote $\can_i=\can^C_{P_i}$, $\can_{ij}=\can^C_{P_{ij}}$.
By Lemma~\ref{surgallem}, for all $i,j\in I$,
\begin{gather}
\can_i\circ(\pi_i\ot_B\pi_i)=(\pi_i\ot C)\circ\can^C_P,\label{cani}\\
\can_{ij}\circ(\pi_{ij}\ot_B\pi_{ij})=(\pi_{ij}\ot C)\circ\can^C_P.\label{canij}
\end{gather}
Hence, as the maps 
$\can^C_P$, $\can_i$, $\can_{ij}$ are bijective and $C$ is flat as a $\mK$-module, it follows 
that
\begin{multline*}
\ker(\pi_{ij}\ot_B\pi_{ij})=(\can^C_P){}^{-1}(\ker(\pi_{ij}\ot C))
=(\can^C_P){}^{-1}(\ker(\pi_{ij})\ot C)\\
=(\can^C_P){}^{-1}(\ker(\pi_i)\ot C+\ker(\pi_j)\ot C)\\
=(\can^C_P){}^{-1}(\ker(\pi_i\ot C))+(\can^C_P){}^{-1}(\ker(\pi_j\ot C))\\
=\ker(\pi_i\ot_B\pi_i)+\ker(\pi_j\ot_B\pi_j).
\end{multline*}

Furthermore, by (\ref{cankcan}), as the maps $\can\cp$, $\cm_P\ot C$ and $\can^C_P$ are bijective, 
  we obtain that
\begin{equation*}
\cm_{P\ot_BP}=(\can\cp){}^{-1}\circ(\cm_P\ot C)\circ\can^C_P
\end{equation*}
is invertible, hence the family $(\ker(\pi_i\ot_B\pi_i))_{i\in I}$ is a complete cover of $P\ot_BP$.
 \end{proof}
 
 \begin{mlem}\label{covsumlem}
 Let $(P(B)^C,(J_i)_{i\in I})$ be a locally $C$-coalgebra Galois extension, and suppose that the ground
 ring $\mK$ is a field. Then the condition (\ref{covsum}) is satisfied. 
 \end{mlem}
 \begin{proof}
 For all $i,j,k\in I$, the  diagrams 
 \begin{equation}
 \xymatrix{
 0 \ar[r] & PdBP \ar[rr] \ar[dd]^-{\left.\pi_i\ot\pi_i\right|_{PdBP}} && P\ot P 
 \ar[rr]^{\Theta_P} \ar[dd]^{\pi_i\ot\pi_i} && P\ot_BP\ar[r]\ar[dd]^{\pi_i\ot_B\pi_i} & 0\\
 & && && &\\
 0\ar[r] & P_idBP_i\ar[rr] && P_i\ot P_i \ar[rr]^{\Theta_{P_i}} && P_i\ot_BP_i\ar[r] & 0
 }\label{largei}
 \end{equation}
 and
 \begin{equation}\xymatrix{
 0 \ar[r] & PdBP \ar[rr] \ar[dd]^-{\left.\pi_{ij}\ot\pi_{ij}\right|_{PdBP}} && P\ot P 
 \ar[rr]^{\Theta_P} \ar[dd]^{\pi_{ij}\ot\pi_{ij}} && P\ot_BP\ar[r]\ar[dd]^{\pi_{ij}\ot_B\pi_{ij}} & 0\\
 & && && &\\
 0\ar[r] & P_{ij}dBP_{ij}\ar[rr] && P_{ij}\ot P_{ij} \ar[rr]^{\Theta_{P_{ij}}} && P_{ij}\ot_BP_{ij}\ar[r] & 0\ ,
 }\label{largeij}
 \end{equation}
 where $\Theta_M:M\ot M\rightarrow M\ot_BM$, $M\in{}_B\Mod_B$ is a natural surjection on the
 quotient space, and $d$ denotes the  universal differential, 
 are clearly commutative and have exact rows.
 
Since the maps $\left.\pi_i\ot\pi_i\right|_{PdBP}$ and $\left.\pi_{ij}\ot\pi_{ij}\right|_{PdBP}$ 
are surjective, the application 
 of the 
 Snake Lemma to the above diagrams yields that the maps
 \begin{equation}
 \left.\Theta_P\right|_{\ker(\pi_i\ot\pi_i)}:\ker(\pi_i\ot\pi_i)\rightarrow \ker(\pi_i\ot_B\pi_i)
 \end{equation} 
 and
  \begin{equation}
 \left.\Theta_P\right|_{\ker(\pi_{ij}\ot\pi_{ij})}:\ker(\pi_{ij}\ot\pi_{ij})\rightarrow \ker(\pi_{ij}\ot_B\pi_{ij})
 \end{equation} 
 are well defined and surjective. Observe that as $\mK$ is a field, for all $i,j\in I$,
 \begin{equation*}
 \ker(\pi_i\ot\pi_i)=\ker(\pi_i)\ot P+P\ot\ker(\pi_i),
 \end{equation*}
 and
 \begin{multline*}
 \ker(\pi_{ij}\ot\pi_{ij})=\ker(\pi_{ij})\ot P+P\ot\ker(\pi_{ij})\\
 =\ker(\pi_i)\ot P+\ker(\pi_j)\ot P+P\ot\ker(\pi_i)+P\ot\ker(\pi_j)\\
 = \ker(\pi_i\ot\pi_i)+ \ker(\pi_j\ot\pi_j).
 \end{multline*}
 Hence, for all $i,j\in I$,
 \begin{multline*}
 \ker(\pi_{ij}\ot_B\pi_{ij})=\Theta_P(\ker(\pi_{ij}\ot\pi_{ij}))\\
 =\Theta_P(\ker(\pi_i\ot\pi_i))+ \Theta_P(\ker(\pi_j\ot\pi_j))
 =\ker(\pi_i\ot_B\pi_i)+\ker(\pi_j\ot_B\pi_j).
 \end{multline*}
 \end{proof}
 
 Therefore, Lemma~\ref{covsumlem} implies that, when working over a field, which is probably
 the most interesting case from a non-commutative geometry point of view, to verify whether
 a locally Galois extension is globally Galois, it suffices to check whether the covering
 $(\ker(\pi_i\ot_B \pi_i))_{i\in I}$ is a complete covering of $P\ot_B P$. More precisely we can state:
 
 \begin{mcor}\label{corfieldcov}
 Suppose that the ground ring $\mK$ is a field and that
 $(P(B)^C,(J_i)_{i\in I})$ is a locally $C$-coalgebra Galois extension. Then $P(B)^C$ is a 
 $C$-coalgebra Galois extension if and only if $(\ker(\pi_i\ot_B\pi_i))_{i\in I}$ is a complete
 cover of $P\ot_BP$.
 \end{mcor}
 
 In view of Corollary~\ref{corfieldcov} it is important to study when a cover is a complete cover.
 
 \begin{mlem}\label{capmullem}
 Let $B$ be an algebra, and let $(K_i)_{i\in I}$ be a family of ideals of $B$. Denote the quotient spaces
  by $B_i=B/K_i$, $B_{ij}=B/(K_i+K_j)$, $i,j\in I$,  and by
 \begin{equation*}
 \pi_i:B\rightarrow B_i,\ \ \pi_{ij}:B\rightarrow B_{ij},\  \ 
 \pi^i_j:B_i\rightarrow B_{ij},\ \ 
 \ i,j\in I
 \end{equation*}
 the canonical surjections.
 Suppose that $M\in{}_B\Mod$, $M_i\in{}_{B_i}\Mod$, $M_{ij}\in {}_{B_{ij}}\Mod$,
 $i,j\in I$, is a family of modules such that $M_{ij}=M_{ji}$, for all $i,j\in I$, and that
 \begin{equation*}
 \chi_i:M\rightarrow M_i,\ \ 
 \chi_{ij}:M\rightarrow M_{ij},\ \
 \chi^i_j:M_i:\rightarrow M_{ij},\ \ \ 
 i,j\in I,
 \end{equation*}
 is a family of surjective $\mK$-linear morphisms with the properties
 \begin{gather}
 \ker(\chi_{ij})=\ker(\chi_i)+\ker(\chi_j),\label{kerchi}\\
 \chi^i_j\circ\chi_i=\chi_{ij}=\chi^j_i\circ\chi_j,\label{chichi}\\
 \chi_i(bm)=\pi_i(b)\chi_i(m),\ \ \chi_{ij}(bm)=\pi_{ij}(b)\chi_{ij}(m),\label{blin} 
 \end{gather}
 for all $i,j\in I$, $b\in B$, $m\in M$.
 
 Let $I=\{1,2,\ldots,n\}$, $n>2$. Suppose that, for all $2<k\leq n$,
 \begin{equation}
 \bigcap_{i\in\{1,2,\ldots,k-1\}}\chi_k(\ker(\chi_i))\subseteq\left(\prod_{i=1}^{k-1}\pi_k(K_i)\right)M_k.
 \end{equation}
 Then, for all $2<k\leq n$,
 \begin{equation}
 \bigcap_{i\in \{1,2,\ldots,k-1\}}(\ker(\chi_i)+\ker(\chi_k))
 =\left(\bigcap_{i\in \{1,2,\ldots,k-1\}}\ker(\chi_i)\right)+\ker(\chi_k).
 \end{equation}
 \end{mlem}
 \begin{proof}
 Note that by (\ref{chichi}) and (\ref{kerchi}), for all $i,j\in I$,
 $
 \ker(\chi^i_j)=\chi_i(\ker(\chi_j))
 $.
 Therefore, for all $2<k\leq n$,
 \begin{multline*}
 \bigcap_{i\in \{1,2,\ldots,k-1\}}(\ker(\chi_i)+\ker(\chi_k))
 =\bigcap_{i\in \{1,2,\ldots,k-1\}}\ker(\chi_{ik})
 =\ker\left(\bigoplus_{i=1}^{k-1}\chi_{ik}\right)\\
 =\ker\left(\left(\bigoplus_{i=1}^{k-1}\chi^k_i \right)\circ\chi_k\right)
 =(\chi_k)^{-1}\left(\ker\left(\bigoplus_{i=1}^{k-1}\chi^k_i \right)\right)\\
 =(\chi_k)^{-1}\left(\bigcap_{i\in \{1,2,\ldots,k-1\}}\ker(\chi^k_i)\right)
 =(\chi_k)^{-1}\left(\bigcap_{i\in \{1,2,\ldots,k-1\}}\chi_k(\ker(\chi_i))\right)\\
 \subseteq (\chi_k)^{-1}\left(\left(\prod_{i=1}^{k-1}\pi_k(K_i)\right)M_k\right)
 =\left(\prod_{i=1}^{k-1}K_i\right)M+\ker(\chi_k)\\
 \subseteq
 \left(\bigcap_{i\in \{1,2,\ldots,k-1\}}\ker(\chi_i)\right)+\ker(\chi_k).
 \end{multline*}
 The inclusion relation in the opposite direction is always satisfied.
 \end{proof}
 
 \begin{mprop} Suppose that the ground ring $\mK$ is a field. 
 Let $(P(B)^C,(J_i)_{i\in I})$ be a locally $C$-coalgebra Galois extension, and let 
 $\pi_i:P\rightarrow P_i=P/J_i$, etc., be the surjections on the quotient spaces. Let 
 $K_i=B\cap J_i$, $i\in I$. Suppose that $I=\{1,2,\ldots,n\}$ and, for all
 $2<k\leq n$,
 \begin{equation}
 \bigcap_{i\in \{1,2,\ldots,k-1\}}\pi_k(J_i)\subseteq\left(\prod_{i=1}^{k-1}\pi_k(K_i)\right)P_k.
 \end{equation}
 Then if the family $(\ker(\pi_i\ot_B\pi_i))_{i\in I}$ is a cover of $P\ot_BP$, then it is a complete cover.
 \end{mprop}
 \begin{proof}
 For all $i,j\in I$, denote $B_i=B/K_i$, $B_{ij}=B/K_{ij}$ and define maps
 \begin{gather}
 \chi_i=\can_i\circ(\pi_i\ot_B\pi_i):P\ot_BP\rightarrow P_i\ot C,\\
 \chi_{ij}=\can_{ij}\circ(\pi_{ij}\ot_B\pi_{ij}):P\ot_BP\rightarrow P_{ij}\ot C,\\
 \chi^i_j=\pi^i_j\ot C:P_i\ot C\rightarrow P_{ij}\ot C,
 \end{gather}
 which clearly satisfy the conditions (\ref{chichi})-(\ref{blin}). Moreover, as the maps 
 $\can_i$ and $\can_{ij}$ are 
 bijective, it follows that
$ \ker(\chi_i)=\ker(\pi_i\ot_B\pi_i)$ and
\begin{equation*}\ker(\chi_{ij})=\ker(\pi_{ij}\ot_B\pi_{ij})=\ker(\pi_i\ot_B\pi_i)+\ker(\pi_j\ot_B\pi_j)
=\ker(\chi_i)+\ker(\chi_j),
\end{equation*}
where the second equality follows from Lemma~\ref{covsumlem}. Moreover, for all
$2<k\leq n$,
\begin{multline*}
\bigcap_{i\in \{1,2,\ldots,k-1\}}\chi_k(\ker(\chi_i))
=\bigcap_{i\in \{1,2,\ldots,k-1\}}\ker(\pi^k_i\ot C)
=\bigcap_{i\in \{1,2,\ldots,k-1\}}\ker(\pi^k_i)\ot C\\
=\ker\left(\bigoplus_{i=1}^{k-1}\pi^k_i\right)\ot C
=\left(\bigcap_{i\in \{1,2,\ldots,k-1\}}\pi_k(J_i)\right)\ot C
\subseteq\left(\prod_{i=1}^{k-1}\pi_k(K_i)\right)(P_k\ot C).
\end{multline*}
Therefore, by Lemma~\ref{capmullem}, for all $2<k\leq n$,
\begin{multline*}
\bigcap_{i\in \{1,2,\ldots,k-1\}}(\ker(\pi_i\ot_B\pi_i)+\ker(\pi_k\ot_B\pi_k))\\
=\left(\bigcap_{i\in \{1,2,\ldots,k-1\}}\ker(\pi_i\ot_B\pi_i)\right)+\ker(\pi_k\ot_B\pi_k),
\end{multline*}
and hence, by Lemma~\ref{calmat3}, the family 
$(\ker(\pi_i\ot_B\pi_i))_{i\in I}$ is a complete cover of $P\ot_BP$.
 \end{proof}
 
 The following lemma is probably well known, but we were unable to find the reference.
 \begin{mlem}
 Let $M$, $M'$, $M''$ be $\mK$-modules and let 
 $K,L\subseteq M$, $K',L'\subseteq M'$, $K'',L''\subseteq M''$ be submodules.
 Suppose that
 \begin{equation*}
 f:K+L\rightarrow K'+L',\ \ 
 g:K'+L'\rightarrow K''+L''
 \end{equation*}
 are $\mK$-linear maps such that the sequences
 \begin{gather}
 \xymatrix{
 0\ar[r] & K\ar[rr]^{\left.f\right|_K} && K' \ar[rr]^{\left.g\right|_{K'}}  && K''\ar[r] & 0,
 }\label{sqk}\\
 \xymatrix{
  0\ar[r] & L\ar[rr]^{\left.f\right|_L} && L' \ar[rr]^{\left.g\right|_{L'}}  && L''\ar[r] & 0,
 }\label{sql}
 \end{gather}
 are well defined and exact. Then the sequence
 \begin{equation}\xymatrix{
 0\ar[r] & K\cap L\ar[rr]^-{\left.f\right|_{K\cap L}} && K'\cap L'
 \ar[rr]^-{\left.g\right|_{K'\cap L'}} && K''\cap L'' \ar[r] & 0
 }\label{stcap}
 \end{equation}
 is exact if and only if the sequence
 \begin{equation}
 \xymatrix{
 0\ar[r] & K+L\ar[rr]^f && K'+L' \ar[rr]^g && K''+L''\ar[r] & 0
 }\label{stsum}
 \end{equation}
 is exact.
 \end{mlem}
 \begin{proof}
 \noindent {\bf (\ref{stsum})$\Rightarrow$(\ref{stcap}).} Suppose that the sequence (\ref{stsum}) is
 exact. Clearly $\left.f\right|_{K\cap L}$ is injective and 
 $\left.g\right|_{K'\cap L'}\circ\left.f\right|_{K\cap L}=0$. Suppose that for some $m'\in K'\cap L'$,
 $g(m')=0$. By the exactness of (\ref{sqk}) and (\ref{sql}), there exist elements $k\in K$ and $l\in L$, 
 such that $f(k)=m'=f(l)$, i.e., $f(k-l)=0$. But $f$ is, by the assumption, injective, therefore
 $k=l\in K\cap L$. Hence the sequence (\ref{stcap}) is exact at $K'\cap L'$.
 
 Finally, let $m''\in K''\cap L''$. By the exactness of (\ref{sqk}) and (\ref{sql}), there exist elements 
 $k\in K'$, $l\in L'$ such that $g(k')=m''=g(l')$, hence $g(k'-l')=0$. By the exactness of (\ref{stsum})
 at $K'+L'$, there exist elements  $k\in K$, $l\in L$ such that $f(k+l)=k'-l'$, i.e.,
 $k'-f(k)=l'+f(l)\in K'\cap L'$ and $g(k'-f(k))=g(k')=m''$, hence $g(K'\cap L')=K''\cap L''$.
 
 \noindent {\bf (\ref{stcap})$\Rightarrow$(\ref{stsum}).} The map $g$ is clearly surjective,
 as $g(K'+L')=g(K')+g(L')=K''+L''$ by the exactness of  (\ref{sqk}) and (\ref{sql}). Similarly,
 for all $k\in K$, $l\in L$, $g\circ f(k+l)=g\circ f(k)+g\circ f(l)=0$.
 
 Suppose that, for some $k'\in K'$, $l'\in L'$, $g(k'+l')=0$, i.e., $g(k')=g(-l')\in K''\cap L''$.
 As, by the assumption $g(K'\cap L')=K''\cap L''$, there exists $m'\in K'\cap L'$ such that
 $g(m')=g(k')=g(-l')$, i.e., $g(k'-m')=0$ and $g(l'+m')=0$. By the exactness of 
 (\ref{sqk}) and (\ref{sql}), there exists $k\in K$, $l\in L$, such that
 $f(k)=k'-m'$ and $f(l)=l'+m'$. Therefore $f(k+l)=k'-m'+l'+m'=k'+l'$, and we have proven that 
 the sequence (\ref{stsum})  is exact at $K'+L'$.
 
 Finally, suppose that $f(k+l)=0$, for some $k\in K$, $l\in L$, i.e,
 $f(k)=f(-l)\in K'\cap L'$. As $g(f(k))=0$, we have by the exactness of (\ref{stcap}) at $K'\cap L'$ that
 there exists $m\in K\cap L$ such that $f(m)=f(k)=f(-l)$, i.e.,
 $f(k-m)=0$ and $f(l+m)=0$. However $f|_K$ and $f|_L$ are by assumption injective,
 hence $k-m=0$ and $l+m=0$. Therefore $k+l=k-m+l+m=0$ and we have proven that $f$ is injective.
 \end{proof}

\begin{mcor}
We keep the notation and assumptions from the above lemma. 
Suppose that in addition we are given exact sequence of
$\mK$-maps
\begin{equation}
\xymatrix{
0\ar[r] & M\ar[r]^s & M'\ar[r]^t & M''\ar[r] & 0,
}
\end{equation}
such that $f=s|_{K+L}$ and $g=t|_{K'+L'}$. Then $g(K'\cap L')=K''\cap L''$ if
and only if  $f(K+L)=\ker(g)$.\label{KLcor}
\end{mcor}
\begin{proof}
It is easy to see that, under the assumptions, the sequences (\ref{stcap}) and (\ref{stsum}) are exact,
apart from the conditions $g(K'\cap L')=K''\cap L''$ and $f(K+L)=\ker(g)$.
\end{proof}
 
 \begin{mlem}\label{small}
 Suppose that $\mK$ is a field. Let $f:M\rightarrow N$, $g:M\rightarrow N'$ be  $\mK$-vector space
 morphisms such that $\ker(f)\cap\ker(g)=\{0\}$. Then 
 \begin{equation}
 \ker(f\ot f)\cap\ker(g\ot g)=\ker(f)\ot \ker(g)+\ker(g)\ot\ker(f).
 \end{equation}
  \end{mlem}
  \begin{proof}
  As $\mK$ is a vector space and $\ker(f)\cap\ker(g)=\{0\}$, we can write
  \begin{gather*}
  M=\bar{M}\oplus\ker(f)\oplus\ker(g),
  \end{gather*}
  and the assertion of the lemma easily follows.
  \end{proof}

 \begin{mprop}
 Suppose that $\mK$ is a field and that $I=\{1,2\}$. Let $(P(B)^C, (J_i)_{i\in I})$ be a
 locally $C$-coalgebra Galois extension. Suppose that 
 \begin{equation}
 J_i=(B\cap J_i)P, \text{ \ \ for all }i\in I. \label{BJPcond}
 \end{equation}
 Then $P(B)^C$ is a $C$-coalgebra Galois extension if and only if
 \begin{multline}
 (\ker(\pi_1\ot\pi_1)+\ker(\pi_2\ot\pi_2))\cap PdBP\nonumber\\
 =\ker(\pi_1\ot\pi_1)\cap PdBP
 +\ker(\pi_2\ot\pi_2)\cap PdBP.
 \label{PdBPcap}
 \end{multline} 

 \end{mprop}
 \begin{proof}
 By the Snake Lemma, for all $i,j\in I$, the commutative diagrams 
 (\ref{largei}) and (\ref{largeij}) with exact rows
 induce  the exact sequences
 \begin{equation}
 \xymatrix{
 0\ar[r] & \ker(\pi_i\ot\pi_i)\cap PdBP \ar[r] & \ker(\pi_i\ot\pi_i)\ar[r] & \ker(\pi_i\ot_B\pi_i)\ar[r] & 0
 }
 \end{equation}
 and 
 \begin{equation}
  \xymatrix{
 0\ar[r] & \ker(\pi_{ij}\ot\pi_{ij})\cap PdBP \ar[r] & \ker(\pi_{ij}\ot\pi_{ij})\ar[r] & \ker(\pi_{ij}\ot_B\pi_{ij})
 \ar[r] & 0.
 }\label{thisij}
 \end{equation}
 By Lemma~\ref{covsumlem}, the sequence (\ref{thisij}) can be written as
 \begin{multline*}
 \xymatrix{
 0\ar[r] & (\ker(\pi_i\ot\pi_i)+\ker(\pi_j\ot\pi_j))\cap PdBP \ar[r] &
 \ker(\pi_i\ot\pi_i)+\ker(\pi_j\ot\pi_j)
 }\\
 \xymatrix{
 { }\ar[r] & \ker(\pi_i\ot_B\pi_i)+\ker(\pi_j\ot_B\pi_j)\ar[r] & 0.
 }
 \end{multline*}
 It follows immediately from Corollary~\ref{KLcor} that 
 \begin{equation}
 \Theta_P(\ker(\pi_i\ot\pi_i)\cap\ker(\pi_j\ot\pi_j))=\ker(\pi_i\ot_B\pi_i)\cap\ker(\pi_j\ot_B\pi_j),
 \label{thetasur}
 \end{equation}
 where $\Theta_P:P\ot P\rightarrow P\ot_B P$ is the natural surjection on the quotient space
 (cf. Lemma~\ref{covsumlem}), if and only if 
 \begin{equation}
 (\ker(\pi_i\ot\pi_i)+\ker(\pi_j\ot\pi_j))\cap PdBP=\ker(\pi_i\ot\pi_i)\cap PdBP+\ker(\pi_j\ot\pi_j)\cap PdBP.
 \label{PdBPcap}
 \end{equation} 
 Let $I=\{1,2\}$. Hence $\ker(\pi_1)\cap\ker(\pi_2)=\{0\}$, and so, by Lemma~\ref{small},
 \begin{equation}
\ker(\pi_1\ot\pi_1)\cap\ker(\pi_2\ot\pi_2)=\ker(\pi_1)\ot\ker(\pi_2)+\ker(\pi_2)\ot\ker(\pi_1).
\label{picappiK}
 \end{equation}
 Suppose that (\ref{BJPcond}) is satisfied.
 Then by (\ref{picappiK}),
$\Theta_P(\ker(\pi_i\ot\pi_i)\cap\ker(\pi_j\ot\pi_j))\Theta_P(\ker(\pi_1)\ot\ker(\pi_2)+\ker(\pi_2)\ot\ker(\pi_1))
=J_1\ot_B(J_2\cap B)P+J_2\ot_B(J_1\cap B)P=0$ as $J_1J_2, J_2J_1\subseteq J_1\cap J_2=\{0\}$.
Therefore, if in addition (\ref{thetasur}) is satisfied, then 
$\ker(\pi_1\ot_B\pi_1)\cap\ker(\pi_2\ot_B\pi_2)=\{0\}$.
 \end{proof}
 
 \section{Locally cleft extensions}\label{locclsect}
 
 \begin{mdef}
 A locally $C$-coalgebra Galois extension $(P(B)^C,(J_i)_{i\in I})$ is called a {\em locally cleft 
 extension}
 if, for all $i\in I$, the quotient modules $P_i$ are cleft $C$-coalgebra Galois extensions. It is called 
 a {\em proper locally cleft extension} if, in addition, for all $i,j\in I$,
 \begin{equation} B\cap (J_i+J_j)=B\cap J_i+B\cap J_j.\label{proper}
 \end{equation}
 \end{mdef}
 
 We adopt the following notation. We denote $P_i=P/J_i$, $P_{ij}=P/(J_i+J_j)$, etc.,  as before. 
 In addition, we have quotient modules $B_i=B/(B\cap J_i)$, $B_{ij}=B/((J_i\cap B)+(J_j\cap B))$, etc.,
 for all $i,j\in I$. We reserve the use of the Greek letter $\pi$ with various subscripts and superscripts
  to surjections
 onto the quotient modules of $B$, i.e.,
 $\pi_i:B\rightarrow B_i$, $\pi^i_j:B_i\rightarrow B_{ij}$, $b+B\cap J_i\mapsto b+B\cap J_i+B\cap J_j$,
 etc. For quotient maps on $P$, we use  sub- and superscripted letter $\chi$, i.e.,
 $\chi_i:P\rightarrow P_i$, etc.
 
 For all $i\in I$, we denote by $\gamma_i:C\rightarrow P_i$ a cleaving map on $P_i$.  Moreover,
 for all $i,j\in I$, we use the notation $\gamma^i_{ij}=\chi^i_j\circ\gamma_i:C\rightarrow P_{ij}$.
 
 Note that, for all $i\in I$, $b\in B$, $\chi_i(b)=\pi_i(b)$ if we identify $B_i$ with $B/J_i\subseteq P_i$.
 Similarly, we can naturally identify $B/(B\cap(J_i+J_j))$ with $B/(J_i+J_j)\subseteq P_{ij}$. 
 In addition, if the relation (\ref{proper}) is satisfied, we can identify $B_{ij}$ with $B/(J_i+J_j)$.
 Note that the condition (\ref{proper}) is equivalent to 
 \begin{equation}
 \chi^i_j(b)=\pi^i_j(b), \text{ for all }i,j\in I, \  b\in B,\label{proper2}
 \end{equation} 
 where we used the above identifications.
 
 The condition (\ref{proper2}) clearly implies that $(B\cap J_i)_{i\in I}$ is a complete cover of $B$.
 
 \begin{mlem}(Cf. Lemma~1 \cite{CalMat:Con}.)\label{idealcomcleft}
 Let $P(B)^C_\gamma$ be a cleft extension, and let $J$ be an ideal in $P$ such that
 $\coact^C(J)\subseteq J\ot C$. Then there exists a left ideal $K$ in $B$ such that 
 $J=K\gamma(C)$. Moreover, if the element $x=1\sw{0}\gamma^{-1}(1\sw{1})$ has a right inverse $y$ 
 in $P$
 (i.e., $xy=1_P$), and
 $Ky\subseteq K$, then $K$ is a two-sided ideal and $K=J\cap B$.
 \end{mlem}
 \begin{proof}
 Let us define
 \begin{equation}
 \label{bideal}
 K=(p\mapsto p\sw{0}\gamma^{-1}(p\sw{1}))(J).
 \end{equation}
 Note that $K\subseteq J$. Therefore, $K\gamma(C)\subseteq J$. On the other hand, for all 
 $p\in J$, $p=p\sw{0}\gamma^{-1}(p\sw{1})\gamma(p\sw{2})\in K\gamma(C)$. Hence
 $J=K\gamma(C)$.
 
 Let $b\in B$, $p\in J$, $b'=p\sw{0}\gamma^{-1}(p\sw{1})\in K$. Then
 $bb'=bp\sw{0}\gamma^{-1}(p\sw{1})=(bp)\sw{0}\gamma^{-1}((bp)\sw{1})\in K$, hence $K$ is
 a left ideal in $B$.
 
 Suppose that the element $x=1\sw{0}\gamma^{-1}(1\sw{1})$ 
 has the right inverse $y$ in $P$, and $Ky\subseteq K$. 
 As shown above, $K\subseteq B\cap J$. On the other hand, let $b\in B\cap J$. Then
 $b=\sum_ik_i\gamma(c_i)$, for some $k_i\in K$, $c_i\in C$. It follows that
 \begin{equation*}
 b1\sw{0}\ot 1\sw{1}=\sum_ik_i\gamma(c_i\sw{1})\ot c_i\sw{2},
 \end{equation*}  
 hence
 \begin{equation*}
 b1\sw{0}\gamma^{-1}(1\sw{1})=\sum_ik_i\gamma(c_i\sw{1})\gamma^{-1}(c_i\sw{2})
 =\sum_ik_i\counit(c_i),
 \end{equation*}
 and therefore $b=\sum_ik_iy\counit(c_i)\in K$. Hence $K=J\cap B$ and it follows
 that $K$ is a two-sided ideal in $B$.
 \end{proof}
 From the proof of the above lemma we immediately obtain
 \begin{mcor}\label{bcapicor}
 Suppose that $P(B)^C_\gamma$ is a cleft $C$-coalgebra Galois extension. Let $K$ be an ideal
 in $B$, and let $J=K\gamma(C)$ be an ideal in $P$. Moreover suppose that the element
$x=1\sw{0}\gamma^{-1}(1\sw{1})\in B$ has a right inverse in $B$. Then $K=J\cap B$.
 \end{mcor}
 
 \newcommand{\splce}{regular locally cleft extension{}}
 
 \begin{mdef}
 Let  $(P(B)^C,(J_i)_{i\in I})$ be a locally cleft extension. Suppose that, for all $i\in I$, the element
 $1\sw{0}\gamma_i^{-1}(1\sw{1})$ has a right inverse in $B_i$. Such a locally cleft extension we shall
 call {\em a \splce}.  
 \end{mdef}
 
 Suppose that $(P(B)^C,(J_i)_{i\in I})$ is a \splce.
 Then, by \linebreak Lemma~\ref{cleftcoinvslem},
 for all $i,j,k\in I$, $P_{ij}^{\co C}=\chi^i_j(B_i)=\chi_{ij}(B)$, 
 $P_{ijk}^{\co C}=\chi^i_{jk}(B_i)=\chi_{ijk}(B)$, etc. 
 
 For all $i,j\in I$, $\ker\chi^i_j$ is an ideal and a right $C$-subcomodule of $P_i$. 
 Therefore, by Lemma~\ref{idealcomcleft}, $\ker\chi^i_j=\bar{K}^i_j\gamma_i(C)$, where
 $\bar{K}^i_j=\ker(\chi^i_j)\cap B_i$ is an ideal in $B_i$. 
 Note that $\ker\pi^i_j=\pi_i(\ker\pi_j)\subseteq \bar{K}^i_j$.
 Define $K^i_j=\pi^i_j(\bar{K}^i_j)$. Observe that
 \begin{equation}\label{bijkij}
 P_{ij}^{\co C}=\chi^i_j(B_i)=B_i/\bar{K}^i_j=\frac{B_i/\ker(\pi^i_j)}{\bar{K}^i_j/\ker(\pi^i_j)} \pi^i_j(B_i)/\pi^i_j(\bar{K}^i_j)=B_{ij}/K^i_j.
 \end{equation}
 A locally cleft extension $(P(B)^C,(J_i)_{i\in I})$ is proper if and only if,
 for all $i,j\in I$, $K^i_j=\{0\}$ (i.e., $\bar{K}^i_j=\ker\pi^i_j$).
 Note that the properness of a locally cleft extension implies that $B_{ij}=P_{ij}^{\co C}$, and then
 \begin{equation}\label{bicapker}\ker(\chi^i_{jk})\cap B_i=\ker(\pi^i_{jk}),\ \ 
 \ker(\chi^{ij}_k)\cap B_{ij}=\ker(\pi^{ij}_k),\ \ 
 \text{for all } i,j,k\in I.
 \end{equation}
 Indeed, 
 \begin{multline*}
 \ker(\chi^i_{jk})=\chi_i(\ker\chi_{jk})=\chi_i(\ker\chi_j)+\chi_i(\ker\chi_k)
 =\ker(\chi^i_j)+\ker(\chi^i_k)\\
 =\ker(\pi^i_j)\gamma_i(C)+\ker(\pi^i_k)\gamma_i(C)=\ker(\pi^i_{jk})\gamma_i(C),
 \end{multline*}
 and
 \begin{equation*}
 \ker(\chi^{ij}_k)=\chi^i_j\circ\chi_i(\ker\chi_k)
 =\chi^i_j(\ker(\pi^i_k)\gamma_i(C))
 =\pi^i_j(\ker\pi^i_k)\gamma^i_{ij}(C)
 =\ker(\pi^{ij}_k)\gamma^i_{ij}(C).
 \end{equation*}
 Then the relations (\ref{bicapker}) follow from Corollary~\ref{bcapicor}. It follows that, for all
 $i,j,k\in I$, $B_{ijk}=\pi^{ij}_k(B_{ij})$ is isomorphic to $P_{ijk}^{\co C}=\chi^{ij}_k(B_{ij})$, and can be 
 identified with it. Under this identification, 
 \begin{equation}
 \left.\chi^i_{jk}\right|_{B_i}=\pi^i_{jk},\ \ 
 \left.\chi^{ij}_k\right|_{B_{ij}}=\pi^{ij}_k.\label{idenchiijk}
 \end{equation}
 
In what follows, we shall examine conditions for a \splce{} to be a proper locally cleft extension. 
This requires the study
of  ideals $K^i_j$, $i,j\in I$. We generalise steps of the proof of Proposition~2 \cite{CalMat:Con}.

For all $i,j\in I$, let us define an isomorphism (cf. (\ref{clthe}))
\begin{equation}
\beta^i_{ij}:P_{ij}\rightarrow P_{ij}^{\co C}\ot C,\ \ 
p_{ij}\mapsto \theta_{\gamma^i_{ij}}(p_{ij})=p_{ij}\sw{0}(\gamma^i_{ij})^{-1}(p_{ij}\sw{1})\ot p_{ij}\sw{2}.
\end{equation}

\begin{mlem}(Cf. the proof of Proposition~2 \cite{CalMat:Con}.)
Suppose that $(P(B)^C,(J_i)_{i\in I})$ is a \splce. Then, for all $i,j\in I$, $K^i_j=K^j_i$.
\end{mlem}
\begin{proof} (Cf. the proof of Proposition~2 \cite{CalMat:Con}.)
For all $i,j\in I$, define the map
\begin{equation}
\tilde{\phi}_{ji}:B_i\ot C\rightarrow P_{ij}^{\co C}\ot  C,\ \ 
b_i\ot c\mapsto\beta^j_{ij}\circ\chi^i_j(b_i\gamma_i(c)).
\end{equation}
It is easy to see that $\ker\tilde{\phi}_{ji}=\bar{K}^i_j\ot C$. Consider the maps
\begin{equation}
Q_{ji}:B_{ij}\rightarrow P_{ij}^{\co C},\ \ 
\pi^i_j(b_i)\mapsto (P_{ij}^{\co C}\ot \counit)\circ
\tilde{\phi}_{ji}(b_i1\sw{0}\gamma_i^{-1}(1\sw{1})\ot 1\sw{2}).
\end{equation}
Note that, as $\pi_i(\ker\pi_j)\subseteq \bar{K}^i_j$, and $\bar{K}^i_j$ is an ideal, 
for all $i,j\in I$,
the maps $Q_{ji}$ are well defined. Suppose that $b_{ij}\in K^i_j$. There exists an
element $b_i\in\bar{K}^i_j$ such that $\pi^i_j(b_i)=b_{ij}$.
It follows that $b_i1\sw{0}\gamma_i^{-1}(1\sw{1})\ot 1\sw{2}\in\bar{K}^i_j\ot C=\ker\tilde{\phi}_{ji}$,
hence $Q_{ji}(b_{ij})=0$, and so $K^i_j\subseteq \ker Q_{ji}$, for all $i,j\in I$.
On the other hand, for all $i,j\in I$ and $b\in B$,
\begin{multline*}
Q_{ji}(\pi_{ij}(b))(P_{ij}^{\co C}\ot \counit)\circ\tilde{\phi}_{ji}(\pi_i(b)1\sw{0}\gamma_i^{-1}(1\sw{1})\ot 1\sw{2})\\
=(P_{ij}^{\co C}\ot \counit)\circ\beta^j_{ij}(\chi_{ij}(b))
=(P_{ij}^{\co C}\ot \counit)(\chi_{ij}(b\sw{0})(\gamma^j_{ij})^{-1}(b\sw{1})\ot b\sw{2})\\
=\chi^j_i(\pi_j(b)1\sw{0}\gamma_j^{-1}(1\sw{1}))
=\pi_{ij}(b)\pi^j_i(1\sw{0}\gamma_j^{-1}(1\sw{1}))+K^j_i.
\end{multline*}
Suppose that, for some element $b_{ij}\in B_{ij}$, $Q_{ji}(b_{ij})=0$, i.e., 
$b_{ij}\pi^i_j(1\sw{0}\gamma_j^{-1}(1\sw{1}))\in K^j_i$. But the element 
$1\sw{0}\gamma_j^{-1}(1\sw{1})$,
by assumption, has a right inverse in $B_j$, and $K^j_i$ is an ideal in $B_{ij}$, hence 
$b_{ij}\in K^j_i$. It follows that $\ker(Q_{ij})=K^j_i$. 

We have proven, that, for all $i,j\in I$, $K^i_j\subseteq K^j_i$, and therefore,
for all $i,j\in I$, $K^i_j=K^j_i$. 
\end{proof}

Let $(P(B)^C,(J_i)_{i\in I})$ be a \splce{} such that the coalgebra $C$ is flat as a $\mK$-module.
Recall from the discussion following equation (\ref{covmap}) that
$P\cp$ is naturally an algebra and right $C$-comodule, and the map 
$\kappa_{P}:P\rightarrow P\cp$, $p\mapsto(\chi_i(p))_{i\in I}$ is an algebra and a right
$C$-comodule isomorphism. It follows that $\kappa_P(B)=(P\cp)^{\co C}$.
Is is clear that
\begin{equation}\label{eqfirstl}
\kappa_P(B)\subseteq\breve{B}=\{(b_i)_{i\in I}\in\bigoplus_{i\in I}B_i\ |\ 
\forall_{i,j\in I}\; \chi^i_j(b_i)=\chi^j_i(b_j)\}.
\end{equation}
On the other hand, let $(b_i)_{i\in I}\in\breve{B}$, then, for all $(p_i)_{i\in I}\in P\cp$,
\begin{equation*}
\coact^C((b_i)_{i\in I}(p_j)_{j\in I})=(\coact^C(b_ip_i))_{i\in I}
=(b_i\coact^C(p_i))_{i\in I}=(b_i)_{i\in I}\coact^C((p_j)_{j\in I}),
\end{equation*}
i.e., $(b_i)_{i\in I}\in(P\cp)^{\co C}$. It follows that $\breve{B}=(P\cp)^{\co C}=\kappa_P(B)$.

Suppose that $(B\cap J_i)_{i\in I}$ is a complete covering of $B$. Then
\begin{equation}\label{eqsecondl}
\kappa_P(B)=B\cp=\{(b_i)_{i\in I}\in\bigoplus_{i\in I}B_i\ |\ 
\forall_{i,j\in I}\; \pi^i_j(b_i)=\pi^j_i(b_j)\}.
\end{equation}
Let $\mu_{ij}:B_{ij}\rightarrow B_{ij}/K^i_j=P_{ij}^{\co C}$ be the canonical surjections. Observe that
$\left.\chi^i_j\right|_{B_i}=\mu_{ij}\circ\pi^i_j$, for all $i,j\in I$. By (\ref{eqfirstl}) and (\ref{eqsecondl}),
for all $(b_i)_{i\in I}\in\bigoplus_{i\in I}B_i$, the condition
\begin{equation}
\pi^i_j(b_i)=\pi^j_i(b_j), \text{ for all }i,j\in I,
\end{equation}
is equivalent to 
\begin{equation}\mu_{ij}(\pi^i_j(b_i)-\pi^j_i(b_j))=0,\text{ for all }i,j\in I.
\end{equation}
In particular, we have the following
\begin{mprop}(cf. Proposition~2, \cite{CalMat:Glu}.)\label{semiispromp}
Let $(P(B)^C,(J_i)_{i\in I})$ be a \splce{} such that the coalgebra $C$ is flat as a $\mK$-module, 
and $I=\{1,2\}$. Then $(P(B)^C,(J_i)_{i\in I})$ is a proper locally cleft extension.
\end{mprop}
\begin{proof}
We prove by contradiction. Suppose that $K^1_2\neq \{0\}$. Then there exists 
an element $r\in \bar{K}^1_2$ such
that $\pi^1_2(r)\neq 0$. Let $(b_1,b_2)\in B\cp$, then $(b_1+r,b_2)\in \breve{B}$ and
$(b_1+r,b_2)\notin B\cp$ which, by the discussion preceding the above proposition, is impossible.
\end{proof}

Suppose that $(P(B)^C,(J_i)_{i\in I})$ is a proper and regular 
locally cleft extension. Let us define the 
family of gauge
transformations
\begin{equation}
\Xi_{ij}:C\rightarrow B_{ij},\ c\mapsto \gamma^i_{ij}(c\sw{1})(\gamma^j_{ij})^{-1}(c\sw{2}),\ \text{ for all }
 i,j\in I. \label{glocdef}
\end{equation}
The gauge transformations $\Xi_{ij}$, $i,j\in I$ satisfy the following conditions. For all $i,j,k\in I$, $c\in C$,
\begin{gather}
\Xi_{ii}(c)=\counit(c),\ \ \ 
\Xi_{ji}=(\Xi_{ij})^{-1},\\
\pi^{ij}_k(\Xi_{ij}(c))\pi^{ik}_j(\Xi_{ik}(c\sw{1}))\pi^{kj}_i(\Xi_{kj}(c\sw{2})). \label{cocyclegauge}
\end{gather}
The first two of the above equalities are obvious, to prove the last one observe that, using 
(\ref{idenchiijk}), we obtain, for all $i,j,k\in I$,
\begin{multline*}
\pi^{ik}_j(\Xi_{ik}(c\sw{1}))\pi^{kj}_i(\Xi_{kj}(c\sw{2}))
=\chi^{ik}_j(\gamma^i_{ik}(c\sw{1})(\gamma^k_{ik})^{-1}(c\sw{2}))
\chi^{kj}_i(\gamma^k_{kj}(c\sw{3})(\gamma^j_{kj})^{-1}(c\sw{4}))\\
=\chi^i_{kj}(\gamma_i(c\sw{1}))\chi^k_{ij}(\gamma_k^{-1}(c\sw{2}))
\chi^k_{ij}(\gamma_k(c\sw{3}))\chi^j_{ik}(\gamma_j^{-1}(c\sw{4}))\\
=\chi^i_{kj}(\gamma_i(c\sw{1}))\chi^j_{ik}(\gamma_j^{-1}(c\sw{2}))
=\chi^{ij}_k(\chi^i_j(\gamma_i(c\sw{1}))\chi^j_i(\gamma_j^{-1}(c\sw{2}))\\
=\chi^{ij}_k(\gamma^i_{ij}(c\sw{1})(\gamma^j_{ij})^{-1}(c\sw{2}))
=\pi^{ij}_k(\Xi_{ij}(c)).
\end{multline*}

Suppose that the ground ring $\mK$ is a field and $I=\{1,2\}$. Let $(P(B)^C,(J_i)_{i\in I})$ be a 
proper locally cleft $C$-coalgebra Galois extension. By Corollary~\ref{corfieldcov},
$P(B)^C$ is a $C$-coalgebra Galois extension if and only if
\begin{equation}
\ker(\chi_1\ot_B\chi_1)\cap\ker(\chi_2\ot_B\chi_2)=\{0\}.
\end{equation}
We shall devote the remainder of the present section to proving that, under not very restrictive
assumptions, the above condition is always satisfied.

\begin{mlem}\label{fibprobaslem}
Let $M_1$, $M_2$, $M_{12}$ be $\mK$-vector spaces, let $\pi^1_2:M_1\rightarrow M_{12}$,
$\pi^2_1:M_2\rightarrow M_{12}$ be surjective linear morphisms. Let 
$M=\{(m,n)\in M_1\oplus M_2|\pi^1_2(m)=\pi^2_1(n)\}$, and denote the projections
on the summands of the direct sum by
$\pi_1:M\rightarrow M_1$, $(m,n)\mapsto m$ and
$\pi_2:M\rightarrow M_2$, $(m,n)\mapsto n$.
As $\mK$ is a field, $\ker\pi^1_2$ and $\ker\pi^2_1$ are direct summands in $M_1$ and $M_2$
respectively, i.e., $M_1=\overline{M}_1\oplus\ker(\pi^1_2)$, $M_2=\overline{M}_2\oplus \ker(\pi^2_1)$,
for some subspaces $\overline{M}_i\subseteq M_i$, $i=1,2$. Let 
\begin{center}
\begin{tabular}{ll}
$\{m_i\}$ be a basis of $\ker\pi^1_2$, & $\{n_i\}$ be a basis of $\ker\pi^2_1$,\\
$\{\bar{m}_i\}$ be a basis of $\overline{M}_1$, & $\{\bar{n_i}\}$ be a basis of $\overline{M}_2$.
\end{tabular}
\end{center}
Suppose that $f:\overline{M}_1\rightarrow M_2$ is a linear map such
that, for all $m\in \overline{M}_1$, $\pi^1_2(m)=\pi^2_1(f(m))$. Then
\begin{flushleft}
the family $\{(0,n_i)\}$ is a basis of $\ker\pi_1$,\\
the family $\{(m_i,0)\}$ is a basis of $\ker\pi_2$,\\
the family $\{(\bar{m}_i,f(\bar{m}_i))\}$ is linearly independent.
\end{flushleft} 
Moreover, denote $\overline{M}=\Span(\{(\bar{m}_i,f(\bar{m}_i))\})$. Then 
$M=\overline{M}\oplus\ker(\pi_1)\oplus\ker(\pi_2)$. 
\end{mlem}
{\bf Remark.} Observe that the map $f$ in Lemma~\ref{fibprobaslem} 
is necessarily injective. Note furthermore that the restriction
$\bar{\pi}^1_2:\overline{M}_1\rightarrow M_{12}$ 
(resp. $\bar{\pi}^2_1:\overline{M}_2\rightarrow M_{12}$) of the map $\pi^1_2$ (resp. $\pi^2_1$) is
an isomorphism, and, in particular, the map 
$f=(\bar{\pi}^2_1)^{-1}\circ\bar{\pi}^1_2:\overline{M}_1\rightarrow M_2$ satisfies the assumptions of
the above lemma. Finally, the restriction $\bar{\pi}_{12}:\overline{M}\rightarrow M_{12}$ of 
the map $\pi^1_2\circ\pi_1$ is an obvious linear isomorphism.
\begin{proof}
Suppose that $x\in \overline{M}\cap(\ker(\pi_1)+\ker(\pi_2))$. Then $x=(m,0)+(0,n)=(m,n)$,
for some elements $m\in \ker(\pi^1_2)$, $n\in \ker(\pi^2_1)$. On the other hand, 
$x\in\Span(\{(\bar{m}_i,f(\bar{m}_i))\})$, hence $x=\sum_i\alpha_i(\bar{m}_i,f(\bar{m}_i))
=(\sum_i\alpha_i\bar{m}_i, f(\sum_i\alpha_i\bar{m}_i))$, for some coefficients $\alpha_i\in\mK$.
 Consequently, 
$\sum_i\alpha_i\bar{m}_i\in\ker(\pi^1_2)$. As $\Span(\{\bar{m}_i\})\cap\ker(\pi^1_2)=\{0\}$, we have $x=0$.

Suppose that $x\in M$. Then there exist unique coefficients $\alpha_i,\beta_i,\gamma_i,\delta_i\in\mK$,
such that 
\begin{multline*}
x=\left(\sum_i\alpha_i\bar{m}_i+\sum_j\beta_jm_j,\sum_k\gamma_k\bar{n}_k+\sum_l\delta_ln_l\right)\\
=\sum_{i}\alpha_i(\bar{m}_i, f(\bar{m}_i))+\left(0,\sum_k\gamma_k\bar{n}_k-
f\left(\sum_i\alpha_i\bar{m}_i\right)\right)+\sum_j\beta_j(m_j,0)+\sum_l\delta_l(0,n_l).
\end{multline*}
Observe that 
\begin{equation*}\pi^2_1\left(\sum_k\gamma_k\bar{n}_k-
f\left(\sum_i\alpha_i\bar{m}_i\right)\right)=\pi^2_1\left(\sum_k\gamma_k\bar{n}_k\right)-
\pi^1_2\left(\sum_i\alpha_i\bar{m}_i\right)=0.
\end{equation*}
 Therefore there exist coefficients $\xi_s\in \mK$, such that
$\sum_k\gamma_k\bar{n}_k-
f(\sum_i\alpha_i\bar{m}_i)=\sum_s\xi_sn_s$. It follows that $M$ is spanned
by  vectors of the form
$(m_i,0)$, $(0,n_i)$, $(\bar{m}_i,f(\bar{m}_i))$. Their linear independence is obvious.
\end{proof}

\begin{mlem}\label{clbasformlem}
Assume that the ground ring $\mK$ is a field, and suppose that $P(B)^C_\gamma$ is a cleft
$C$-coalgebra Galois extension. Let $\{b_i\}$ be a linear basis of $B$ and let 
$\{h_i\}$ be a  linear basis of $C$. Then
\begin{equation}\label{pbpclbas}
\{b_i\gamma(h_j)\ot_B \gamma(h_k)\}
\end{equation}
is a linear basis of $P\ot_B P$.
\end{mlem}
\begin{proof}
 The map
 \begin{equation}F:P\ot_BP\rightarrow P\ot C,\ \ 
 p\ot_B p'\mapsto pp'\sw{0}\gamma^{-1}(p'\sw{1})\ot p'\sw{2},
 \end{equation}
  is a linear isomorphism (cf. Proposition~\ref{brzcl}), such that, for all $i,j,k$,
  \begin{equation}
  F(b_i\gamma(h_j)\ot_B \gamma(h_k))=b_i\gamma(h_j)\ot h_k.
  \end{equation}
  Vectors $\{b_i\gamma(h_j)\ot h_k\}$ form a basis of $P\ot C$.
  We conclude that the family (\ref{pbpclbas}) is a basis of $P\ot_B P$. 
\end{proof}

\begin{mprop}\label{glutwogallem}
Suppose that the ground ring $\mK$ is a field and that $I=\{1,2\}$.
Let $(P(B)^C,(J_i)_{i\in I})$ be a regular locally cleft $C$-coalgebra Galois extension.
Assume that
\begin{equation}\label{comuclcond} 
\gamma_2(C)\ker(\pi^2_1)=\ker(\pi^2_1)\gamma_2(C).
\end{equation}
Then
$P(B)^C$ is a $C$-coalgebra Galois extension.
\end{mprop}
{\bf Remark.} Note that if $C$ is a Hopf algebra and $P_2=B_2\ot C$, 
then one can choose $\gamma_2:C\rightarrow P_2$, $c\mapsto 1\ot c$, and then the
condition (\ref{comuclcond}) is automatically satisfied.
\begin{proof}
Suppose that $B_1=\overline{B}_1\oplus \ker(\pi^1_2)$, $B_2=\overline{B}_2\oplus \ker(\pi^2_1)$
and 
\begin{center}
\begin{tabular}{lll}
$\{\bar{x}_i\}$ is a basis of $\overline{B}_1$, & $\{\bar{y}_i\}$ is a basis of
 $\overline{B}_2$, & $\{h_i\}$ is a basis of $C$,\\
 $\{x_i\}$ is a basis of $\ker\pi^1_2$, & $\{y_i\}$ is a basis of $\ker\pi^2_1$. & { }
\end{tabular}
\end{center}By Proposition~\ref{brzcl}, 
\begin{equation*}
P_1=\overline{B}_1\gamma_1(C)\oplus\ker(\pi^1_2)\gamma_1(C),\ \ 
P_2=\overline{B}_2\gamma_2(C)\oplus\ker(\pi^2_1)\gamma_2(C).
\end{equation*}
Let $\overline{P}_1=\overline{B}_1\gamma_1(C)$, 
$\overline{P}_2=\overline{B}_2\gamma_2(C)$. 
By Proposition~\ref{semiispromp}, $(P(B)^C,(J_i)_{i\in I})$ is a proper locally cleft extension.
Therefore $\ker\chi^1_2=\ker(\pi^1_2)\gamma_1(C)$ and $\ker\chi^2_1=\ker(\pi^2_1)\gamma_2(C)$.
It follows that
\begin{equation}
P_1=\overline{P}_1\oplus\ker(\chi^1_2),\ \ 
P_2=\overline{P}_2\oplus\ker(\chi^2_1),
\end{equation}
and then,
\begin{equation}
\begin{array}{ll}
\{\bar{x}_i\gamma_1(h_j)\} \text{ is a basis of } \overline{P}_1, &
 \{\bar{y}_i\gamma_2(h_j)\} \text{ is a basis of } \overline{P}_2 \\
 \{x_i\gamma_1(h_j)\} \text{ is a basis of } \ker(\chi^1_2), &
 \{y_i\gamma_2(h_j)\}\text{ is a basis of } \ker(\chi^2_1).
\end{array}
\end{equation}
For all $i$, let $\breve{x}_i$ denote the unique element of $\overline{B}_2$ such that
$\pi^1_2(\bar{x}_i)=\pi^2_1(\breve{x}_i)$.
Similarily, for all $i$, let $\breve{h}_i$ (resp. $\breve{g}_i$) 
be the unique element of $\overline{P}_2$ (resp $\overline{P}_1$) with the property
$\chi^1_2(\gamma_1(h_i))=\chi^2_1(\breve{h}_i)$
(resp. $\chi^2_1(\gamma_2(h_i))=\chi^1_2(\breve{g}_i)$).  
Let us define the following elements of $P$:
\begin{equation}
\begin{array}{lll}
X_i=\kappa_P^{-1}((x_i,0)), & 
Y_j=\kappa_P^{-1}((0,y_j)), &
\overline{X}_k=\kappa_P^{-1}((\bar{x}_k,\breve{x}_k)),\\
H_s=\kappa_P^{-1}((\gamma_1(h_s),\breve{h}_s)),&
G_t=\kappa_P^{-1}((\breve{g}_t,\gamma_2(h_t))), & { } 
\end{array}
\end{equation} 
for all $i,j,k,s,t$.
Note that condition (\ref{comuclcond}) implies that
\begin{equation}\label{yggy}
\Span(\{Y_jG_t\})=\Span(\{G_tY_j\}).
\end{equation}
Let us define $\overline{B}=\Span(\{\overline{X}_k\})$,
$\overline{P}=\Span(\{\overline{X}_kH_s\})$.
It follows immediately from Lemma~\ref{fibprobaslem} that
\begin{equation}
B=\overline{B}\oplus\ker(\pi_1)\oplus\ker(\pi_2),\ \ 
P=\overline{P}\oplus\ker(\chi_1)\oplus\ker(\chi_2),
\end{equation}
and 
\begin{equation}
\begin{array}{llllll}
\{X_i\} & \text{is a basis of} & \ker(\pi_2),\ \ & 
\{X_iH_s\} & \text{is a basis of} & \ker(\chi_2),\\ 
\{Y_j\} & \text{is a basis of} & \ker(\pi_1),\ \ &
\{Y_jG_t\} & \text{is a basis of} & \ker(\chi_1),\\
\{\overline{X}_k\} & \text{is a basis of} & \overline{B},\ \ &
\{\overline{X}_kH_s\} & \text{is a basis of } & \overline{P}.
\end{array}\label{basofp}
\end{equation}
Denote $\overline{P\ot_BP}=\Span(\{\overline{X}_kH_s\ot_B H_t\})$.
We claim that
\begin{equation}
P\ot_BP=\overline{P\ot_BP}\oplus\ker(\chi_1\ot_B\chi_1)\oplus
\ker(\chi_2\ot_B\chi_2),
\end{equation}
and
\begin{equation}\label{basofpbp}
\begin{array}{lll}
\{\overline{X}_kH_s\ot_B H_t\} & \text{is a basis of} & \overline{P\ot_BP},\\
\{X_iH_s\ot_B H_t\} & \text{is a basis of} & \ker(\chi_2\ot_B\chi_2),\\
\{Y_jG_s\ot_B G_t\} & \text{is a basis of} & \ker(\chi_1\ot_B\chi_1).
\end{array}
\end{equation}
First we prove that the above vectors span $P\ot_BP$. Denote
$(P\ot_BP)'=($Span of the vectors (\ref{basofpbp})$)$. By statements (\ref{basofp}), it is obvious that
the vectors
\begin{center}
\begin{tabular}{||r|l||r|l||r|l||}
\hline
{\bf 1)} & $X_iH_s\ot_B X_jH_t$,  &
{\bf 2)} & $X_iH_s\ot_B Y_mG_t$, &
{\bf 3)} & $X_iH_s\ot_B \overline{X}_kH_t$,\\
\hline

{\bf 4)} & $Y_mG_s\ot_B X_iH_t$, &
{\bf 5)} & $Y_mG_s\ot_B Y_nG_t$, &
{\bf 6)} & $Y_mG_s\ot_B\overline{X}_kH_t$,\\ \hline
{\bf 7)} & $\overline{X}_kH_s\ot_B X_iH_t$, &
{\bf 8)} & $\overline{X}_kH_s\ot_B Y_nG_t$, &
{\bf 9)} & $\overline{X}_kH_s\ot_B\overline{X}_lH_t$,\\
\hline 
\end{tabular} 
\end{center}
for all $i,j,k,l,m,n,s,t$,
span $P\ot_BP$.
Moving $B$-factors from  the right to the left  leg in each of the above tensor products, and
then using the list (\ref{basofp}), we obtain
the following results.  

\noindent{\bf The tensor products of types 2) and 4)} are simply equal to zero, as
\[\ker(\chi_1)\ker(\pi_2)=\ker(\chi_2)\ker(\pi_1)=\{0\}.\]

\noindent {\bf The tensor products of types 1), 3), 7)} clearly belong to
\begin{equation*}
\ker(\chi_2)\ot_B\Span(\{H_t\})\subseteq\Span(\{X_iH_s\ot_B H_t\})
\subseteq (P\ot_BP)'.
\end{equation*}

\noindent {\bf The tensor products of types 5) and 8)} clearly belong to
\begin{equation*}
\ker(\chi_1)\ot_B\Span(\{G_t\})\subseteq\Span(\{Y_jG_s\ot_B G_t\})
\subseteq (P\ot_BP)'.
\end{equation*}

\noindent {\bf The tensor products of type 6)} belong to
 \begin{multline*}
 \ker(\chi_1)\ot_B\Span(\{H_t\})
 =\Span(\{Y_jG_s\})\ot_B\Span(\{H_t\})\\
 =\Span(\{G_sY_j\})\ot_B\Span(\{H_t\})
 \subseteq\Span(\{G_s\})\ot_B\ker(\chi_1)\\
 =\Span(\{G_s\})\ot_B\Span(\{Y_jG_t\})\\
 \subseteq \Span(\{Y_jG_s\ot_B G_t\})\subseteq (P\ot_BP)',
 \end{multline*}
 where in the second equality we used eq.~(\ref{yggy}).
 
 \noindent {\bf The tensor products of type 9)} belong to
 \begin{equation*}
 \begin{split}
 P\ot_B\Span(\{H_t\})& \Span(\{\overline{X}_kH_s\ot_B H_t\})\oplus\Span(\{X_iH_s\ot_B H_t\})\\
 {}&\ \ \ \ \ \ \oplus\Span(\{Y_jG_s\ot_B H_t\}).
 \end{split}
 \end{equation*}
 In the previous point, we have proven that the last direct summand in the above expression 
 is also contained in $(P\ot_B P)'$. Therefore tensors of type 9) belong to $(P\ot_BP)'$.
 It follows that the tensors (\ref{basofpbp}) span $P\ot_BP$. 
 
 Suppose that the element $z\in \ker(\chi_1\ot_B\chi_1)\cap\ker(\chi_2\ot_B\chi_2)$.
 As tensors (\ref{basofpbp}) span $P\ot_BP$, there exists a family of coefficients
 $a_{ist},b_{jst},c_{kst}\in \mK$, 
 such that
 \begin{equation}
 z=\sum_{i,s,t}a_{ist}X_iH_s\ot_B H_t+\sum_{j,s,t}b_{jst}Y_jG_s\ot_B G_t
 +\sum_{k,s,t}c_{kst}\overline{X}_{k}H_s\ot_BH_t.
 \end{equation}
 Note that, for all $i,j,k,s$, 
 \begin{gather*}\chi_1(X_i)=x_i,\ \chi_1(Y_j)=0,\ \chi_1(\overline{X}_k)=\bar{x}_k,\ 
 \chi_1(H_s)=\gamma_1(h_s),\\ 
 \chi_2(Y_j)=y_j,\ \chi_2(G_s)=\gamma_2(h_s).
 \end{gather*}
 It follows that
 \begin{equation*}
0=(\chi_1\ot_B\chi_1)(z)=\sum_{i,s,t}a_{ist}x_i\gamma_1(h_s)\ot_B\gamma_1(h_t)
+\sum_{k,s,t}c_{kst}\bar{x}_k\gamma_1(h_s)\ot_B\gamma_1(h_t).
\end{equation*}
 By Lemma~\ref{clbasformlem}, this implies that $a_{ist}=c_{kst}=0$, for all $i,k,s,t$.
 Then
 \begin{equation*}
 0=(\chi_2\ot_B\chi_2)(z)=\sum_{j,s,t}b_{jst}y_j\gamma_2(h_s)\ot_B\gamma_2(h_t).
 \end{equation*}
 It follows, by Lemma~\ref{clbasformlem}, that $b_{jst}=0$, for all $j,s,t$, 
 i.e., $z=0$ and $\ker(\chi_1\ot_B\chi_1)\cap\ker(\chi_2\ot_B\chi_2)=\{0\}$. 
 By Corollary~\ref{corfieldcov}, $P(B)^C$ is a $C$-coalgebra Galois extension.
 Note that we have also proven that the tensors
 (\ref{basofpbp}) are linearly independent, and, consequently, they form a basis of $P\ot_B P$.  
\end{proof}

\section{Gluing cleft extensions}\label{gluclextsectm}

The quantum geometry situation, which corresponds to the most usual setting for the classical method of constructing 
principal bundles by patching together trivial principal bundles, is as follows.
We are given an algebra $B$, which has a complete covering $(K_i)_{i\in I}$, and a coalgebra $C$.
For each of the quotient spaces $B_i=B/K_i$, $i\in I$ (resp. $B_{ij}=B/(K_i+K_j)$, $i,j\in I$), we 
construct a cleft $C$-coalgebra Galois extension $P_i(B_i)^C_{\gamma_i}$
(resp. $P_{ij}(B_{ij})^C_{\gamma_{ij}}$). 
Let us denote by $\pi_i:B\rightarrow B_i$, $\pi^i_j:B_i\rightarrow B_{ij}$, $i,j\in I$, etc.,  the
canonical surjections. 
Then we choose surjective algebra and right $C$-comodule
morphisms 
\begin{equation}
\chi^i_j:P_i\rightarrow P_{ij},\ \  i,j\in I,
\label{mapschiij} 
\end{equation}
such that
$\chi^i_j(b_i)=\pi^i_j(b_i)$, for all $i,j\in I$, $b_i\in B_i$, and we use them for gluing
(cf. \ref{gluingdef})
\begin{equation}\label{princgludef}
P=\bigoplus_{\chi^i_j}P_i\{(p_i)_{i\in I}\in\bigoplus_{i\in I}P_i\ |\ \forall_{i,j\in I}\;\chi^i_j(p_i)=\chi^j_i(p_j)\}.
\end{equation}
If the coalgebra $C$ is flat as a $\mK$-module then (cf. the discussion in Section~\ref{glusect})
$P$ is naturally a right $C$-comodule. Then,
for each $n\in I$, we define the algebra and right $C$-comodule map 
\begin{equation}
\chi_n:P\rightarrow P_n,\ \ (p_i)_{i\in I}\mapsto p_n. \label{mapschin}
\end{equation}
If all the maps $\chi_i$, $i\in I$, are surjective then
$(P(B)^C,(\ker\chi_i)_{i\in I})$ is
 a proper locally cleft $C$-coalgebra Galois  extension. Moreover,
 for all $i\in I$, $P_i\simeq P/\ker(\chi_i)$. The following lemma gives necessary and sufficient
 conditions for the maps $\chi_i$, $i\in I$, to be surjective.
 
 \begin{mlem}\label{priglucllem}
 Suppose that, for each $i\in I$, the element $1\sw{0}\gamma_i^{-1}(1\sw{1})\in B_i$ has a 
 right inverse in $B_i$. 
 The algebra and right $C$-comodule maps (\ref{mapschiij}) satisfy
 the condition $\chi^i_j(b_i)=\pi^i_j(b_i)$, for all $i,j\in I$, $b_i\in B_i$, if and only if there
  exists a family of convolution invertible maps (gauge transformations) 
  $\Gamma^i_j:C\rightarrow B_{ij}$, $i,j\in I$, such that
  \begin{equation}
  \chi^i_j(p_i)=\pi^i_j(p_i\sw{0}\gamma_i^{-1}(p_i\sw{1}))\Gamma^i_j(p_i\sw{2})\gamma_{ij}(p_i\sw{3}),
  \label{chiijgammaij}
  \end{equation}
  for all $i,j\in I$, $p_i\in P_i$. Furthermore, assume that the coalgebra 
  $C$ is flat as a $\mK$-module. Let
   $I=\{1,2,\ldots,N\}$, and
   suppose that either $N\leq 3$, or the algebra $B$ and its complete covering
  $(K_i)_{i\in I}$ satisfy the condition (cf. eq.~(\ref{covdeepcap})),
  \begin{equation}
\bigcap_{1\leq j\leq i}(\ker\pi^{k+1}_j+\ker\pi^{k+1}_{i+1})
=\left(\bigcap_{1\leq j\leq i}\ker\pi^{k+1}_j\right)+\ker\pi^{k+1}_{i+1},\label{kerchiijkcap}
\end{equation} 
for all
$1\leq k<N$ and $1\leq i<k$. Define gauge transformations $\Xi_{ij}:C\rightarrow B_{ij}$, 
$c\mapsto\Gamma^i_j(c\sw{1})(\Gamma^j_i)^{-1}(c\sw{2})$, $i,j\in I$ (cf. eq.~(\ref{glocdef})). 
Then
the maps (\ref{mapschin}) are surjective if and only if the condition (\ref{cocyclegauge}) is satisfied:
\begin{equation}
\pi^{ij}_k(\Xi_{ij}(c))\pi^{ik}_j(\Xi_{ik}(c\sw{1}))\pi^{kj}_i(\Xi_{kj}(c\sw{2})). \label{cocxicon}
\end{equation}  
 \end{mlem}
 \noindent {\bf Remark.} It is clear that, while maps $\Gamma^i_j$, $i,j\in I$, define surjections
 $\chi^i_j$, the space $P=\bigoplus_{\chi^i_j}P_i$ is fully defined by the maps $\Xi_{ij}$, $i,j\in I$. Indeed,
 for all $i,j\in I$, $p_i\in P_i$, $p_j\in P_j$, the condition $\chi^i_j(p_i)=\chi^j_i(p_j)$ is equivalent
 to
 \begin{equation}
 \pi^i_j(p_i\sw{0}\gamma_i^{-1}(p_i\sw{1}))\Xi_{ij}(p_i\sw{2})\ot p_i\sw{3}
 =\pi^j_i(p_j\sw{0}\gamma_j^{-1}(p_j\sw{1}))\ot p_j\sw{2}.
 \end{equation}
 \begin{proof}
 Suppose that the algebra and right $C$-comodule maps (\ref{mapschiij}) satisfy
 the condition $\chi^i_j(b_i)=\pi^i_j(b_i)$, for all $i,j\in I$, $b_i\in B_i$. Then,
 for all $p_i\in P_i$, 
 \begin{multline*}
 \chi^i_j(p_i)=\chi^i_j(p_i\sw{0}\gamma_i^{-1}(p_i\sw{1})\gamma_i(p_i\sw{2}))\\
 =\pi^i_j(p_i\sw{0}\gamma_i^{-1}(p_i\sw{1}))\chi^i_j(\gamma_i(p_i\sw{2}))\gamma_{ij}^{-1}(p_i\sw{3})
 \gamma_{ij}(p_i\sw{4})).
 \end{multline*}
 Defining the maps 
 $\Gamma^i_j(c)=\chi^i_j(\gamma_i(c\sw{1}))\gamma_{ij}^{-1}(c\sw{2})$, $i,j\in I$ yields
 eq.~(\ref{chiijgammaij}).
 Conversely, let the maps (\ref{mapschiij}) have the form 
 (\ref{chiijgammaij}). Then, for all $i,j\in I$, $b_i\in B_i$,
 \begin{equation*}
 \chi^i_j(b_i) \pi^i_j(b_i)\pi^i_j(1\sw{0}\gamma_i^{-1}(1\sw{1}))\Gamma^i_j(1\sw{2})\gamma_{ij}(1\sw{3})
 =\pi^i_j(b_i)\chi^i_j(1)=\pi^i_j(b_i). 
 \end{equation*}
 
 Assume that the coalgebra $C$ is flat as a $\mK$-module.
 We will check that the maps (\ref{mapschiij}) satisfy the assumptions of 
 Proposition~\ref{chinsurpro}, which in turn will prove that the maps
 (\ref{mapschin}) are surjective.
 
 Note that 
 \begin{equation}
 \ker\chi^i_j=\ker(\pi^i_j)\gamma_i(C), \text{ for all }i,j\in I.\label{kerchiij}
 \end{equation}
  Indeed, by (\ref{clthe}), 
 $\chi^i_j=\theta^{-1}_{\gamma^i_{ij}}\circ(\pi^i_j\ot C)\circ\theta_{\gamma_i}$,
 which means that $\ker(\chi^i_j)=\theta^{-1}_{\gamma_i}(\ker(\pi^i_j\ot C))
 =\theta^{-1}_{\gamma_i}(\ker(\pi^i_j)\ot C)=\ker(\pi^i_j)\gamma_i(C)$,
 where $\gamma^i_{ij}=\chi^i_j\circ\gamma_i$, and the second equality follows from the flatness of $C$.
 Observe that (condition~(\ref{calmat1as}) of Proposition~\ref{chinsurpro})
 \begin{equation}
 \chi^i_j(\ker\chi^i_k)=\chi^j_i(\ker\chi^j_k), \text{ for all }i,j,k\in I.
 \end{equation}
 Indeed, for all $i,j,k\in I$, $\pi^i_j(\ker\pi^i_k)=\pi^j_i(\ker\pi^j_k)$, and then, as $\mK$-modules
  $\pi^i_j(\ker\pi^i_k)$
 are ideals, for all $c\in C$, $\pi^i_j(\ker\pi^i_k)\Xi_{ij}(c\sw{1})\ot c\sw{2}\subseteq
 \pi^i_j(\ker\pi^i_k)\ot C=\pi^j_i(\ker\pi^j_k)\ot C$. Consequently,
 $\pi^i_j(\ker\pi^i_k)\gamma^i_{ij}(C)\subseteq \pi^j_i(\ker\pi^j_k)\gamma^j_{ij}(C)$. Furthermore
  by (\ref{kerchiij}), $\chi^i_j(\ker\chi^i_k)=\pi^i_j(\ker\pi^i_k)\gamma^i_{ij}(C)$,
  hence
 it follows that, for all $i,j,k\in I$, $\chi^i_j(\ker\chi^i_k)\subseteq\chi^j_i(\ker\chi^j_k)$.
 
 For each $i,j,k\in I$, the map
 \begin{gather}
 W^i_{jk}:P_i/(\ker(\chi^i_j)+\ker(\chi^i_k))\rightarrow B_{ijk}\ot C,\nonumber\\ 
 p_i+\ker(\chi^i_j)+\ker(\chi^i_k)\mapsto\pi^i_{jk}(p_i\sw{0}\gamma_i^{-1}(p_i\sw{1}))\ot p_i\sw{2},
 \end{gather} 
 is an isomorphism. Indeed, since $\ker(\chi^i_j)+\ker(\chi^i_k)=\ker(\pi^i_{jk})\gamma_i(C)$,
 $W^i_{jk}$ is well defined. It is also obviously surjective. Moreover, suppose that 
 $W^i_{jk}(p_i+\ker(\chi^i_j)+\ker(\chi^i_k))=0$. This means that
 $\theta_{\gamma_i}(p_i)\in\ker(\pi^i_{jk}\ot C)=\ker(\pi^i_{jk})\ot C$, hence 
 $p_i\in \ker(\chi^i_j)+\ker(\chi^i_k)$. Note that, for all
 $b_i\in B_i$ and $c\in C$, 
 \[
 (W^i_{jk})^{-1}(\pi^i_{jk}(b_i)\ot C)=b_i\gamma_i(c)+\ker(\chi^i_j)+\ker(\chi^i_k).
 \]
 
 Suppose that the maps 
 \begin{equation*}
 \phi^k_{ij}:P_j/(\ker\chi^j_i+\ker\chi^j_k)\rightarrow 
P_i/(\ker\chi^i_j+\ker\chi^i_k),\ \ i,j,k\in I
\end{equation*}
 are the isomorphisms (\ref{phiijkdef}), i.e.,
for all $p_j\in P_j$,
\begin{equation*}
\phi^k_{ij}(p_j+\ker(\chi^j_i)+\ker(\chi^j_k))=p_i+\ker(\chi^i_j)+\ker(\chi^i_k),
\end{equation*}
 where
$p_i$ is any element of $P_i$ such that $\chi^i_j(p_i)=\chi^j_i(p_j)$ (cf. Remark after
Proposition~\ref{chinsurpro}). For each $i,j,k\in I$, let us define the isomorphisms
\begin{equation}
\bar{\phi}^k_{ij}=W^i_{jk}\circ\phi^k_{ij}\circ(W^j_{ik})^{-1}:B_{ijk}\ot C\rightarrow B_{ijk}\ot C.
\end{equation}
It is easy to see that explicitly, for all $b_{ijk}\in B_{ijk}$, $c\in C$, 
\begin{equation}
\bar{\phi}^k_{ij}(b_{ijk}\ot c)=b_{ijk}\pi^{ij}_k(\Xi_{ji}(c\sw{1}))\ot c\sw{2}.\label{xiphirel}
\end{equation}
Clearly, the condition (\ref{phiijkcoc}) is equivalent to 
\begin{equation}
\bar{\phi}_{ik}^j=\bar{\phi}_{ij}^k\circ\bar{\phi}_{jk}^i, \text{ for all }i,j,k\in I,
\end{equation}
and this in turn is, by eq.~(\ref{xiphirel}), equivalent to the condition (\ref{cocxicon}).

Finally, in view of the flatness of $C$ and the condition (\ref{kerchiijkcap}) we obtain,
for all $1\leq k<N$ and $1\leq i<k$,
\begin{multline*}\bigcap_{1\leq j\leq i}(\ker\chi^{k+1}_j+\ker\chi^{k+1}_{i+1})
=\theta^{-1}_{\gamma_{k+1}}\left(
\bigcap_{1\leq j\leq i}\theta_{\gamma_{k+1}}(\ker\chi^{k+1}_j+\ker\chi^{k+1}_{i+1})
\right)\\
=\theta^{-1}_{\gamma_{k+1}}\left(
\bigcap_{1\leq j\leq i}(\ker(\pi^{k+1}_{j,i+1})\ot C)\right)
=\theta^{-1}_{\gamma_{k+1}}\left(\ker\left(
\bigoplus_{1\leq j\leq i}\pi^{k+1}_{j,i+1}\ot C\right)\right)\\
=\theta^{-1}_{\gamma_{k+1}}\left(\ker\left(
\bigoplus_{1\leq j\leq i}\pi^{k+1}_{j,i+1}\right)\ot C\right)
=\theta^{-1}_{\gamma_{k+1}}\left(
\bigcap_{1\leq j\leq i}(\ker\pi^{k+1}_j+\ker\pi^{k+1}_{i+1})\ot C\right)\\
=\theta^{-1}_{\gamma_{k+1}}\left(
\bigcap_{1\leq j\leq i}(\ker\pi^{k+1}_j)\ot C\right)+\ker\chi^{k+1}_{i+1}\\
=\theta^{-1}_{\gamma_{k+1}}\left(\ker\left(
\bigoplus_{1\leq j\leq i}\pi^{k+1}_j\ot C\right)\right)+\ker\chi^{k+1}_{i+1}
=\bigcap_{1\leq j\leq i}\ker(\chi^{k+1}_j)+\ker(\chi^{k+1}_{i+1}).
\end{multline*}
Thus all the assumptions of Proposition~\ref{chinsurpro} are satisfied, and hence the maps
(\ref{mapschin}) are surjective.
 \end{proof}

\section{Example: The quantum lens spaces}\label{lensexsectm}

It was shown in \cite{HaMaSz:Pro} that by gluing two quantum discs, $D_p$ and $D_q$, one can obtain
the quantum 2-sphere 
$S^2_{pq}$, 
and that the universal $C^\ast$ algebra of functions on $S_{pq}$ is 
isomorphic to equatorial or latitudinal Podle\'s spheres
 (\cite{Pod:Sph}). In \cite{CalMat:Con} 
a quantum sphere $S^3_{pq}$ was obtained, by gluing quantum solid tori 
(cf. Subsection~\ref{quantumsolidtorus}) $D_p\times S^1$ and $D_q\times S^1$, 
as an example of a locally trivial $U(1)$-quantum principal bundle with the base space
$S^2_{pq}$. It was also shown
that $\ffun(S^3_{pq})(\ffun(S^2_{pq}))^{\ffun(U(1))}$ is a principal Hopf-Galois extension.

As an illustration of methods described earlier in this chapter, we will construct a locally
cleft Hopf Galois extension of $\ffun(S^2_{pq})$ by gluing two quantum solid
tori (Subsection~\ref{quantumsolidtorus}) $D_p\times_\theta S^1$ and $D_q\times_\theta S^1$,
obtaining this way quantum lens spaces $L^{p,q,\theta}_\beta$ of charge $\beta$, for all $\beta\in\mZ$.
 As a special case,
for $\beta=1$, this gives the Heegaard quantum sphere (\cite{BaHaMaSz:He}).

Another example of a construction of   quantum lens spaces by gluing
two quantum solid tori of type $D\times_\theta S^1$  can be found in  \cite{MatTom:Lens}.

\subsection{The quantum unit disc}\label{quantumdisc}

A two-parameter family of quantum unit discs was defined in \cite{Kl:Disc}. Here we consider the one 
parameter subfamily studied in \cite{Kli:Disc}. We start with a coordinate $\ast$-algebra $\ffun(D_p)$ generated by 
$x$ and the relation
\begin{equation}
x^\ast x-pxx^\ast=1-p,\ \ 0<p<1.\label{disc}
\end{equation} 
The spectrum of $xx^\ast$ (in any $C^\ast$-algebra completion) is
\begin{equation}
\spec(xx^\ast)=\{1-p^n\ |\ n=0,1,2,\ldots\}\cup \{1\},
\end{equation}
i.e. $xx^\ast\leq 1$, where $\leq$ is understood as an order relation  between positive operators.
This justifies the name `unit disc'. Furthermore, this relation means also that $\|x\|=1$ in any $C^\ast$
completion of $\ffun(D_p)$
(Theorem 2.1.1 \cite{Mur:Cs}).

Observe that (\ref{disc}) has the following useful symmetry. Let $x_{-}$ be the generator of 
$\ffun(D_{p^{-1}})$, (where we consciously abuse the notation as $p^{-1}\geq 1$), i.e.,
 \begin{equation*}
x_{-}^\ast x_{-}-p^{-1}x_{-}x_{-}^\ast=1-p^{-1}.
\end{equation*} 
Then assignment $x\mapsto x_{-}^{\ast}$ can be extended to a $\ast$-algebra isomorphism
\begin{equation}
\kappa_p:\ffun(D_p)\rightarrow\ffun(D_{p^{-1}}).\label{kappa}
\end{equation}

The coordinate algebra $\ffun(D_p)$ can be completed to the $C^\ast$-algebra
$C(D_p)$ with the norm 
\begin{equation}
\|a\|=\sup_{\varrho}\|a\|_{\varrho},\ \ a\in\ffun(D_p),\label{uninorm}
\end{equation}
where the  supremum is taken over all 
bounded representations 
$\varrho:\ffun(D_p)\rightarrow \Bnd(\Hil{H}_\varrho)$ 
of $\ffun(D_p)$ and 
$\|\cdot\|_\varrho$ denotes the operator norm in the representation $\varrho$. The $C^\ast$ algebra
$C(D_p)$ is called a {\em universal enveloping algebra} of $\ffun(D_p)$. Note that the norm 
(\ref{uninorm}) is well defined because $\|x\|_\varrho=1$ for all $\varrho$.  

Irreducible bounded representations of $\ffun(D_p)$
are unitarily equivalent to one of the following representations.
\begin{enumerate}
\item For all $0\leq\phi<2\pi$, there is a one dimensional representation
$\varrho_\phi:\ffun(D_p)\rightarrow \mC$,
\begin{equation}
\varrho_\phi(x)=e^{i\phi},\ \ \varrho_\phi(x^\ast)=e^{-i\phi}.
\end{equation}
\item There is also an infinitely dimensional representation
$\varrho_\infty:\ffun(D_p)\rightarrow \Bnd(\Hil{H}_\infty)$, where $\Hil{H}_\infty$ is generated
by  orthonormal vectors $\Psi_n$, $n=0,1,2,\ldots$, and
\begin{equation}
\varrho_\infty(x)\Psi_n=\sqrt{1-p^{n+1}}\Psi_{n+1},\ \ 
\varrho_\infty(x^\ast)\Psi_{n+1}=\sqrt{1-p^{n+1}}\Psi_n,\ \ 
\varrho_\infty(x^\ast)\Psi_0=0.\label{discinf}
\end{equation}
\end{enumerate}
The infinite dimensional representation $\varrho_\infty$ 
is faithful. Nonfaithful representations correspond to subsets 
of quantum discs. In particular, one-dimensional representations are also characters, that is, they 
correspond to the classical points of the quantum space. In the case of the 
quantum disc, one-dimensional
representations describe  the classical unit circle.    

Let us adopt notational convention $x^{-n}\equiv (x^\ast)^n$. It follows from (\ref{disc}) that
$
(1-xx^\ast)x=px(1-xx^\ast)
$,
and, therefore,
\begin{equation}
(1-xx^\ast)^nx^m=p^{mn}x^m(1-xx^\ast)^n,\ \ 
n\in\mN,m\in\mZ.\label{1mxx}
\end{equation} 
Moreover, it is easy to prove by induction that, for all $n\geq 0$,
\begin{subequations}\label{xxsxsx}
\begin{gather}
(x^\ast)^nx^n=1+\sum_{m=1}^n
(-1)^mp^{nm-\frac{m(m-1)}{2}}\pchoose{n}{m}{p^{-1}}
(1-xx^\ast)^m,\label{xsx}\\
x^n(x^\ast)^n=1+\sum_{m=1}^n
(-1)^mp^{-nm+\frac{m(m+1)}{2}}\pchoose{n}{m}{p}
(1-xx^\ast)^m,\label{xxs}
\end{gather}
\end{subequations}
where we  used the standard $p$-deformed binomial coefficients
\begin{gather}
\pchoose{n}{m}{p}\frac{[n]_p!}{[m]_p![n-m]_p!},\nonumber\\
[n]_p!=[1]_p[2]_p\ldots[n-1]_p[n]_p \text{ for } n\in\mN,\ 
[0]_p!=1,\nonumber\\
[n]_p=1+p+p^2+\ldots +p^{n-2}+p^{n-1} \text{ for } n\in\mN, \ 
[0]_p=0.\label{pchoosedef}
\end{gather}
Note that equation (\ref{xxs}) follows from (\ref{xsx}) by the application of the isomorphism (\ref{kappa}). 

It is now obvious that, for all $n,m\in \mZ$,
\begin{equation}
x^nx^m=x^{n+m}(1+Q^p_{n;m}(1-xx^\ast)),\label{xmxn} 
\end{equation}
where $Q^p_{n;m}$ is a polynomial of degree at most $\min\{|m|,|n|\}$ and such that $Q^p_{n;m}(0)=0$.
For example, if $m\geq n\geq 0$, then
$x^{-n}x^{m}=(x^{-n}x^n)x^{m-n}$, and then use (\ref{xsx}) and (\ref{1mxx}).

As elements $x^nx^m$, $n,m\in\mZ$ obviously span $\ffun(D_p)$ as a vector space, (\ref{xmxn}) means
that also elements of the form 
\begin{equation}
x^n(1-xx^\ast)^m,\ \  n\in \mZ, \ m\in\mN_0,\label{discbase}
\end{equation}
span $\ffun(D_p)$. In fact, using the  infinite dimensional representation (\ref{discinf}) one can prove that
the above family forms a basis of $\ffun(D_p)$. Indeed, suppose that
\[\sum_{n=-\infty}^\infty\sum_{m=0}^\infty A_{nm}x^n(1-xx^\ast)^m=0,\] for some coefficients
$A_{mn}\in\mC$, only a finite number of which are different from zero. Therefore, for any $k\geq 0$,
\begin{multline*}
0=\varrho_\infty(\sum_{n=-\infty}^\infty\sum_{m=0}^\infty A_{nm}x^n(1-xx^\ast)^m)\Psi_k
=\sum_{n=-\infty}^\infty\sum_{m=0}^\infty p^{mk}A_{nm}\varrho_\infty(x^n)\Psi_k\\
=\sum_{m=0}^\infty p^{mk}(
\sum_{n=1}^k A_{-nm}\sqrt{1-p^k}\ldots\sqrt{1-p^{k-n+1}}\Psi_{k-n}
+A_{0m}\Psi_k\\
+\sum_{n=1}^\infty A_{nm}
\sqrt{1-p^{k+1}}\ldots\sqrt{1-p^{k+n}}\Psi_{k+n}.
)
\end{multline*}
As vectors $\Psi_k$ are linearly independent and $p<1$, 
\begin{equation*}
\sum_{m=0}^\infty p^{mk}A_{nm}=0 \text{ for  all } k\geq 0,\ -k\leq n,
\end{equation*}
i.e., for each $n\in \mZ$, the polynomial $\sum_{m=0}^\infty A_{nm}t^m$ has an infinite 
number of distinct roots
$p^{k}$, $k\geq 0$ if $n\geq 0$ or $k\geq-n$ otherwise. 
This is possible only if $A_{nm}=0$ for all $n\in \mZ$, $m\geq 0$. In addition,
we have proven that the representation (\ref{discinf}) is faithful.

\subsection{The quantum torus}\label{quantumtorus}

Quantum torus was defined in \cite{Rie:TV}. The coordinate algebra of the quantum torus $\ffun(T_\phi)$, 
$\phi\in [0,2\pi)$, 
is generated by unitary elements $U$, $V$ which satisfy the following commutation relations
\begin{equation}
UV=e^{i\phi}VU,\ \ 
UV^\ast=e^{-i\phi}V^\ast U.\label{torus}
\end{equation}
Obviously, the elements
\begin{equation}
V^nU^m,\ \ \ 
n,m\in\mZ,
\end{equation}
form a basis of $\ffun(T_\phi)$.
The coordinate algebra $\ffun(T_\phi)$ can be completed to the enveloping 
$C^\ast$-algebra $C(T_\phi)$ using 
representations of $\ffun(T_\phi)$. The representation theory of $\ffun(T_\phi)$ depends on whether
$\phi$ is a rational or irrational multiple of $2\pi$.

Suppose that $\phi=2\pi\frac{M}{N}$, where $M,N\in \mN$, $M<N$, and $M$ and $N$ 
are relatively prime.
Then $U^N$ and $V^N$ are central in $\ffun(T_\phi)$, and we can classify irreducible
representations of $\ffun(T_\phi)$ according to their eigenvalues. It turns out that
 irreducible representations of $\ffun(T_\phi)$ include the ones  
 isomorphic to one of the representations
$\varrho_{\alpha\beta}:\ffun(T_\phi)\rightarrow \Bnd(\Hil{H}^{\alpha\beta})$, 
$\alpha,\beta\in [0,2\pi)$, where
\newcommand{\vc}{\Psi^{\alpha\beta}}
\newcommand{\rp}{\varrho_{\alpha\beta}}
$\Hil{H}^{\alpha\beta}$ is spanned by orthonormal vectors $\vc_n$, $n\in\mZ_N$, and
\begin{equation}
\label{rational}
\rp(U^{\pm 1})\vc_n=e^{\pm i\frac{\alpha}{N}}e^{\pm 2\pi in\frac{M}{N}}\vc_n,\ \ \ \
\rp(V^{\pm 1})\vc_n=e^{\pm i\frac{\beta}{N}}\vc_{n\pm 1}.
\end{equation}

On the other hand, if $\phi$ is an irrational multiple of $2\pi$, then  irreducible representations 
include the representations
unitarily isomorphic to one of the
$\varrho_\alpha:\ffun(T_\phi)\rightarrow\Bnd(\Hil{H}^\alpha)$, 
$\alpha\in [0,2\pi)$, where $\Hil{H}^\alpha$ is
the closure of the linear span of the family of orthonormal vectors $\Psi^\alpha_n$, $n\in\mZ$, and
\begin{equation}
\varrho_\alpha(U^{\pm 1})\Psi^\alpha_n=e^{\pm i\alpha}e^{\pm in\phi}\Psi^\alpha_n,\ \ \ \ 
\varrho_\alpha(V^{\pm 1})\Psi^\alpha_n=\Psi^\alpha_{n\pm 1}.\label{irrational}
\end{equation}

Note, that the quantum torus is not a type I $C^\ast$-algebra, therefore
(cf. the discussion at the end of the chapter 3 in \cite{Arv:Cs}) 
its irreducible representations cannot be explicitly listed.

Representations (\ref{irrational}) are faithful. Representations (\ref{rational}) are not faithful. However,
choose two sequences $(\alpha_n)_{n\in \mN}$, $(\beta_n)_{n\in\mN}$, such that
$\alpha_n, \beta_m\in [0,2\pi)$, $m,n\in\mN$ and $\alpha_n\neq\alpha_m$,
$\beta_n\neq\beta_m$ if $n\neq m$. Then the representation
\begin{equation}
\varrho^{(\alpha_n)_{n\in\mN}(\beta_n)_{n\in\mN}}=\bigoplus_{m,n\in\mN}
\varrho^{\alpha_m\beta_n}:\ffun(T_\phi)\longrightarrow\bigoplus_{m,n\in\mN}
\Bnd(\Hil{H}^{\alpha_m\beta_n})
\end{equation}
is faithful. Indeed, suppose that for some coefficients $A_{mn}\in\mC$, $m,n\in\mZ$, a 
finitely many of which are different from \newcommand{\rmn}{\varrho^{(\alpha_n)_{n\in\mN}(\beta_n)_{n\in\mN}}}
zero, $\rmn(\sum_{m,n\in\mZ}A_{mn}V^mU^n)=0$, i.e., for all
$k,l\in\mZ$ and $s\in\mZ_N$,
\newcommand{\rmns}{\varrho^{\alpha_{k}\beta_{l}}}
\newcommand{\vcn}{\Psi^{\alpha_k\beta_l}}
\begin{equation*}
0=\rmns(\sum_{m,n\in\mZ}A_{mn}V^mU^n)\vcn_s
=\sum_{m,n\in\mZ}A_{mn}e^{i n\frac{\alpha_k}{N}}e^{2\pi is\frac{M}{N}}e^{i m\frac{\beta_l}{N}}
\vcn_{s+[m]_N}.
\end{equation*}
Then, for all $s\in\mZ_N$, $m\in\{0,1,2,\ldots, N-1\}$,
\begin{equation*}
\sum_{j,n\in\mZ}A_{(m+jN)n}
e^{i n\frac{\alpha_k}{N}}e^{i\beta_l j}=0 \text{ for all } k,l\in\mN.
\end{equation*}
This, by the argument on the number of distinct roots of a finite polynomial, implies that
$A_{mn}=0$, for all $m,n\in\mZ$. 

\subsection{The quantum solid torus}
 \label{quantumsolidtorus}
 The quantum solid torus is an example of a cleft Hopf-Galois extension. 
 
 Denote by $h$ the unitary and central generator of the coordinate algebra of the unit circle $\ffun(S^1)$.
 
Solid torus is the Cartesian product $D\times S^1$ of the unit disc and the unit circle. 
Therefore, one can define a
 coordinate algebra of quantum torus as the tensor product 
 $\ffun(D_p)\ot \ffun(S^1)$ of the coordinate algebra of the quantum unit disc eq. (\ref{disc}) and
 the coordinate algebra of the unit circle. One can 
 introduce a further quantisation parameter
  by making the tensor product noncommutative, i.e., by requesting that the subalgebras 
 $\ffun(D_p)\ot 1$ and $1\ot\ffun(S^1)$ do not
 mutually commute. We identify the generators $x\ot 1$ and $1\ot h$ with $x$ and $h$ respectively.
 The coordinate algebra of the quantum solid torus $\ffun(D_p\times_\theta S^1)$, $0<p<1$, 
 $0\leq \theta<2\pi$,
 is generated as a ${}^\ast$-algebra by $x$ and $h$, subject to the relations    
 \begin{gather}
 hh^\ast=1=h^\ast h,\ \ 
 x^\ast x-pxx^\ast=1-p,\nonumber\\
 hx=e^{i\theta}xh,\ \ 
 hx^\ast=e^{-i\theta}x^\ast h.\label{solidtorus}
 \end{gather}
 The linear basis of $\ffun(D_p\times_\theta S^1)$ consists of the elements of the form
 \begin{equation}
 (1-xx^\ast)^kx^m h^n,\ \ \ 
 k\in\mN_0,\ m,n\in\mZ.\label{soltorbas}
 \end{equation}
 Observe that there exists a surjective ${}^\ast$-algebra morphism
 \begin{equation}
 \pi_\partial:\ffun(D_p\times_\theta S^1)\rightarrow \ffun(T_\theta),\ \pi_\partial(h)=U,\ 
 \pi_\partial(x)=V,
 \end{equation}
 of the quantum solid torus onto the quantum torus (eq. (\ref{torus})), i.e., the 
  quantum torus is the boundary
 of the quantum solid torus.
 
 Enveloping $C^\ast$-algebra $C(D_p\times_\theta S^1)$ 
 can be obtained from $\ffun(D_p\times_\theta S^1)$ 
  using $C^\ast$ representations.  Irreducible representations of 
 $\ffun(D_p\times_\theta S^1)$ include those unitarily isomorphic either to the representation obtained
 by composing one of the irreducible representations of the quantum torus $\ffun(T_\theta)$ with
 the map $\pi_\partial$ as well one of the representations
 $\varrho_{p,\theta}^\alpha:\ffun(D_p\times_\theta S^1)\rightarrow \Bnd(\Hil{H}^\alpha_{p,\theta})$,
 where $\Hil{H}^\alpha_{p,\theta}$ is generated by orthonormal vectors
 $\Psi_n$, $n\in \mN_0$, and
 \begin{gather}
 \varrho_{p,\theta}^\alpha(x)\Psi_n=\sqrt{1-p^{n+1}}\Psi_{n+1},\ \ 
 \varrho_{p,\theta}^\alpha(x^\ast)\Psi_n=\sqrt{1-p^n}\Psi_{n-1} \text{ if } n>0,\nonumber\\
 \varrho_{p,\theta}^\alpha(x^\ast)\Psi_0=0,\ \ 
 \varrho_{p,\theta}^\alpha(h^{\pm 1})\Psi_n=e^{\pm i(\alpha+n\theta)}\Psi_n.\label{soltorrep}
 \end{gather}
 If $\theta$ is irrational then the representation (\ref{soltorrep}) is faithful. If $\theta$ is rational,
then for any sequence $(\alpha_n)_{n\in \mN}$, such that $0\leq\alpha_n<2\pi$ and $\alpha_i\neq\alpha_j$ if 
$i\neq j$, for all $i,j,n\in\mN$, the representation 
\[\bigoplus_{n\in \mN}\varrho_{p,\theta}^{\alpha_n}:
\ffun(D_p\times_\theta S^1)\rightarrow \Bnd(\bigoplus_{n\in\mN}\Hil{H}^{\alpha_n}_{p,\theta})\]
is faithful. 

 Denote by $H=\ffun(U(1))$ the coordinate algebra of $U(1)$ and let $u$ be the unitary generator of 
 $H$. The algebra $\ffun(D_p\times_\theta S^1)$ is clearly a right  $H$-comodule $*$-algebra, 
 with the coaction
 defined on the generators by 
 \begin{equation}
 \coact^H(x)=x\ot 1, \ \ \coact^H(h)=h\ot u.\label{solidtco}
 \end{equation}
 It is easy to see that, $\ffun(D_p\times_\theta S^1)^{\co H}=\ffun(D_p)$ and
 $\ffun(D_p\times_\theta S^1)(\ffun(D_p))^H_{\gamma_T}$ is a cleft $H$-Hopf Galois extension, where,
 for all $n\in \mZ$,
 \begin{gather}
 \gamma_T:H\rightarrow \ffun(D_p\times_\theta S^1),\ \ 
 u^n\mapsto h^n,\nonumber\\
 \gamma_T^{-1}:H\rightarrow \ffun(D_p\times_\theta S^1),\ \ 
 u^n\mapsto h^{-n},\\
 \end{gather}
 are the cleaving map and its convolution inverse, respectively.

\subsection{Gluing of two quantum solid tori}

Let $H=\ffun(U(1))$ be the Hopf algebra generated by a unitary and group-like element $u$.

Let the deformation parameters $p,q\in (0,1)$, $\theta,\theta',\theta''\in \mR$. We define
$P_1=\ffun(D_p\times_\theta S^1)$, $P_2=\ffun(D_q\times_{\theta'} S^1)$ 
(Subsection~\ref{quantumsolidtorus}), $P_{12}=\ffun(T_{\theta''})$ 
(Subsection~\ref{quantumtorus}).  A $\ast$-algebra $P_1$ is generated  by 
the elements  $x$, $h$, which satisfy relations (\ref{solidtorus}), and it is a right 
$H$-comodule algebra with the coaction defined by (\ref{solidtco}).
The corresponding generators of $P_2$, $y$ and $g$, satisfy the relations 
 \begin{gather}
 gg^\ast=1=g^\ast g,\ \ 
 y^\ast y-qyy^\ast=1-q,\nonumber\\
 gy=e^{i\theta'}yg,\ \ 
 gy^\ast=e^{-i\theta'}y^\ast g.\label{ygstcomu}
 \end{gather}
 $P_2$ is a right $H$-comodule $\ast$-algebra with the coaction defined on generators
 by $\coact^H(y)=y\ot 1$, $\coact^H(g)=g\ot u$. Finally, $P_{12}$ is a right $H$-comodule 
 $\ast$-algebra generated by unitary elements $U$ and $V$ satisfying $UV=e^{i\theta''}VU$, 
 with the right $H$-coaction defined by the relations $\coact^H(V)=V\ot 1$, $\coact^H(U)=U\ot u$.
 Note that $B_1=P_1^{\co H}\simeq \ffun(D_p)$ (see Section~\ref{quantumdisc}), 
 is generated as a $\ast$-algebra by $x$, $B_2=P_2^{\co H}\simeq \ffun(D_q)$ is generated by $y$ and
 $B_{12}=P_{12}\simeq\ffun(S^1)$ is generated by $V$.
 Let the algebra surjections $\pi^1_2:B_1\rightarrow B_{12}$, $\pi^2_1:B_2\rightarrow B_{12}$ be 
 defined on generators by
 \begin{equation}
 \pi^1_2(x)=V,\ \ \pi^2_1(y)=V.\label{pijstbas}
 \end{equation}
 We define cleaving maps by
 \begin{equation}
 \gamma_1^{\pm1}(u^n)=h^{\pm n},\ \ 
 \gamma_2^{\pm1}(u^n)=g^{\pm n},\ \ 
 \gamma_{12}^{\pm1}(u^n)=U^{\pm n},\ \ \ 
 \text{for all }n\in \mZ.\label{gammastcl}
 \end{equation}  
 Then $P_1(B_1)^H_{\gamma_1}$, $P_2(B_2)^H_{\gamma_2}$, $P_{12}(B_{12})_{\gamma_{12}}^H$
 are cleft extensions. By Lemma~\ref{priglucllem}, in order to define 
 gluing surjections (\ref{mapschiij}), $\chi^1_2:P_1\rightarrow P_{12}$, 
 $\chi^2_1:P_2\rightarrow P_{12}$, we need to find appropriate convolution invertible
 maps $\Gamma^1_2,\Gamma^2_1:H\rightarrow B_{12}$. 
 To fix the notation, without losing the generality, we shall only consider $\Gamma^2_1$.
  
 For all $n\in\mZ$, $\Gamma^2_1(u^n)$ and $(\Gamma^2_1)^{-1}(u^n)$ are Laurent polynomials
 in $V$ such that $\Gamma^2_1(u^n)(\Gamma^2_1)^{-1}(u^n)=1$. By the standard argument
 about degree counting, this implies that
 \begin{equation}
 (\Gamma^2_1)^{\pm 1}(u^n)=\mu(n)^{\pm 1}V^{\pm \nu(n)},
 \text{ where } \mu:\mZ\rightarrow\mC\setminus\{0\},\ \nu:\mZ\rightarrow\mZ.\label{bgdef}
 \end{equation}
 The map $\chi^2_1$ must be algebraic, hence in paricular, 
 \begin{equation*}
 \chi^2_1(y^mg^ny^kg^l)=\chi^2_1(y^mg^n)\chi^2_1(y^kg^l),\text{  for all } m,n,k,l\in \mZ.
 \end{equation*}
 Substituting (\ref{pijstbas}), (\ref{gammastcl}), (\ref{bgdef}) and
 (\ref{chiijgammaij}) yields
 \begin{equation*}
\chi^2_1(y^mg^ny^kg^l)=\chi^2_1(e^{ink\theta'}y^my^kg^{n+l})
=e^{ink\theta'}\mu(n+l)V^{m+k+\nu(n+l)}U^{n+l},
 \end{equation*}
 and
 \begin{multline*}
 \chi^2_1(y^mg^n)\chi^2_1(y^kg^l)=\mu(n)\mu(l)V^{m+\nu(n)}U^nV^{k+\nu(l)}U^l\\
 =e^{i\theta''n(k+\nu(l))}\mu(n)\mu(l)V^{m+k+\nu(n)+\nu(l)}U^{n+l}.
 \end{multline*}
 It follows that, for all $n,l,k\in \mZ$,
 \begin{gather}
 \nu(n+l)=\nu(n)+\nu(l),\label{reclen1}\\
 e^{ink(\theta'-\theta'')}\mu(n+l)=e^{i\theta''n\nu(l)}\mu(n)\mu(l).\label{reclen2}
 \end{gather}
 Condition (\ref{reclen1}) implies that, for all $n\in \mZ$, $\nu(n)=\beta n$, where $\beta=\nu(1)$.
 Only left hand side of 
 condition (\ref{reclen2}) depends on $k$, therefore it  can be satisfied for all $k\in\mZ$ only if
 $\theta'=\theta''$. We assume this, and then we have the following recursive relation
 \begin{equation}
 \mu(n+l)=e^{i\theta'n\beta l}\mu(n)\mu(l),\text{ for all }n,l\in\mZ,
 \end{equation} 
 which has a family of solutions
 \begin{equation}
 \mu(n)=\alpha^ne^{i\beta\theta'\frac{n^2}{2}},\text{ for all }n\in\mZ,
 \end{equation}
 where $\alpha\in\mC\setminus\{0\}$. It follows that, for all $n\in \mZ$,
 $\Gamma^2_1(u^n)=\alpha^ne^{i\beta\theta'\frac{n^2}{2}}V^{\beta n}$.
 
 Similarly we prove that $\theta$ must equal $\theta''$ and then, for all $n\in \mZ$,
 $\Gamma^1_2(u^n)=(\alpha')^ne^{i\beta'\theta\frac{n^2}{2}}V^{\beta' n}$,
 for some $\alpha'\in\mC\setminus\{0\}$ and $\beta'\in\mZ$. In particular
 $\Gamma^1_2(u^n)=1$, for all $n\in \mZ$, is an admissible gauge transformation
 and, by the Remark after
 the Lemma~\ref{priglucllem}, we can assume just that without losing any generality.
 Accordingly, the most general form of gluing maps for two quantum solid tori can be defined 
 on the basis elements (cf. (\ref{soltorbas})) of respective solid tori  as 
 \begin{gather}
 \chi^1_2((1-xx^\ast)^kx^mh^n)=\delta_{k0}V^mU^n,\nonumber\\
 \chi^2_1((1-yy^\ast)^ky^mg^n)=\delta_{k0}\alpha^ne^{i\beta\theta\frac{n^2}{2}}
 V^{m+\beta n}U^n,\label{chi12def}
 \end{gather}
 for all $m,n\in \mZ$, $k\in\mN_0$.
 Note that $\chi^1_2$ is a $\ast$-algebra map, and $\chi^2_1$ is a $\ast$-algebra map if $|\alpha|=1$.
Observe that the glued algebra $P=\bigoplus_{\chi^i_j}P_i$ is  a $\ast$-algebra in a natural way
 (i.e., with a $\ast$-operation defined by starring the  components of the direct sum) if and only if
the maps $\chi^i_j$ are $\ast$-algebra morphisms. On the other hand, scaling of $g$ by a 
number of modulus
one is an $H$-comodule $\ast$-algebra isomorphism of $P_2$. It follows that,
if  $|\alpha|=1$, the parameter $\alpha$
  can be absorbed, up to an isomorphism, by the redefinition $g\mapsto \alpha^{-1} g$.
Accordingly, in what follows, we shall only consider the case $\alpha=1$. 

Let us denote the generators of the  algebra $P_1^{-}=\ffun(D_p\times_{-\theta}S^1)$
(resp. $P_2^{-}=\ffun(D_q\times_{-\theta}S^1)$,
$P_{12}^{-}=\ffun(T_{-\theta})$) with the same symbols as the generators of $P_1$ (resp.
$P_2$, $P_{12}$). We define, by the action on generators, the $\ast$-algebra isomorphisms
\begin{gather}
\eta_1:P_1\rightarrow P_1^{-},\ \ x\mapsto x,\ h\mapsto h^\ast,\nonumber\\
\eta_2:P_2\rightarrow P_2^{-},\ \ y\mapsto y,\ g\mapsto g^\ast,\nonumber\\
\eta_{12}:P_{12}\rightarrow P_{12}^{-},\ \ V\mapsto V,\ U\mapsto U^\ast.
\end{gather}
Clearly, for all $m,n\in\mZ$, $k\in\mN_0$,
\begin{gather}
\bar{\chi}^1_2=\eta_{12}\circ\chi^1_2\circ\eta_1^{-1}:P_1^{-}\rightarrow P_{12}^{-},\ \ 
(1-xx^\ast)^kx^mh^n\mapsto \delta_{k0}V^mU^n,\nonumber\\
\bar{\chi}^2_1=\eta_{12}\circ\chi^2_1\circ\eta_2^{-1}:P_2^{-}\rightarrow P_{12}^{-},\ \ 
(1-xx^\ast)^kx^mh^n\mapsto \delta_{k0}e^{i\beta\theta\frac{n^2}{2}}V^{m-\beta n}U^n.
\end{gather}
Thus maps $\bar{\chi}^1_2$, $\bar{\chi}^2_1$ have
 the same form as maps (\ref{chi12def}) after substituting $\theta\mapsto -\theta$,
$\beta\mapsto-\beta$. Denote $P^{-}=P^{-}_1\oplus_{\bar{\chi}^i_j}P^{-}_2$. It follows that
the map 
\begin{equation}
\eta=\eta_1\oplus\eta_2:P\rightarrow P^{-}
\end{equation}
is a $\ast$-algebra isomorphism. Consequently, without losing generality, 
in what follows we shall confine ourselves
 to the 
case $\beta\geq 0$.

Using (\ref{chi12def}) and Lemma~\ref{fibprobaslem}, it is easy to see, that the vectors
\begin{gather}
((1-xx^\ast)^kx^mh^n,0),\ \ 
(0,(1-yy^\ast)^ky^mg^n),\ \ 
(x^mh^n,e^{-i\beta\theta\frac{n^2}{2}}y^{m-\beta n}g^n),\label{basislens}
\end{gather}
$m,n\in\mZ$, $k>0$, form a basis of $P$.
\begin{mlem}\label{lensgenslem}
The elements
\begin{equation}
\xi=(1-xx^\ast,0),\ 
z=(x,y),\ a=(e^{\frac{i\beta\theta}{2}}x^\beta h,g),\ 
b=(e^{\frac{i\beta\theta}{2}}h^{-1},y^{\beta}g^{-1})\label{lensgens}
\end{equation}
of $P$ generate $P$ as a $\ast$-algebra. 
\end{mlem}
\begin{proof}
It is enough to show that the basis vectors
(\ref{basislens}) are expressible in terms of elements (\ref{lensgens}). Observe that
\begin{equation}
(0,1-yy^\ast)=1-zz^\ast-\xi.
\end{equation} 
It follows immediately that, for all $k>0$, $m,n\in \mZ$,
\begin{equation}
((1-xx^\ast)^kx^mh^n,0)=e^{\frac{i\beta\theta}{2}n}\xi^kz^mb^{-n},
\end{equation}
and
\begin{equation}
(0,(1-yy^\ast)^ky^mg^n)=(1-zz^\ast-\xi)^kz^ma^n.
\end{equation}
Furthermore, using equation (\ref{xmxn}), for all $m,n\in \mZ$,
\begin{multline*}
(x^mh^n,e^{-i\beta\theta\frac{n^2}{2}}y^{m-\beta n}g^n)\\
=(x^mh^n,e^{-i\beta\theta\frac{n^2}{2}}y^my^{-\beta n}g^n)
-(0,e^{-i\beta\theta\frac{n^2}{2}}y^{m-\beta n}Q^q_{m;-\beta n}(1-yy^\ast)g^n)\\
=z^m(h^n,e^{-i\beta\theta\frac{n^2}{2}}y^{-\beta n}g^n)-
e^{-i\beta\theta\frac{n^2}{2}}z^{m-\beta n}Q^q_{m;\beta n}(1-zz^\ast-\xi)a^n.
\end{multline*}
As $y^{-\beta n}g^n=e^{i\beta\theta\frac{n(n-1)}{2}}(y^{-\beta}g)^n$, it follows that
\begin{equation*}
(h^n,e^{-i\beta\theta\frac{n^2}{2}}y^{-\beta n}g^n)(h,e^{-\frac{i\beta\theta}{2}}y^{-\beta}g)^n=e^{\frac{i\beta\theta}{2}n}b^{-n}.
\end{equation*}
\end{proof}

\begin{mlem}\label{lensrellem}
The generators $\xi$, $z$, $a$, $b$ of $P$ satisfy the following relations.
\begin{subequations}\label{lens}
\begin{gather}\xi^\ast=\xi,\ \ \xi z=pz\xi,\ \ z^\ast z-qzz^\ast=1-q-(p-q)\xi,\label{lensa}\\
(1-zz^\ast-\xi)\xi=0,\label{lensb}\\
\xi a=p^\beta a\xi,\ \ \xi b=b\xi,\ \ 
za=e^{-i\theta}az,\ \ zb=e^{i\theta}bz,\label{lensc}\\
za^\ast-e^{i\theta}a^\ast z=(p^\beta-1)\xi z^{1-\beta}b,\label{lensd}\\
z^\ast b-e^{-i\theta}bz^\ast=(1-q^\beta)z^{\beta-1}(1-zz^\ast-\xi)a^\ast,\label{lense}\\
ab=e^{i\beta\theta}ba,\ \ ab^\ast=e^{-i\beta\theta}b^\ast a,\label{lensf}\\
ba=z^\beta,\label{lensg}\\
a^\ast a=\sum_{m=0}^\beta (-1)^mp^{\beta m-\frac{m(m-1)}{2}}\pchoose{\beta}{m}{p^{-1}}\xi^m,
\label{lensh}\\
aa^\ast=\sum_{m=0}^\beta (-1)^mp^{-\beta m+\frac{m(m+1)}{2}}\pchoose{\beta}{m}{p}\xi^m,\label{lensi}\\ 
b^\ast b=\sum_{m=0}^\beta (-1)^mq^{\beta m-\frac{m(m-1)}{2}}\pchoose{\beta}{m}{q^{-1}}
(1-zz^\ast-\xi)^m,\label{lensj}\\
bb^\ast=\sum_{m=0}^\beta (-1)^mq^{-\beta m+\frac{m(m+1)}{2}}\pchoose{\beta}{m}{q}
(1-zz^\ast-\xi)^m.\label{lensk}
\end{gather}
\end{subequations}
where $\pchoose{n}{m}{p}$ are deformed binomial coefficients defined in (\ref{pchoosedef}) 
\end{mlem}
\begin{proof}
Easy if tedious proof is left to the reader.
\end{proof}

By the discussion in Section~\ref{glusect}, the algebra $P$ is naturally an $H$-comodule
$\ast$-algebra. The coaction $\coact^H:P\rightarrow P\ot H$ is defined on generators by
\begin{equation}\label{lenscoact}
\coact^H(\xi)=\xi\ot 1,\ \ 
\coact^H(z)=z\ot 1,\ \ 
\coact^H(a)=a\ot u,\ \ 
\coact^H(b)=b\ot u^{-1}.
\end{equation}
It is clear (cf. discussion around eq.~(\ref{eqfirstl})) that
$P^{\co H}=B=B_1\oplus_{\pi^i_j}B_2$. It follows that $P^{\co H}$ is generated by 
the elements $\xi,z\in P$.

\subsection{Lens spaces of positive charge}

Let $p,q\in(0,1)$, $\theta\in[0,2\pi)$, $\beta\in \mN_0$, and let
 $\ffun(L_\beta^{p,q,\theta})$ be the quotient of a free $\ast$-algebra generated by the elements
$\xi$, $z$, $a$, $b$, modulo the relations (\ref{lens}). We will call $\ffun(L^{p,q,\theta}_\beta)$ 
a {\em coordinate algebra of functions on a quantum lens space $L^{p,q,\theta}_\beta$
of positive charge} $\beta$.

\newcommand{\msv}{{\mathcal{A}}}
Consider the family  
\begin{equation}
\label{lensbasis}
\{\xi^kz^mb^n\;|\;k>0,\; m,n\in\mZ\},\ \ 
\{(1-zz^\ast-\xi)^kz^ma^n\;|\;k\geq 0,\; m,n\in \mZ\}
\end{equation}
of vectors in $\ffun(L^{p,q,\theta}_\beta)$. We will prove that it is a basis of 
$\ffun(L^{p,q,\theta}_\beta)$. First we need to prove several technical lemmas.
Let
\begin{equation}\label{msvdef}
\msv=(\text{ Span of the family (\ref{lensbasis}) }).
\end{equation}
\begin{mlem}\label{ellem}
The elements $1_\msv$, $\xi$, $z$, $z^\ast$, $a$, $a^\ast$, $b$, $b^\ast$ belong to $\msv$.
\end{mlem}
\begin{proof}
The assertion is obvious in the case of the elements
$1_\msv$, $\xi$, $z$, $z^\ast$, $a$, $a^\ast$. Furthermore,
\begin{multline*}
b=b\left(aa^\ast-\sum_{m=1}^\beta 
(-1)^mp^{-\beta m+\frac{m(m+1)}{2}}\pchoose{\beta}{m}{p}\xi^m\right)\\
=z^\beta a^\ast-\sum_{m=1}^\beta (-1)^mp^{-\beta m+\frac{m(m+1)}{2}}\pchoose{\beta}{m}{p}\xi^mb,
\end{multline*} 
where we used eq.~(\ref{lensi}) in the first equality, and eq.~(\ref{lensg}) in the second.  Therefore
$b\in \msv$.
Similarly, using eq.~(\ref{lensh}) and eq.~(\ref{lensg}), we obtain
\begin{equation*}b^\ast=e^{-i\beta\theta}z^{-\beta}a-
\sum_{m=1}^\beta (-1)^mp^{\beta m-\frac{m(m-1)}{2}}\pchoose{\beta}{m}{p^{-1}}\xi^mb^\ast\in\msv.
\end{equation*}
\end{proof}

\begin{mlem}
The following relations are satisfied in
$\ffun(L^{p,q,\theta}_\beta)$:
\begin{subequations}\label{lensprop}
\begin{gather}(1-zz^\ast-\xi)z=qz(1-zz^\ast-\xi),\label{lenspropa}\\ 
(1-zz^\ast-\xi)a=a(1-zz^\ast-\xi),\label{lenspropb}\\
(1-zz^\ast-\xi)b=q^\beta b(1-zz^\ast-\xi),\label{lenspropc}\\
\xi a^n=e^{i\beta\theta\frac{n(n+1)}{2}}\xi z^{\beta n}b^{-n},\text{ for all }n\in\mZ.
\label{lenspropd}
\end{gather}
\end{subequations}
\end{mlem}
\begin{proof}
Using (\ref{lensa}) yields
\begin{multline*}(1-zz^\ast-\xi)z=z-z(z^\ast z)-pz\xi\\
=z-z(1-q(1-zz^\ast-\xi)-p\xi)-pz\xi=qz(1-zz^\ast-\xi).
\end{multline*}
Furthermore,
\begin{multline*}(1-zz^\ast-\xi)a=a-z(z^\ast a)-p^\beta a\xi
=a-e^{i\theta}z(az^\ast+(1-p^\beta)b^\ast z^{\beta-1}\xi)-p^\beta a\xi\\
=a-azz^\ast-e^{i\beta\theta}(1-p^\beta)z^\beta b^\ast\xi-p^\beta a\xi\\
=a-azz^\ast-(1-p^\beta)abb^\ast\xi-p^\beta a\xi
=a(1-zz^\ast-\xi),
\end{multline*}
where in the forth equality we used eq.~(\ref{lensg}). Similar proof of the equation (\ref{lenspropc}) is
left to the reader as an exercise. 

To prove the property (\ref{lenspropd}), we note that, by (\ref{lenspropc}),
(\ref{lensj}) and (\ref{lensk}),
\begin{equation}
\xi b^nb^{-n}=1,\text{ for all }n\in\mZ.
\end{equation}
Therefore, for all $n\in\mZ$,
\begin{equation*}
e^{i\beta\theta\frac{n(n+1)}{2}}\xi z^{\beta n}b^{-n}
=e^{i\beta\theta\frac{n(n+1)}{2}}\xi (ba)^nb^{-n}
=\xi a^nb^nb^{-n}=\xi a^n.
\end{equation*}
\end{proof}

\begin{mlem}\label{mullcolem}
The vector subspace $\msv\subseteq \ffun(L^{p,q,\theta}_\beta)$ (eq.~(\ref{msvdef})) is closed under multiplication.
\end{mlem}
\begin{proof}
It is enough to consider products of basis vectors (\ref{lensbasis}). Observe that,
by equations (\ref{lenspropa})--(\ref{lenspropc}),
\begin{equation*}
\{\xi^kz^mb^n\;|\;k\in\mN,\; m,n\in\mZ\}\cdot\{(1-zz^\ast-\xi)^lz^sa^t\;|\;l\in\mN,\; s,t\in\mZ\}=\{0\},
\end{equation*}
and
\begin{equation*}
\{(1-zz^\ast-\xi)^lz^sa^t\;|\;l\in\mN,\; s,t\in\mZ\}\cdot\{\xi^kz^mb^n\;|\;k\in\mN,\; m,n\in\mZ\}=\{0\}.
\end{equation*}

Note that
$zz^\ast=1-(1-zz^\ast-\xi)-\xi$, and, by (\ref{lensa}),
\[z^\ast z=1-q(1-zz^\ast-\xi)-p\xi.\]
 It follows, using (\ref{lensa}), (\ref{lensb}) and (\ref{lenspropa}),
that, for all $n,m\in\mZ$,
\begin{subequations}\label{zabs}
\begin{equation}\label{zabsa}
z^nz^m=(1+\mathcal{P}_{n,m}(\xi)+\mathcal{Q}_{n,m}(1-zz^\ast-\xi))z^{n+m},
\end{equation}
where $\mathcal{P}_{n,m}$ and $\mathcal{Q}_{n,m}$ are polynomials such that
$\mathcal{P}_{n,m}(0)=\mathcal{Q}_{n,m}(0)=0$. Similarly by 
relations (\ref{lensh})--(\ref{lensk}), (\ref{lensc}) and (\ref{lenspropb})--(\ref{lenspropc}),
for all $n,m\in \mZ$,
\begin{gather}
a^na^m=(1+\mathcal{P}^{'}_{n,m}(\xi))a^{n+m},\label{zabsb}\\
b^nb^m=(1+\mathcal{Q}^{'}_{n,m}(1-zz^\ast-\xi))b^{n+m},\label{zabsc}
\end{gather}
where polynomials $\mathcal{P}^{'}_{n,m}$,  $\mathcal{Q}^{'}_{n,m}$,
 satisfy
$\mathcal{P}^{'}_{n,m}(0)=\mathcal{Q}^{'}_{n,m}(0)=0$. 
\end{subequations}

It follows that, for all $m,n,s,t\in\mZ$ and $k,l\in\mN$,
\begin{multline}\label{ksizbksizb}
(\xi^k z^mb^n)(\xi^l z^sb^t)
=p^{-lm}e^{-ins\theta}\xi^{k+l}z^mz^sb^{n+t}\\
=p^{-lm}e^{-ins\theta}\xi^{k+l}(1+\mathcal{P}_{m,s}(\xi))z^{m+s}b^{n+t}\\
\in\Span(\{\xi^kz^mb^n\;|\;k\in\mN,m,n\in\mZ\})\subseteq\msv.
\end{multline}
Using eq.~(\ref{lenspropd}) yields, for all
$m,n,s,t\in \mZ$, $k\in\mN$,
\begin{equation*}
(\xi^kz^mb^n)(z^sa^t)=e^{i\beta\theta\frac{t(t+1)}{2}}\xi^kz^mb^nz^sz^{\beta t}b^{-t},
\end{equation*}
hence, by eq.~(\ref{ksizbksizb}), 
\begin{equation}\label{xizbza}(\xi^kz^mb^n)(z^sa^t)\in\Span(\{\xi^kz^mb^n\;|\;k\in\mN,m,n\in\mZ\})\subseteq\msv.
\end{equation}
Analogously, 
\begin{multline}
(z^sa^t)(\xi^kz^mb^n)=p^{-\beta tk} z^s\xi^k a^t z^m b^n
=e^{i\beta\theta\frac{t(t+1)}{2}}p^{-\beta tk} z^s\xi^kz^{\beta t} b^{-t} z^mb^n\\
\in\Span(\{\xi^kz^mb^n\;|\;k\in\mN,m,n\in\mZ\})\subseteq\msv.\label{zaxizb}
\end{multline}
In the remainder of the proof we need the following observation. For all $m,n\in\mZ$,
\begin{equation}\label{azzaspan}
a^nz^m\in e^{imn\theta}z^ma^n+\Span(\{\xi^kz^mb^n\;|\;k\in\mN,m,n\in\mZ\}).
\end{equation}
We use induction on $m,n\in\mZ$. The above formula is obviously true for $m$ or $n$ equal to zero.
By eq.~(\ref{lensc}) and eq.~(\ref{lensd}), it is also true for $m,n=\pm 1$. 
For brevity write  $\msv'=\Span(\{\xi^kz^mb^n\;|\;k\in\mN,m,n\in\mZ\})$.
Let $k=\pm 1$, $kn\geq 0$.
Then, using equations (\ref{ksizbksizb}), (\ref{xizbza}), (\ref{zaxizb}), we obtain
\begin{multline*}
a^{n+k}z^m=a^k(a^nz^m)\in e^{imn\theta}(a^kz^m)a^n+a^k\msv'\subseteq
e^{imn\theta}(a^kz^m)a^n+\msv'\\
\subseteq e^{im(n+k)\theta}z^ma^{n+k}+e^{imn\theta}\msv' a^n+\msv'
\subseteq e^{im(n+k)\theta}z^ma^{n+k}+\msv'. 
\end{multline*}
Similarly, for $k=\pm 1$, $km\geq 0$,
\begin{multline*}
a^nz^{m+k}=(a^nz^m)z^k\in e^{imn\theta}z^m(a^nz^k)+\msv' z^k\subseteq
e^{imn\theta}z^m(e^{ink\theta}z^ka^n+\msv')+\msv' z^k\\
\subseteq e^{i(m+k)n\theta}z^{m+k}a^n+\msv'.
\end{multline*}
Using (\ref{azzaspan}) and then (\ref{zabs}) and (\ref{lensprop}), we obtain, for all $m,n,s,t\in \mZ$,
\begin{multline*}
z^ma^nz^sa^t
\in e^{ins\theta}z^mz^sa^na^t+z^m\msv' a^t\\
\subseteq e^{ins\theta}(1+\mathcal{Q}_{m,s}(1-zz^\ast-\xi)+\mathcal{P}_{m,s}(\xi))z^{m+s}
(1+\mathcal{P}^{'}_{n,t}(\xi))a^{n+t}+\msv'\\
=e^{ins\theta}(1+\mathcal{Q}_{m,s}(1-zz^\ast-\xi))z^{m+s}a^{n+t}
+\mathcal{P}^{''}(\xi)z^{m+s+\beta(n+t)}b^{-(n+t)}+\msv'\subseteq\msv,
\end{multline*} 
where $\mathcal{P}^{''}$ is a polynomial such that $\mathcal{P}^{''}(0)=0$, and in the last equality
we used eq.~(\ref{lenspropd}) and then eq.~(\ref{zabsa}). This shows that,
for all $m,n,s,t\in\mZ$ and $k,l\in\mN_0$,
$((1-zz^\ast-\xi)^kz^ma^n)((1-zz^\ast-\xi)^lz^sa^t)\in\msv$, which ends the proof.
\end{proof}

\begin{mprop}
Vectors (\ref{lensbasis}) form a  basis of $\ffun(L^{p,q,\theta}_\beta)$. The algebras
 $\ffun(L^{p,q,\theta}_\beta)$ and $P=P_1\oplus_{\chi^i_j}P_2
 =\ffun(D_p\times_\theta S^1)\oplus_{\chi^i_j}\ffun(D_q\times_\theta S^1)$
are mutually isomorphic. Here the maps
 $\chi^1_2$ and $\chi^2_1$ are defined in eq.~(\ref{chi12def}) with
 $\alpha=1$ and $\beta\geq 0$.
\end{mprop}
\begin{proof}
Let the algebra maps $\chi_i:\ffun(L^{p,q,\theta}_\beta)\rightarrow P_i$, $i=1,2$, be defined 
on generators by
\begin{gather}
\chi_1(\xi)=1-xx^\ast,\ \ \chi_1(z)=x,\ \ \chi_1(a)=e^{\frac{i\beta\theta}{2}}x^\beta h,\ \  
\chi_1(b)=e^{\frac{i\beta\theta}{2}}h^{-1},\nonumber\\
\chi_2(\xi)=0,\ \ \chi_2(z)=y,\ \ \chi_2(a)=g,\ \ \chi_2(b)=y^\beta g^{-1}.\label{chilensdef}
\end{gather}
By Lemmas~\ref{lensgenslem} and 
\ref{lensrellem} these maps are well defined, and by Lemma~\ref{lensgenslem},
the map $\chi=\chi_1\oplus\chi_2:\ffun(L^{p,q,\theta}_\beta)\rightarrow P$ is surjective.
Let $w\in\ker\chi$.
By Lemmas~\ref{ellem} and \ref{mullcolem}, vectors (\ref{lensbasis}) span $\ffun(L^{p,q,\theta}_\beta)$,
hence
\begin{equation}
w=\sum_{m,n\in\mZ\atop k>0}\mu_{kmn}\xi^kz^mb^n+\sum_{s,t\in\mZ\atop l\geq 0}\nu_{lst}
(1-zz^\ast-\xi)^lz^sa^t,
\end{equation}
for some coefficients $\mu_{kmn}, \nu_{lst}\in\mC$, where $m,n,s,t\in\mZ$, $k>0$, $l\geq 0$. 
By assumption, $\chi_1(w)=0$ and $\chi_2(w)=0$. It follows that
\begin{equation*}
0=\chi_2(w)=\sum_{s,t\in\mZ\atop l\geq 0}\nu_{lst}(1-yy^\ast)^ly^sg^t.
\end{equation*}
Since the elements $(1-yy^\ast)^ly^sg^t$, $l\in\mN_0$, $s,t\in\mZ$,  form a linear basis
of $P_2$, this implies that, for all $l\in\mN_0$, $s,t\in\mZ$, $\nu_{lst}=0$.  But then
\begin{equation*}
0=\chi_1(w)=\sum_{m,n\in\mZ\atop k>0}\mu_{kmn}e^{\frac{i\beta\theta}{2}n}(1-xx^\ast)^kx^mh^n,
\end{equation*}
which implies that, for all $k\in\mN$ and $m,n\in\mZ$,
$\mu_{kmn}e^{\frac{i\beta\theta}{2}n}=0$ and so $\mu_{kmn}=0$.
Hence $w=0$ and therefore $\ker\chi=\{0\}$ and so $\chi:\ffun(L^{p,q,\theta}_\beta)\rightarrow P$
is a $\ast$-algebra isomorphism. 
It follows that we can identify $\ffun(L^{p,q,\theta}_\beta)$ with $P$.
We have also proven that vectors (\ref{lensbasis}) are linearly
independent and hence they form a linear basis of  $\ffun(L^{p,q,\theta}_\beta)$.

Let us define a right  $H$-coaction 
$\coact^H:\ffun(L^{p,q,\theta}_\beta)\rightarrow\ffun(L^{p,q,\theta}_\beta)\ot H$ by (eq.~\ref{lenscoact}),
which makes $\chi:\ffun(L^{p,q,\theta}_\beta)\rightarrow P$ an $H$-comodule isomorphism.
It follows, by the discussion after eq.~(\ref{lenscoact}), that $B=P^{\co H}$ is isomorphic to
the quotient of a free algebra, generated by elements, $\xi$, $z$, by the relations
(\ref{lensa})-(\ref{lensb}). This, in turn, is the coordinate algebra $\ffun(S^2_{pq})$ on the quantum 
2-sphere $S^2_{pq}$ defined, by gluing two quantum discs $D_p$ and $D_q$, in \cite{CalMat:Glu} and \cite{HaMaSz:Pro}.
\end{proof}

\subsection{The inverse of the canonical map on $\ffun(L^{p,q,\theta}_\beta)$}

By Proposition~\ref{glutwogallem}, $P(B)^H$ is an $H$-Hopf Galois extension.
The translation maps on $P_1$ and $P_2$ are given explicitly,  by the formulae,
for all $n\in \mZ$,
\begin{gather}
\trans_1:H\rightarrow P_1\ot_BP_1,\ \ u^n\mapsto h^{-n}\ot_B h^n,\nonumber\\
\trans_2:H\rightarrow P_2\ot_BP_2,\ \ u^n\mapsto g^{-n}\ot_B g^n.
\end{gather}
By eq.~(\ref{caninvlocgalexplfor}), the translation map on $P$ is given explictly as, for all $n\in\mZ$, 
\begin{equation}
\trans:H\rightarrow P\ot_BP,\ \ u^n\mapsto\cm^{-1}_{P\ot_BP}(\trans_1(u^n), \trans_2(u^n)),
\end{equation}
i.e., for all $n\in\mZ$, the element $\tau(u^n)\in P\ot_BP$ is uniquely determined by the property
\begin{equation}
(\chi_1\ot_B\chi_1)(\tau(u^n))=h^{-n}\ot_Bh^n,\ \ 
(\chi_2\ot_B\chi_2)(\tau(u^n))=g^{-n}\ot_Bg^n,\label{taucondlens}
\end{equation}
where maps $\chi_1$, $\chi_2$ were defined in (\ref{chilensdef}).
In order to find $\tau(u^n)$, first note that, for all $n\in\mZ$,
\begin{equation}
(\chi_1\ot_B\chi_1)(b^n\ot_B b^{-n})=h^{-n}\ot_Bh^n,\ \ 
(\chi_2\ot_B\chi_2)(a^{-n}\ot_Ba^n)=g^{-n}\ot_Bg^n. 
\end{equation}
Then, for all $n\in\mZ$,
\begin{multline}
b^n\ot_Bb^{-n}
=a^{-n}a^nb^n\ot_Bb^{-n}+(1-a^{-n}a^n)b^n\ot_Bb^{-n}\\
=a^{-n}\ot_Ba^nb^nb^{-n}+(1-a^{-n}a^n)b^n\ot_Bb^{-n}\\
=a^{-n}\ot_Ba^n+a^{-n}\ot_Ba^n(b^nb^{-n}-1)+(1-a^{-n}a^n)b^n\ot_Bb^{-n}.
\end{multline}
Observe that, for all $n\in\mZ$, $1-a^{-n}a^n\in\ker\chi_2$ and $b^nb^{-n}-1\in\ker\chi_1$, therefore
the elements
\begin{equation}
\begin{split}
\trans(u^n)&=b^n\ot_Bb^{-n}+a^{-n}\ot_Ba^n(1-b^nb^{-n})\\
{}&=a^{-n}\ot_Ba^n+(1-a^{-n}a^n)b^n\ot_Bb^{-n},
\end{split}
\end{equation}
$n\in\mZ$,
satisfy conditions (\ref{taucondlens}), and hence define the translation map on $P$. 

The above method of computation was inspired by the proof of Proposition~1 in \cite{CalMat:Glu}.

\subsection{Representations of $\ffun(L^{p,q,\theta}_\beta)$}
\newcommand{\repr}{\varrho}

To find  representations of  $\ffun(L^{p,q,\theta}_\beta)$
we use the same method as was used in \cite{HaMaSz:LocHopFib}.
Let $\repr:\ffun(L^{p,q,\theta}_\beta)\rightarrow\Bnd(\Hil{H})$ be any representation of 
$\ffun(L^{p,q,\theta}_\beta)$ as a subalgebra of the algebra of bounded operators 
$\Bnd(\Hil{H})$ on a Hilbert space $\Hil{H}$. Note that, by the relations
(\ref{lensprop}), (\ref{lensa}) and (\ref{lensc}), the subspaces $\ker\repr(\xi)$ and
$\ker\repr(1-zz^\ast-\xi)$ are invariant. For any pair of closed subspaces 
$\mathcal{S}\subseteq\mathcal{S}'\subseteq \Hil{H}$, let
$\mathcal{S}^{\perp_{\mathcal{S}'}}$ denote the closure of the orthogonal complement of $\mathcal{S}$ in $\mathcal{S}'$. 
For brevity, we denote $\mathcal{S}^\perp=\mathcal{S}^{\perp_{\Hil{H}}}$.
In this section symbol `$\oplus$'  denotes the orthogonal direct sum of Hilbert spaces. 
Hilbert space $\Hil{H}$ can be decomposed into a direct sum
\begin{equation*}
\begin{split}
\Hil{H}=\ker\repr(\xi)\oplus(\ker\repr(\xi))^\perp
&=(\ker\repr(\xi)\cap\ker\repr(1-zz^\ast-\xi))\\
{}&\ \ \ \oplus
(\ker\repr(\xi)\cap\ker\repr(1-zz^\ast-\xi))^{\perp_{\ker\repr(\xi)}}\oplus(\ker\repr(\xi))^\perp.
\end{split}
\end{equation*}
Suppose that  $\Psi\in(\ker\repr(\xi))^\perp$ is such that
$\repr(1-zz^\ast-\xi)\Psi\neq 0$. Since $(\ker\repr(\xi))^\perp$ is an invariant subspace, we obtain,
by the relation (\ref{lensb}),
\begin{equation*}
0=\repr(\xi(1-zz^\ast-\xi))\Psi=\repr(\xi)\repr(1-zz^\ast-\xi)\Psi\neq 0,
\end{equation*}
which is a contradiction. It follows that $(\ker\repr(\xi))^\perp\subseteq\ker\repr(1-zz^\ast-\xi)$.
For brevity, let us denote 
$\Hil{H}_0=\ker\repr(\xi)\cap\ker\repr(1-zz^\ast-\xi)$,
$\Hil{H}'=(\ker\repr(\xi)\cap\ker\repr(1-zz^\ast-\xi))^{\perp_{\ker\repr(\xi)}}$,
$\Hil{H}''=(\ker\repr(\xi))^\perp$. It follows that 
$\Hil{H}=\Hil{H}_0\oplus\Hil{H}'\oplus\Hil{H}''$, and we have an orthogonal
 direct sum decomposition of
the representation $\repr$ into subrepresentations
\begin{equation*}
\repr_0:\ffun(L^{p,q,\theta}_\beta)\rightarrow \Bnd(\Hil{H}_0),\ \ 
\repr':\ffun(L^{p,q,\theta}_\beta)\rightarrow \Bnd(\Hil{H}'),\ \ 
\repr'':\ffun(L^{p,q,\theta}_\beta)\rightarrow \Bnd(\Hil{H}'').
\end{equation*}
In the  representation $\repr_0$, the relations (\ref{lens}) are reduced to
\begin{subequations}
\label{lens0}
\begin{gather}
\repr_0(\xi)=0,\ \ \repr_0(b)=\repr_0(z)^\beta\repr_0(a^\ast),\\
\repr_0(z^\ast)\repr_0(z)=1=\repr_0(z)\repr_0(z^\ast),\ \ 
\repr_0(a^\ast)\repr_0(a)=1=\repr_0(a)\repr_0(a^\ast),\\
\repr_0(a)\repr_0(z^{\pm 1})=e^{\pm i\theta}\repr_0(z^{\pm 1})\repr_0(a).
\end{gather}
\end{subequations}
Similarly, in the representation $\repr'$, the relations (\ref{lens}) assume the form
\begin{subequations}\label{lensp}
\begin{gather}
\repr'(\xi)=0,\ \ \repr'(b)=\repr'(z)^\beta\repr'(a^\ast),\\
\repr'(z^\ast)\repr'(z)-q\repr'(z)\repr'(z^\ast)=1-q,\ \ 
\repr'(a^\ast)\repr'(a)=1=\repr'(a)\repr'(a^\ast),\\
\repr'(a)\repr'(z^{\pm 1})=e^{\pm i\theta}\repr'(z^{\pm 1})\repr'(a).
\end{gather}
\end{subequations}
Finally, the representation $\repr''$ reduces the relations (\ref{lens}) to the form
\begin{subequations}
\label{lenspp}
\begin{gather}
\repr''(\xi)=1-\repr''(z)\repr''(z^\ast),\ \ \repr''(a)=\repr''(b^\ast)\repr''(z)^\beta,\\
\repr''(z^\ast)\repr''(z)-p\repr''(z)\repr''(z^\ast)=1-p,\ \ 
\repr''(b^\ast)\repr''(b)=1=\repr''(b)\repr''(b^\ast),\\
\repr''(z)\repr''(b^{\pm 1})=e^{\pm i\theta}\repr''(b^{\pm 1})\repr''(z).
\end{gather}
\end{subequations}
>From the representation theory of the quantum solid torus 
(Section~\ref{quantumsolidtorus}), it follows  that irreducible representations of
$\ffun(L^{p,q,\theta}_\beta)$ include the ones
 unitarily equivalent to one of the following representations.

For all $0\leq\mu<2\pi$, there exists a representation 
$\repr^{'}_\mu:\ffun(L^{p,q,\theta}_\beta)\rightarrow\Bnd(\Hil{H}^{'}_\mu)$, where 
$\Hil{H}^{'}_\mu$ has an orthonormal Hilbert basis $\Psi^{'}_n$, $n\in\mN_0$, such that,
for all $n\in\mN_0$,
\begin{gather}\repr^{'}_\mu(z)\Psi^{'}_n=\sqrt{1-q^{n+1}}\Psi^{'}_{n+1},\ \ 
\repr^{'}_\mu(z^\ast)\Psi^{'}_n=\sqrt{1-q^{n}}\Psi^{'}_{n-1}
\text{ if }n>0,\ \ \repr^{'}_\mu(z^\ast)\Psi^{'}_0=0,\nonumber\\
\repr^{'}_\mu(a^{\pm 1})\Psi^{'}_n=e^{\pm i(\mu+n\theta)}\Psi^{'}_n,\ \ 
\repr^{'}_\mu(\xi)\Psi^{'}_n=0,\nonumber\\ 
\repr^{'}_\mu(b)\Psi^{'}_n=e^{-i(\mu+n\theta)}\sqrt{1-q^{n+1}}\ldots\sqrt{1-q^{n+\beta}}
\Psi^{'}_{n+\beta},\nonumber\\
\repr^{'}_\mu(b^\ast)\Psi^{'}_n\left\{\begin{array}{ll}
0 &\text{if }n<\beta,\\
e^{i(\mu+(n-\beta)\theta)}\sqrt{1-q^{n}}\ldots\sqrt{1-q^{n-\beta+1}}
\Psi^{'}_{n-\beta} & \text{otherwise}.
\end{array}
\right.
\end{gather}

Similarly, for all $0\leq\mu<2\pi$, there exists a representation 
$\repr^{''}_\mu:\ffun(L^{p,q,\theta}_\beta)\rightarrow\Bnd(\Hil{H}^{''}_\mu)$, where 
$\Hil{H}^{''}_\mu$ has an orthonormal Hilbert basis $\Psi^{'}_n$, $n\in\mN_0$, such that,
for all $n\in\mN_0$,
\begin{gather}
\repr^{''}_\mu(z)\Psi^{''}_n=\sqrt{1-p^{n+1}}\Psi^{''}_{n+1},\ \ 
\repr^{''}_\mu(z^\ast)\Psi^{''}_n=\sqrt{1-p^{n}}\Psi^{''}_{n-1}
\text{ if }n>0,\ \ \repr^{''}_\mu(z^\ast)\Psi^{''}_0=0,\nonumber\\
\repr^{''}_\mu(b^{\pm 1})\Psi^{''}_n=e^{\pm i(\mu-n\theta)}\Psi^{''}_n,\ \ 
\repr^{''}_\mu(\xi)\Psi^{''}_n=p^n\Psi^{''}_n,\nonumber\\
\repr^{''}_\mu(a)\Psi^{''}_n=e^{-i(\mu-(n+\beta)\theta)}\sqrt{1-p^{n+1}}\ldots\sqrt{1-p^{n+\beta}}
\Psi^{''}_{n+\beta},\nonumber\\
\repr^{''}_\mu(a^\ast)\Psi^{''}_n\left\{
\begin{array}{ll}
0 &\text{if }n<\beta,\\
e^{i(\mu-n\theta)}\sqrt{1-p^{n}}\ldots\sqrt{1-p^{n-\beta+1}}
\Psi^{''}_{n-\beta} & \text{otherwise}.
\end{array}
\right.
\end{gather}

Finally, depending on whether deformation parameter $\theta$ is a rational or irrational
multiple of $2\pi$, we have one of the following families of irreducible representations.

Suppose that $\theta=2\pi\frac{M}{N}$, where $M,N\in\mZ$, $N> 0$,
and $M$ and $N$ are relatively prime. Then, for all $0\leq\mu,\nu<2\pi$, there exists a representation 
$\repr^{\mu\nu}_0:\ffun(L^{p,q,\theta}_\beta)\rightarrow\Bnd(\Hil{H}^{\mu\nu}_0)$, where 
$\Hil{H}^{\mu\nu}_0$ has an orthonormal Hilbert basis $\Psi_n$, $n\in\mZ_N$, and,
for all $n\in\mZ_N$,
\begin{gather}
\repr^{\mu\nu}_0(a^{\pm 1})\Psi_ne^{\pm i\frac{\mu}{N}\pm in\theta}\Psi_n,\ \ 
\repr^{\mu\nu}_0(z^{\pm 1})\Psi_n=e^{\pm i\frac{\nu}{N}}\Psi_{n\pm 1},\ \ 
\repr^{\mu\nu}_0(\xi)\Psi_n=0,\nonumber\\
\repr^{\mu\nu}_0(b)\Psi_ne^{i\frac{\nu\beta-\mu}{N}-in\theta}\Psi_{n+\beta},\ \ 
\repr^{\mu\nu}_0(b^\ast)\Psi_ne^{i\frac{\mu-\nu\beta}{N}+i(n-\beta)\theta}\Psi_{n-\beta}.
\end{gather}

If $\theta$ is an irrational multiple of $2\pi$, then, for any $0\leq \mu<2\pi$,
we have a representation
$\repr^\mu_0:\ffun(L^{p,q,\theta}_\beta)\rightarrow\Bnd(\Hil{H}^\mu_0)$. The Hilbert space
$\Hil{H}^\mu_0$ has an orthonormal basis
$\Psi_n$, $n\in \mZ$. For all $n\in\mZ$,
\begin{gather}
\repr^\mu_0(a^{\pm 1})\Psi_n=e^{\pm i\mu\pm in\theta}\Psi_n,\ \ 
\repr^\mu_0(z^{\pm 1})\Psi_n=\Psi_{n\pm 1},\ \ 
\repr^\mu_0(\xi)\Psi_n=0,\nonumber\\
\repr^\mu_0(b)\Psi_n=e^{-i\mu-in\theta}\Psi_{n+\beta},\ \ 
\repr^\mu_0(b^\ast)\Psi_n=e^{i\mu+i(n-\beta)\theta}\Psi_{n-\beta}.
\end{gather}

If $\theta$ is irrational, then for any $0\leq\mu,\nu<2\pi$, the representation
$\repr^{'}_\mu\oplus\repr^{''}_\nu:\ffun(L^{p,q,\theta}_\beta)\rightarrow\Bnd(\Hil{H}^{'}_\mu\oplus\Hil{H}^{''}_\nu)$
is faithful.

\subsection{Final remarks}

We conclude this section with a number of remarks about the structure of quantum lens spaces.
We also comment on the $K$-theory of quantum lens spaces.
The detailed developments of the topics discussed here, can be considered as directions for future work.

The relations (\ref{lens}) defining $\ffun(L^{p,q,\theta}_\beta)$ assume a particularly simple form
in the case $\beta=1$.  Namely, if $\beta=1$, then $z=ba$, $\xi=1-aa^\ast$, and the relations
(\ref{lens}) are reduced to 
\begin{gather}
\label{heegard}
a^\ast a-paa^\ast=1-p,\ \ \ b^\ast b-qbb^\ast=1-q,\ \ \ 
ab=e^{i\theta}ba,\ \ \ ab^\ast=e^{-i\theta}b^\ast a,\nonumber\\
(1-aa^\ast)(1-bb^\ast)=0.
\end{gather}
It follows that the quantum lens space of charge $1$, $L^{p,q,\theta}_1$, can be identified with the 
Heegaard quantum sphere $S^3_{p,q,\theta}$ considered in \cite{BaHaMaSz:He}.

For brevity, let us denote $A=\ffun(S^3_{p,q,\theta})$. Let us define a $\mZ$-grading on $A$
with
\begin{equation}
\mdeg(a)=1,\ \ \ \mdeg(a^\ast)=-1,\ \ \ \mdeg(b)=-1,\ \ \ \mdeg(b^\ast)=1,
\end{equation}
and, for all $n\in\mZ$, let $A_n=\{w\in A\;|\;\mdeg(w)=n\}$.  For any $\beta\in\mN$, 
define the subalgebra 
$A(\beta)=\bigoplus_{n\in\mZ}A_{\beta n}$. It can be shown that the $\ast$-algebra map
$f_\beta:\ffun(L^{p,q,\beta\theta}_\beta)\rightarrow A(\beta)$, given on generators by
\begin{equation}
f_\beta(\xi)=1-aa^\ast,\ \ \ 
f_\beta(z)=ba,\ \ \ 
f_\beta(a)=a^\beta,\ \ \ 
f_\beta(b)=e^{i\theta\frac{\beta(\beta-1)}{2}}b^\beta,
\end{equation}
is a well defined $\ast$-algebra isomorphism.

In \cite{BaHaMaSz:He}, 
it was  demonstrated that $C(S^3_{p,q,\theta})$ is isomorphic as a $C^\ast$-algebra 
to a fibre product of two $C^\ast$-algebras isomorphic to quantum solid tori,
and then the Mayer-Vietoris sequence was used to compute the $K$-theory of $C(S^3_{p,q,\theta})$ 
using the  $K$-theory of quantum solid tori. 
We expect that this method can be adapted to compute the $K$-theory of 
$C(L^{p,q,\theta}_\beta)$. This is a direction for future work.

\hspace{12pt}

{\noindent \Large\bf Acknowledgements}

I would like to thank prof. Tomasz Brzezi\'nski for many fruitful discussions.
I would also like to thank prof. Reiner Matthes for pointing out error in the first
version of my paper, concerning
 statements about representations of 
quantum torus and quantum solid torus.
This research was supported by the EPSRC grant GR/S0107801.

\bibliographystyle{plain}

\end{document}